%
%
\input amstex
\documentstyle{amsppt}

\def\a{\alpha}

\def\Ann{\text{\rom{Ann}}}
\def\atv{\wt{a}_{\vv}}
\def\av{a_{\vv}}

\def\b{\beta}
\def\bh{\wh{\b}}

\def\brel{\buildrel}

\def\Ba{B_{\a}}

\def\Bt{\wt{B}}
\def\Bta{\wt{B}_{\a}}
\def\Btv{\wt{B}_{\vv}}
\def\bu{\bullet}

\def\CC{{\bold C}}
\def\cA{{\Cal A}}

\def\cD{{\Cal D}}
\def\cE{{\Cal E}}
\def\cEa{\cE^{\a}}
\def\cF{{\Cal F}}
\def\cG{{\Cal G}}
\def\cGa{\cG^{\a}}
\def\cH{{\Cal H}}

\def\cl{\colon}
\def\cL{{\Cal L}}
\def\cM{{\Cal M}}
\def\cMt{\wt{\Cal M}}
\def\cN{{\Cal N}}
\def\cO{{\Cal O}}
\def\cod{\text{\rom{cod}}}

\def\cQ{{\Cal Q}}
\def\cR{{\Cal R}}
\def\cS{{\Cal S}}
\def\cSa{{\Cal S}^{\a}}
\def\cT{{\Cal T}}
\def\cTa{{\Cal T}^{\a}}
\def\cU{{\Cal U}}
\def\cV{{\Cal V}}
\def\cW{{\Cal W}}

\def\cZ{{\Cal Z}}
\def\d{\delta}
\def\D{\Delta}

\def\deg{\text{\rom{deg}}}
\def\Def{\text{\rom{Def}}}
\def\det{\text{\rom{det}}}
\def\Det{\text{\rom{Det}}}
\def\e{\epsilon}

\def\End{\text{\rom{End}}}
\def\eba{{\bar\e}_{\a}}
\def\epa{\e^{+}_{\a}}
\def\Ext{\text{\rom{Ext}}}
\def\es{\emptyset}
\def\Ev{E_{\vv}}
\def\FF{{\Bbb F}}
\def\g{\gamma}
\def\G{\Gamma}

\def\Hom{\text{\rom{Hom}}}
\def\hra{\hookrightarrow}

\def\ia{i_{\a}}
\def\Ia{\wh{I}_{\a}}
\def\id{\text{\rom{id}}}

\def\ih{{\hat\iota}}
\def\im{\text{\rom{im}}}
\def\Im{\text{\rom{Im}}}
\def\isa{i_{\a}}
\def\isyh{i_{\hat y}}
\def\Jh{\wh{J}}
\def\k{\kappa}
\def\ker{\text{\rom{ker}}}
\def\Ker{\text{\rom{Ker}}}
\def\l{\lambda}

\def\L{\Lambda}
\def\lra{\longrightarrow}

\def\Mf{\cM^{fl}}
\def\Mhv{\wh{\cM}_{\vv}}
\def\Mt{\wt{\cM}}
\def\Mtv{\wt{\cM}_{\vv}}
\def\mua{\mu_{\a}}

\def\muma{\mu^{-}_{\a}}

\def\Mv{\cM_{\vv}}
\def\Mw{\cM_{\ww}}
\def\n{\noindent}
\def\nua{\nu_{\a}}

\def\om{\omega}
\def\omt{\wt{\omega}}
\def\omtv{\wt{\omega}_{\vv}}
\def\op{\oplus}
\def\ot{\otimes}
\def\ov{\overline}
\def\O{\Omega}
\def\Oh{\widehat{\Omega}}
\def\Ot{\widetilde{\Omega}}
\def\Ov{\Omega_{\vv}}
\def\Pa{P_{\a}}
\def\Pab{\ov{P}_{\a}}
\def\Pam{P^{-}_{\a}}
\def\Pap{P^{+}_{\a}}

\def\ph{\wh{p}}
\def\pih{\wh{\pi}}
\def\pit{\wt{\pi}}
\def\pitv{\wt{\pi}_{\vv}}
\def\phih{\wh{\phi}}

\def\PGL{\text{\rom{PGL}}}
\def\PP{{\Bbb P}}

\def\psit{\wt{\psi}}

\def\QQ{{\Bbb Q}}

\def\rk{\text{\rom{rk}}}

\def\Ra{R_{\a}}

\def\Rat{P^{-}_{\a}}

\def\s{\sigma}

\def\Sia{\Si_{\a}}
\def\Si{\Sigma}
\def\Sih{\wh{\Sigma}}
\def\Sib{\ov{\G}}

\def\Sit{\wt{\Sigma}}
\def\Siv{\Sigma_{\vv}}
\def\Sta{\wt{\Si}_{\a}}
\def\Stab{\text{\rom{Stab}}}
\def\sm{\setminus}
\def\ss{\subset}

\def\t{\theta}
\def\T{\Theta}
\def\Th{\wh{\Theta}}
\def\Ta{{\Theta}_{\a}}
\def\Tt{\T_{\tau}}

\def\Thb{\wh{\Theta}_{\bh}}

\def\tm{\times}

\def\tU{\wt{U}}

\def\tX{\widetilde{X}}

\def\Tr{\text{\rom{Tr}}}

\def\ul{\underline}

\def\vf{\varphi}
\def\vv{{\bold v}}

\def\Va{V_{\a}}

\def\Vta{\wt{V}_{\a}}
\def\Wa{W_{\a}}

\def\wh{\widehat}
\def\wt{\widetilde}
\def\ww{{\bold w}}
\def\xh{\wh{x}}

\def\Y{\Upsilon}
\def\Ya{Y_{\a}}

\def\Yab{\ov{Y}_{\a}}
\def\Yap{Y^{+}_{\a}}
\def\yh{\wh{y}}

\def\z{\zeta}
\def\Zta{\wt{Z}_{\a}}
\def\ZZ{{\Bbb Z}}
\topmatter
\title
A new six dimensional irreducible symplectic variety. 
\endtitle
\rightheadtext{A new six dimensional irreducible symplectic variety.}
\author 
Kieran G.~O'Grady
\endauthor
\address
Universit\`a La Sapienza, Dip.to di Matematica G.~Castelnuovo,  
Piazz.le A.~Moro 5, 00185 Roma 
\endaddress
\email
ogrady\@ mat.uniroma1.it
\endemail
\date 
October 18 2000
\enddate
\dedicatory
Dedicato a Riccardino
\enddedicatory
\thanks
Supported by Cofinanziamento MURST 1999-2001
\endthanks
\endtopmatter
\document
\head
1. Introduction.
\endhead
There are three types of ``building blocks'' in the Bogomolov
decomposition~\cite{B, Th.2} of compact K\"ahlerian manifolds with  torsion
$c_1$, namely complex tori, Calabi-Yau varieties, and irreducible symplectic
manifolds.
 We are interested in the last
type, i.e.~simply-connected compact  
K\"ahlerian manifolds carrying a holomorphic symplectic form which spans
$H^{2,0}$. (The holonomy of a Ricci-flat K\"ahler metric is equal to
$Sp(r)$, hence these manifolds are hyperk\"ahler~\cite{B}.) The stock of
available irreducible symplectic manifolds appears to be quite scarce,
expecially if we think of the many examples of CalabiYau's. Every known
irreducible symplectic manifold is a deformation of one of the following
varieties: the Hilbert scheme parametrizing zero-dimensional subschemes of a
$K3$ of fixed length~\cite{B},  the generalized Kummer variety parametrizing
zero-dimensional subschemes of a complex torus of fixed length and whose
associated $0$-cycle sums up to $0$~\cite{B}, the ($10$-dimensional)
desingularization of the moduli space of rank-two semistable torsion-free
sheaves on a $K3$ with $c_1=0$, $c_2=4$ constructed by the
author~\cite{O1}.  Briefly: all known examples are deformations of an
irreducible factor in the Bogomolov decomposition of a
moduli space of semistable sheaves on a surface with trivial canonical bundle
or, as in the last case, of a symplectic desingularization of such a moduli
space. This paper provides a new example in dimension $6$. To put our result
in perspective we recall some results on moduli spaces of sheaves on a
surface with trivial canonical bundle. Let $X$ be such a surface and $D$ an
ample divisor on it: given a vector $\ww\in H^{*}(X;\ZZ)$, we let
$\cM_{\ww}(X,D)$ be the moduli space of $D$-semistable torsion-free
sheaves $F$ on $X$ with Mukai vector  
$$v(F):=ch(F)\sqrt{Td(X)}=\ww.\tag{1.1}$$  
Mukai~\cite{Muk2} proved that the open subset  $\cM_{\ww}^{s}(X,D)$
parametrizing stable sheaves, if non empty, is smooth symplectic (i.e.~it has
a regular symplectic from) of pure dimension equal to 
$$\dim\cM^{s}_{\ww}(X,D)=2+(\ww,\ww),\tag{1.2}$$
where $(\cdot,\cdot)$ is the quadratic form on $H^{*}(X;\ZZ)$ defined by
$$(\ww,\ww):=\int\limits_{X}(\ww_1\wedge\ww_1-2\ww_0\wedge\ww_2).$$
(Here $\ww_i$ is the component of $\ww$ belonging to $H^{2i}(X;\ZZ)$.) 
Yoshioka~\cite{Y1,Y2}, extending previous partial results~\cite{HG,O2},
showed that under a technical condition on $(D,\ww)$   ensuring that all
semistable sheaves are stable, the moduli space $\cM_{\ww}(X,D)$ is either a
deformation of a Hilbert scheme of points on a $K3$, or else its Bogomolov
factors are abelian surfaces and a deformation of a generalized Kummer
variety. In view of this result we search for a moduli  space containing
points parametrizing strictly semistable sheaves (i.e.~non stable), and
singular at these points, admitting a symplectic
desingularization, in the hope that one of the Bogomolov factors of the
desingularization is a new irreducible symplectic variety. This is what was
done to produce the new $10$-dimensional example mentioned above, the
moduli space being that of certain sheaves on $K3$. In this paper we will
carry out successfully this program with a moduli space of sheaves on an
abelian surface, describe as follows. Let $C$ be a smooth
irreducible projective curve of genus two and $J:=Pic^0(C)$. We set
$\vv:=2-2\eta_J$, where $\eta_J\in H^4(J;\ZZ)$ is the orientation class of
$J$. Let $\Mv$ be the moduli  space $\cM_{\vv}(J,\T)$, where  $\T$ is a
Theta divisor. Many of the results proved
in~\cite{O1} for the moduli space $\cM_4$ of torsion-free semistable
rank-two sheaves on a $K3$ with $c_1=0$, $c_2=4$, remain valid for $\Mv$,
provided one makes the following technical assumption (similar
to~\cite{O1,(0.2)}):  
$$\text{There is no divisor $A$ on $J$ such that
$A\cdot\T=0$ and  $A\cdot A=(-2)$.}\tag{1.3}$$   
One such result says
that the singular locus of $\Mv$ coincides with the set of
S-equivalence classes of  strictly semistable sheaves, i.e.~equivalent to
$I_{p_1}\ot\xi_1\op I_{p_2}\ot\xi_2$, where $p_i\in J$ and $\xi_i\in\Jh$
($\Jh:=Pic^0(J)$).
 Most importantly,  the procedure of~\cite{O1} carries over to 
give a symplectic desingularization $\pitv\cl\Mtv\to\Mv$; we let $\omtv$ be the
symplectic form on $\Mtv$.    The variety $\Mtv$ is of pure dimension $10$
(see~(2.1.4)). It is not symplectically irreducible: consider the map     
$$\matrix
\Mv & \brel \av\over\lra & J\tm \Jh \\
[F] & \mapsto & (\sum c_2(F),[\det F]),\\
\endmatrix$$
where $\sum c_2(F)$ (the Albanese map) is the sum of the points (with
multiplicities) of any representative of $c_2(F)\in CH_0(J)$.  Set
$\atv:=\av\circ\pitv$. As is easily checked $\atv$ is surjective, hence
$\Mtv$ is not symplectically irreducible.  To ``factor out $J\tm\Jh$'' we set 
$$\cMt:=\atv^{-1}(0,\wh{0}),\qquad \omt:=\omtv|_{\cMt}.$$
The main result of this paper is the following.

\proclaim{(1.4)Theorem}
Keep assumptions as above: $\Mt$ is a $6$-dimensional  irreducible
symplectic variety, i.e.~one-connected and with $H^{2,0}(\cMt)$ spanned
by the symplectic form $\omt$. Furthermore $b_2(\cMt)=8$.     
\endproclaim

Since the known $6$-dimensional irreducible symplectic varieties  have 
$b_2=7$
or $b_2=23$~\cite{B}, it follows that $\cMt$ cannot be deformed into one of
the known symplectic varieties, not even up to birational equivalence.

The paper is organized as follows. In \S 2  we construct $\Mtv$: this is
straightforward, all the work was done in~\cite{O1}. An important feature of
this desingularization is that there is a moduli-theoretic interpretation of   the
complement of $\pitv^{-1}\{[I_p(\xi)\op I_p(\xi)]\}$ as the set of
$\wt{S}$-equivalence classes of simple semistable sheaves with the given
Chern character, where $\wt{S}$-equivalence is a relation finer than
$S$-equivalence and (slightly) coarser than isomorphism. Finally we check that
$\Mt$ is $6$-dimensional and  symplectic. In \S 3 we notice that the
generalized Lefschetz Hyperplane Section Theorem, together with methods
of~\cite{L1,O1}, gives the following topological result. Let
$\cM:=\av^{-1}(0,\wh{0})$, and $\pit\cl\cMt\to\cM$ be the restriction of
$\pitv$.  Let $\a\in J$, and $\Ta$ be the translate by $\a$ of a symmetric
theta-divisor $\T$. The inclusion   
$$\pit^{-1}\{[F]\in\cM|\ \text{$F|_{\Ta}$ is not locally-free semistable}\}
\hra\cMt\tag{1.5}$$
induces isomorphisms on $\pi_0$, $\pi_1$, and a surjection on $H_2$. 
This is useful because while  
$\cMt$ is mysterious, the left-hand side can be described in terms of other
well-known moduli spaces. In fact if $\a$ is generic, the
left-hand side of~(1.5) decomposes into three irreducible
components $\Sta\cup\Bta\cup\Vta$, which have the following 
description. The component $\Sta$, which  is mapped by
$\pit$ to the locus parametrizing strictly semistable sheaves
singular at some point of $\Ta$, is a $\PP^1$-bundle over the quotient of $C\tm \Jh$
by the equivalence relation~(4.2.4). 
The component $\Bta$, which is mapped by $\pit$ to
the closure of the locus parametrizing stable sheaves which are singular at some
point of $\Ta$, is a $\PP^1$-bundle over $\ov{C}\tm K^{[2]}\Jh$, where $\ov{C}$ is
obtained from $C$ by identifying the points $q_1,q_2$ given by~(4.2.2), and
$K^{[2]}\Jh$ is the (smooth) Kummer surface of $\Jh$. Finally
$\Vta$, which is mapped by $\pit$ to the closure of the locus parametrizing  sheaves
whose restriction to
$\Ta$ is a non-semistable locally-free sheaf, is obtained from a $\PP^1$-bundle
over $\Jh$ by flopping sixteen $(-1,-1)$-curves.  These results are contained
in Sections (4) through (6). These sections also contain descriptions of the
double and triple intersections of the sets $\Sta$, $\Bta$, $\Vta$, and most
of the topological results which are needed in the proof of
Theorem~(1.4). 
 In \S 7 we give the proof of Theorem~(1.4). We show that the left-hand
side of Inclusion~(1.5) is simply connected and that its second Betti number
is at most equal to $8$: this implies that $\Mt$ is simply connected and
that $b_2(\Mt)=8$. Then we produce a Hodge substructure of $H^2(\Mt)$
of dimension $8$, with $h^{2,0}=1$; this finishes the proof of~(1.4).
This shows that Inclusion~(1.5) induces an isomorphism on
$H_2$, while a priori the Lefschetz Hyperplane Theorem only gives that it
is a surjection; we comment on this in~(3.8).  

The reader may understand the logical structure of the proof   by going
through the statements whose numbering is decorated by a superscript $*$.
These are the results which are invoked in Section~(7) to prove
Theorem~(1.4); most of them are contained in the last subsection of Sections
(2) through (6) (we list the exceptions in the introduction to each
section).        
\subhead
Notation 
\endsubhead 
We choose once and for all a Weierstrass point $p_0\in C$. 
Let $u$ be the Abel-Jacobi embedding 
$$\matrix
C & \brel u\over\hra & J \\
p & \mapsto & [p-p_0].\\
\endmatrix
\tag{1.6}$$
Then $\T:=u(C)$ is a symmetric theta-divisor. For $\a\in J$ let
$\T_{\a}:=\T+\a$. We have the isomorphism
$$\matrix
J & \brel\d\over\lra & \Jh \\
x & \mapsto & [\T_x-\T].\\
\endmatrix\tag{1.7}$$
We set $\wh{x}=\d(x)$; in general points of $\Jh$ are
denoted by $\xh,\yh,\ph$, etc. Let $\Th:=\d(\T)$; for $\bh\in\Jh$ let
$\Th_{\bh}:=\Th+\bh$. We will identify both $\Ta$ and
$\Thb$ with $C$ via the maps
$$C\ni p\mapsto
\cases
\ia(p):=\a+u(p-p_0)\in\Ta & \\
i_{\bh}(p):=\bh+\d\circ u(p-p_0)\in\Thb.
\endcases
\tag{1.8}$$ 
Set $j_{\a}:=\ia^{-1}$, $j_{\bh}:=i_{\bh}^{-1}$. 

We let 
$$\phi\cl J\tm\Jh\to J,\qquad \phih\cl J\tm\Jh\to \Jh\tag{1.9}$$
be the projections. Let $\cL$ be the  normalized Poincar\'e line-bundle on
$J\tm\Jh$, i.e.~the tautological line-bundle such that
$\cL|_{0\tm\Jh}\cong\cO_{\Jh}$.   If $Z\ss J$, $\wh{Z}\in\Jh$ are zero-dimensional
subschemes we set   
$$\cL_Z:=\phih_{*}(\cL\ot\phi^{*}\cO_Z)\quad 
\cL_{\wh{Z}}:=\phi_{*}(\cL\ot\phih^{*}\cO_{\wh{Z}}).$$
In particular if $Z=x$ is a single point, $\cL_x=\cL|_{x\tm\Jh}$, and similarly for
$\wh{Z}=\xh$. We will often use the following formulae:
$$\cL_{\wh{x}}|_{\Ta}\cong j_{\a}^{*}\cO_C(x)\quad
\cL_{x}|_{\Thb}\cong j_{\bh}^{*}\cO_C(x).\tag{1.10}$$
Here we denote by $x$ both a divisor on $C$ and its linear equivalence class 
(an element of $J$); this will be a habit throughout the paper. We reserve the
notation $[x]$ for the line-bundle on $C$ corresponding to the invertible sheaf
$\cO_C(x)$. 

For $\a\in J$ let 
$$\text{$\Jh[2]_{-{\hat\a}}$ be the translation of $\Jh[2]$ by
$(-{\hat\a})$,}\tag{1.11}$$
 where $\Jh[2]$ is the kernel of multiplication
by two, and 
$$\text{$\nua\cl\Ia\to\Jh$ be the blow up of
$\Jh[2]_{-{\hat\a}}$.}\tag{1.12}$$ 
We let $E$ be the exceptional divisor of $\nua$, and  $E_1,\ldots,E_{16}$ 
be its irreducible components. 

Let $X$ be a projective variety and $D$ an ample divisor on $X$. 
A torsion-free sheaf $F$ on $X$ is {\it $D$-semistable} if it is
Gieseker-Maruyama semistable with respect to $D$, i.e.~if for all proper
subsheaves $E\ss F$   
$$\rk F\cdot\chi(E(nD))\le\rk E\cdot\chi(F(nD)),\quad \text{for all $n\gg 0$}.
\tag{1.13}$$
If there exists $E\ss F$ such that the inequality is an equality then $F$ is
{\it strictly semistable}, otherwise it is stable. There is also the notion of
{\it slope-(semi)stability}: if for all $E\ss F$ with $0<\rk E<\rk F$
$$\mu(E):={1\over \rk E}c_1(E)\cdot D^{k-1}\le
{1\over \rk F}c_1(F)\cdot D^{k-1}=:\mu(F),\quad k:=\dim X,$$
$F$ is $D$-slope semistable. It is $D$-slope stable if the inequality is always
strict. Writing out the polynomials appearing in~(1.13) one shows that
$D$-semistability implies $D$-slope semistability, and $D$-slope
stability implies $D$-stability. Throughout the paper we fix the ample divisor
$\T$ on $J$: (semi)stability of a sheaf will be $\T$-(semi)stability, and
similarly for slope-(semi)stability.

We recall that moduli spaces of semistable torsion-free sheaves parametrize
{\it S-equivalence} classes of such sheaves~\cite{G}. To define
$S$-equivalence one associates to a semistable sheaf $F$ a direct sum of stable
sheaves $Gr(F)$, and then declares that $F_1$ is S-equivalent to $F_2$ if
$Gr(F_1)\cong Gr(F_2)$. If $\rk F=2$ (the only case to be considered in this
paper)  
$$Gr(F)=
\cases
F & \text{if $F$ is stable}\\
L\op(F/L) & \text{if $F$ is strictly semistable, $L\ss F$ destabilizes.}
\endcases
$$
If $F$ is a semistable sheaf we let $[F]$ be its $S$-equivalence class.

A {\it family of sheaves on $X$ parametrized by $B$} is a sheaf $\cF$ on $X\tm
B$, flat over $\cO_B$. For $b\in B$ we set $\cF_b:=\cF|_{X\tm\{b\}}$. 
\head
2. Construction of $\cMt$
\endhead
This subsection must be read with the aid of~\cite{O1}.
\subhead
2.1. Symplectic desingularization of $\Mv$
\endsubhead
We will show that the construction  of a symplectic
desingularization of $\cM_4$, the moduli space of rank-two torsion-free
sheaves on a $K3$ surface $X$ with $c_1=0$ and $c_2=4$ (the ample divisor
defining (semi)stability must satisfy~(0.2) of~\cite{O1}, an assumption
analogous to~(1.3)), carries over to give a symplectic desingularization of
$\Mv$.  As in~\cite{O1} we start by classifying  strictly semistable sheaves.  Let
$[F]\in\Mv$, and assume  
$$0\to L_1\to F\to L_2\to 0\tag{2.1.1}$$
is a destabilizing sequence. 

\proclaim{(2.1.2)Lemma}
Keep notation as above. Then $L_i\cong I_{x_i}\ot\xi_i$, where $I_{x_i}$ is the
ideal sheaf of a point $x_i\in J$, and $\xi_i\in\Jh$. Conversely, if $F$ fits into
Exact Sequence~(2.1.1) with $L_i$ of this form, then $[F]\in\Mv$ and $F$ is
strictly semistable.
\endproclaim
\demo{Proof}
Since $\rk F=2$, both $L_1$ and $L_2$ are rank-one sheaves. They are both
torsion-free. This is clear for $L_1$, because it is a subsheaf of the torsion-free
sheaf $F$. If $L_2$ had non-zero torsion $T$, the surjection of $F$ to $L_2/T$
would desemistabilize $F$. Thus $L_i\cong I_{Z_i}\ot\xi_i$, where $Z_i$ is a
zero-dimensional subscheme of $J$, and $\xi_i$ is a line-bundle on $J$.
Since~(2.1.1) is destabilizing (but $F$ is semistable),
$$2n^2+2(\xi_i\cdot\T)n+\xi_i\cdot\xi_i-2\ell(Z_i)=
2\chi(L_i(n\T))=\chi(F(n\T))=2n^2-2,\quad\text{all $n$},$$
From $\xi_i\cdot\T=0$ and the Hodge index Theorem we get that
$\xi_i\cdot\xi_i\le 0$. Equating the constant coefficients of the above
polynomials we get that 
$$0\le 2\ell(Z_i)=\xi_i\cdot\xi_i+2.\tag{2.1.3}$$
Hence $-2\le\xi_i\cdot\xi_i$; by Assumption~(1.3) equality cannot hold, and
since the intersection form is even we get that $\xi_i\cdot\xi_i=0$. By the
Hodge index Theorem $\xi_i$ is algebraically equivalent to zero,
i.e.~$\xi_i\in\Jh$. From~(2.1.3) we also get that $\ell(Z_i)=1$. This proves
the first assertion of the lemma. The converse is immediate.   
\qed
\enddemo
If in the above lemma we replace $I_{x_i}(\xi_i)$ by $I_{Z_i}$, where $Z_i$ is a
length-two subscheme of the $K3$ surface $X$, we get Lemma~(1.1.5)
of~\cite{O1}, in the case $n=2$. Moreover, all the results in
Subsections~(1.4)-(1.8) of~\cite{O1} remain valid
if one makes the substitution 
$$\text{($I_Z$ on $X$ with $\ell(Z)=2$)}\longmapsto
\text{($I_x(\xi)$ on $J$ where $x\in J$ and $\xi\in\Jh$)}.$$
Thus we get, as in~(1.8) of~\cite{O1}, a desingularization
$\pih_{\vv}\cl\wh{\cM}_{\vv}\to\Mv$, by blowing up the locus parametrizing the
``worst'' strictly semistable sheaves   
$$\Ov:=\{[I_x(\xi)\op I_x(\xi)]|\ x\in J,\ \xi\in\Jh\},$$
and then blowing up the strict transform of 
$$\Siv:=\{[I_{x_1}(\xi_1)\op I_{x_2}(\xi_2)]|\ x_i\in J,\ \xi_i\in\Jh\}.$$ 
Associated to a non-zero two-form $\om$ on $J$ there is a regular 
Mukai two-form $\wh{\om}_{\vv}$ on $\wh{\cM}_{\vv}$, which degenerates on
$\Oh_{\vv}:=\pih^{-1}(\Ov)$ (copy the proof of~\cite{O1, (2.2.3)}). Thus 
$\wh{\cM}_{\vv}$ is not symplectic. One verifies that the proofs of
Propositions~(2.0.1)-(2.0.3) of~\cite{O1} remain valid if we replace
everywhere $\wh{\cM}_4$, $\Oh_4$, $X^{[2]}$, $I_Z$, by
$\wh{\cM}_{\vv}$, $\Oh_{\vv}$, $(J\tm\Jh)$, $I_x(\xi)$ respectively. Thus
Proposition~(2.0.3) tells us that we can blow down $\Oh_{\vv}$ and get a
projective desingularization $\Mtv$ of $\Mv$; let $\pitv\cl\Mtv\to\Mv$ be
the desingularization map. The two-form $\omtv$ on $\Mtv$ induced by 
$\wh{\om}_{\vv}$  is symplectic.   We set
$\Sit_{\vv}:=\pitv^{-1}(\Si_{\vv})$,
$\wt{\O}_{\vv}:=\pitv^{-1}(\O_{\vv})$. Since $\Sit_{\vv}$ is a Cartier
divisor in $\Mtv$, and $\Mtv\sm\Sit_{\vv}=\pitv^{-1}(\Mv^s)$,
Equation~(1.2) gives that  
$$\text{$\Mtv$ is of pure dimension $10$.}\tag{2.1.4}$$ 
\subhead
2.2. Moduli-theoretic interpretation of $\Mtv\sm\wt{\O}_{\vv}$
\endsubhead
We will show that points of $\Mtv\sm\wt{\O}_{\vv}$ are in one-to-one
correspondence with simple semistable sheaves on $J$, modulo $\wt{\text{\rom
S}}$-equivalence, a relation finer than S-equivalence. Let $F$ be a  torsion-free
simple semistable sheaf $F$ on $J$ with $v(F)=\vv$, where $v(F)$ is as
in~(1.1). If $F$ is stable, then $E$ is $\wt{\text{\rom S}}$-equivalent to $F$
if and only if it is isomorphic to it. If $F$ is strictly semistable then by
Lemma~(2.1.2) there is an exact sequence
$$0\to L_1\to F\to L_2\to 0,$$
where $L_i\cong I_{x_i}(\xi_i)$. We associate to $F$ the extension class
of the above exact sequence 
$$e_F\in\Ext^1(L_2,L_1).$$
This is non-zero because $F$ is simple, and is well-defined modulo $\CC^{*}$,
because the destabilizing surjection $F\to L_2$ is determined up to $\CC^{*}$. 
Yoneda multiplication 
$$\Y\cl\Ext^1(L_1,L_2)\tm\Ext^1(L_2,L_1)\to\Ext^2(L_1,L_1)
\tag{2.2.1}$$
can be identified with Serre duality, because the trace
$\Tr\cl\Ext^2(L_1,L_1)\to H^2(\cO_J)$  is an isomorphism. Thus $\Y$ is a perfect
pairing.  Since $F$ is simple, $L_1\not\cong L_2$ and hence Riemann-Roch gives
$\dim\Ext^1(L_1,L_2)=2$. Thus the annihilator (with respect to $\Y$) of $e_F$ is
one-dimensional; let  $e^{\bot}_F\in\Ext^1(L_1,L_2)$ be a generator. 
A sheaf $E$ is $\wt{\text{\rom S}}$-equivalent to $F$ if it is isomorphic to
$F$ or to the extension  
$$0\to L_2\to E\to L_1\to 0\tag{2.2.2}$$
determined by $e^{\bot}_F$. Let $\Ev$ be the set of $\wt{S}$-equivalence
classes of torsion-free simple semistable sheaves $F$ on $J$ with
$v(F)=\vv$. We will show that there is a natural bijection  
$$\psit_{\vv}\cl(\Mtv\sm\Ot_{\vv})\brel\sim\over\lra\Ev.\tag{2.2.3}$$
The desingularization map $\pitv\cl\Mtv\to\Mv$ identifies
$(\Mtv\sm\Sit_{\vv})$ with $\Mv^{st}$. The latter is the set of isomorphism
classes of torsion-free stable sheaves $F$ with $v(F)=\vv$, which injects into
$\Ev$. This injection defines $\psit_{\vv}$ outside $\Sit_{\vv}$. Now we define
$\psit_{\vv}$ on $(\Sit_{\vv}\sm\Ot_{\vv})$. Since the map $\Mhv\to\Mtv$ is the
contraction of $\Oh_{\vv}$, it defines an isomorphism 
$$(\Sih_{\vv}\sm\Oh_{\vv})\brel\sim\over\lra
(\Sit_{\vv}\sm\Ot_{\vv})\tag{2.2.4}$$
commuting with $\pih_{\vv}$, $\pit_{\vv}$. We will give an injection of
$(\Sih_{\vv}\sm\Oh_{\vv})$ into $\Ev$. We use the notation of ~\cite{O1,\S 1},
adapted to our moduli space. Thus $Q$ is the Quot-scheme of which $\Mv$ is the
G.I.T. $\PGL(N)$-quotient. For $x\in Q$ we let  
$$\cO_J(-k)^{N}\to F_x$$
be the quotient parametrized by $x$ (see~\cite{O1,(1.1)}). Let
$$\Si_Q^0:= \{x\in Q^{ss}|\ \text{$F_x\cong L_1\op L_2$, $L_1\not\cong L_2$}\},
\tag{2.2.5}$$
and $\Si_Q$ be its closure in $Q$. For $x\in\Si_Q^0$ with $F_x\cong L_1\op L_2$
we have by~\cite{O1,(1.4.1)} a canonical isomorphism 
$$(C_{\Si_Q}Q)_x\cong
\{(\e,\eta)\in\Ext^1(L_1,L_2)\op\Ext^1(L_2,L_1)|\ 
\e\cup\eta=0\},\tag{2.2.6}$$ 
where $C_{\Si_Q}Q$ is the normal cone to $\Si_Q$ in $Q$. 
Since $\Mhv$ is the $\PGL(N)$-quotient of the variety $S$ obtained from $Q$ by
first blowing up 
$$\O_Q:=\text{closure of $\{x\in Q^{ss}|\ F_x\cong L\op L\}$},$$  
and then the strict transform of $\Si_Q$, we have
$$\pih^{-1}([L_1\op L_2])\cong\PP N(L_1,L_2)//\CC^{*},
\qquad L_1\not\cong L_2\tag{2.2.7}$$ 
where $N(L_1,L_2)$ is the right-hand side of~(2.2.6),  and the action
of $\l\in\CC^{*}$ is given by $\l(\e,\eta)=(\l\e,\l^{-1}\eta)$
(see~\cite{O1,(1.4.1)}).  Let 
$$C(L_1,L_2):=
\{\vf\in\Ext^1(L_2,L_1)^{*}\ot\Ext^1(L_2,L_1)|\ 0=\Tr\vf=\Det\vf\}.$$  
As is easily verified
$$\PP N(L_1,L_2)//\CC^{*}\cong\PP C(L_1,L_2),\tag{2.2.8}$$
with quotient map  given by
$$\matrix
\PP N(L_1,L_2)^{s} & \lra & \PP C(L_1,L_2),\\
[(\e,\eta)] & \mapsto & [\e\ot \eta].
\endmatrix
\tag{2.2.9}$$
(Here $\e\in\Ext^1(L_1,L_2)\cong\Ext^1(L_2,L_1)^{*}$, the isomorphism being
given by~(2.2.1).) 
We have an injection
$$\matrix
\PP C(L_1,L_2) & \hra & \Ev\\
[\vf] & \mapsto & ([\Ann(\Ker\vf)],[\Im\vf]),
\endmatrix$$
and this defines $\wt{\psi}_{\vv}$ on $\pitv^{-1}([L_1\op
L_2])\cong\pih_{\vv}^{-1}([L_1\op L_2])$. We have defined $\psit_{\vv}$ on all of
$(\Mtv\sm\Ot_{\vv})$; one checks immediately that it  is a bijection. Our next
task is to show that $\psit_{\vv}$ is regular on parameter spaces for simple
sheaves. 
 Let $\cE$ be a family of torsion-free simple 
semistable sheaves on $J$ parametrized by a scheme $T$, with
$v(\cE_t)=\vv$ for all $t\in T$, and let
$$\matrix
T & \brel\psi_{\cE}\over\lra & \Mv\\
t & \mapsto & [\cE_t]
\endmatrix$$
 be the modular map. Since $\cE$ is a family of simple semistable sheaves, the
image of $\psi_{\cE}$ is contained in $(\Mv\sm\Ov)$.

\proclaim{(2.2.10)Proposition}
Keep notation as above. There exists a lift
$$\psit_{\cE}\cl T\to(\Mtv\sm\Ot_{\vv})$$
of $\psi_{\cE}$ with the following properties.
\roster
\item
The induced map of sets is (given Identification~(2.2.3)):
$$\psit_{\cE}(t)= \text{$\wt{S}$-equivalence class of $\cE_t$}.$$
\item
Let $0\in T$. If the map from the germ of $T$ at $0$ to the
 deformation space of $\cE_0$ is an isomorphism, then
$\psit_{\cE}$ is an isomorphism in a neighborhood of $0$ (for the classical
topology).  
\endroster
\endproclaim

\demo{Proof}
Clearly it suffices to prove the
proposition for $\cE$ satisfying the hypothesis of Item~(2); in this case the germ
$(T,t)$ is the deformation space of $\cE_t$ (after shrinking $T$ if necessary),
hence it suffices to prove that $\psi_{\cE}$ lifts in  a neighborhood of $0$.  Let
$E:=\cE_0$. Assume $E$ is stable: it is well-known that $\psi_{\cE}$ is an
isomorphism near $0$, and furthermore since $\pit\cl\Mtv\to\Mv$ is an
isomorphism over $\Mv^{st}$,  the map $\psit_{\cE}$ lifts trivially and it has the
required properties. Now assume $E$ is strictly semistable, and let
$$0\to L_1\to E\to L_2\to 0\tag{2.2.11}$$
be the destabilizing sequence of $E$. 
By Serre duality and simpleness of $E$,     
$$\Ext^2(E,E)\cong\Ext^0(E,E)^{*}\cong\CC,$$
hence the deformation space of $E$ is smooth by Mukai~\cite{Muk2,(0.1)};
thus $T$ is also  smooth   (after shrinking $T$ around $0$ if necessary).  Let 
$$\Si_S^0:=(\pi_R\circ\pi_S)^{-1}(\Si_Q^0),$$
where $\pi_R\cl R\to Q$ is the blow up of $\O_Q$, $\pi_S\cl S\to R$ is the
blow up of  the strict transform of $\Si_Q$ (see~\cite{O1,(1.1)}), and $\Si_Q^0$ is
as in~(2.2.5). Let $x_0\in Q^{ss}$ be such that $F_{x_0}\cong L_1\op
L_2$, where $L_1$, $L_2$ are the sheaves fitting into~(2.2.11). Let
$$y_0:=[e_E^{\bot},e_E]\in
\PP(N(L_1,L_2))=(\pi_R\circ\pi_S)^{-1}(x_0),\tag{2.2.12}$$ 
where $e_E$ is the extension class of~(2.2.11) (the second equality follows
from~(2.2.6)).   The point $y_0$ is stable for the $\CC^{*}$-action on
$\PP(N(L_1,L_2))$,  i.e.~it is stable with respect to the $\PGL(N)$-action on $S$.
Let $U\ss S^{s}$ be ``Luna's \'etale slice'' normal to the orbit $\PGL(N)y_0$. Then
(see~\cite{O1,(1.8.8)}) 
$$\Stab({y_0})\cong\ZZ/(2)$$
acts on $U$. Let $\wh{\cF}$ be the sheaf on $J\tm U$ obtained pulling back the
tautological family of quotients parametrized by $Q$ via the map $U\to Q$.
Shrinking $U$ around $y_0$  (in the classical
topology), we may assume that it satisfies the following properties. Let 
$\Si_S$ be the closure of $\Si_S^0$, and $\Si_U:=\Si_S\cap U$; there exist families
of rank-one torsion-free sheaves $\cL_1$, $\cL_2$ on $J$ parametrized by $\Si_U$
and an exact sequence 
$$0\to\cL_2\to\wh{\cF}|_{J\tm\Si_U}\to\cL_1\to 0,\tag{2.2.13}$$
such that:
\roster   
\item"{(I)}"
the restriction of~(2.2.13) to $J\tm\{y\}$ 
$$0\to L_2(y)\brel f_2(y)\over\lra \wh{\cF}
\brel f_1(y)\over\lra L_1(y)\to 0$$
is a destabilizing sequence of $\wh{\cF}_y$ for all $y\in\Si_U$,
\item"{(II)}"
for $y=y_0$ we have $L_2(y_0)=L_2$, $L_1(y_0)=L_1$.
\endroster
Since $\Si_U$ is a Cartier divisor we can construct the elementary
modification of $\wh{\cF}$ associated to~(2.2.13), i.e.~the sheaf $\cG$ on
$J\tm U$ fitting into the exact sequence  
$$0\to\cG\to\wh{\cF}\to i_{*}(\cL_1)\to 0,\tag{2.2.14}$$
where $i\cl J\tm\Si_U\hra J\tm U$ is the inclusion. Then $\cG$ is flat over
$U$. If $y\in (U\sm\Si_U)$ we have $\cG_y\cong
\wh{\cF}_y$. 
Now assume $y\in\Si_U$, and let $x:=\pi_R\circ\pi_S(y)$. By~(2.2.6) the
point $y$ is an element of $\PP(N(L_1,L_2))$, which we write explicitely as
$$y=[\e(y),\eta(y)]\in\PP(\Ext^1(L_1(y),L_2(y))\op\Ext^1(L_2(y),L_1(y))).$$
Tensorizing~(2.2.14)
by $\cO_{U}/m_y$ ($m_y$ is the maximal ideal of $y\in U$) one gets an exact
sequence 
$$0\to L_1(y)\to \cG_y\to L_2(y)\to 0,\tag{2.2.15}$$
We claim that 
$$\text{$\eta(y)$ is the extension class of~(2.2.15).}\tag{2.2.16}$$
 To prove this, let 
$$\k\cl T_yS\to \Ext^1(\wh{\cF}_y,\wh{\cF}_y)$$
be the Kodaira-Spencer map, and $v(y)\in T_y S$ be a vector normal to 
$\Si_U$; according to~\cite{O4,(1.11)} the extension class of~(2.2.15) is
equal to 
$$f_1(y)\circ\k(v(y))\circ f_2(y).$$
Clearly this equals $\eta(y)$. Since $y\in U\ss S^{s}$, we have
$\eta(y)\not=0$, and hence $\cG$ is a family of simple sheaves. Furthermore
by~(2.2.16) and~(2.2.12) $\cG_{y_0}\cong E$.  Hence, shrinking $U$ if
necessary, the family $\cG$ is the pull back of $\cE$ by a map $g\cl U\to T$ with
$g(y_0)=0$. Since $E$ is simple $g(y_1)=g(y_2)$ if and only if
$\cG_{y_1}\cong\cG_{y_2}$. As is easily verified this is the case if and only if
$y_1$ and $y_2$ are equivalent for the $\ZZ/(2)$-action on $U$. Thus we get an
injective map 
$$\ov{g}\cl U/\ZZ(2)\to T,\qquad \ov{g}(y_0)=0.$$
Both $U$ and $T$ have the same dimension as $\Mhv$, so $\dim U=\dim T$. Since
$T$ is smooth we get that $\ov{g}$ is an isomorphism near $y_0$. On the other
hand, restricting the quotient map $S^{s}\to\Mhv$ to $U$ we get 
$$\rho\cl U/\ZZ(2)\hra(\Mhv\sm\Oh_{\vv})\cong(\Mtv\sm\Ot_{\vv}),$$
which is an isomorphism onto an open neighborhood of $\rho(y_0)$ (by Luna's
\'etale slice Theorem). Furthermore $\rho(y_0)\in(\Sit_{\vv}\sm\Ot_{\vv})$
corresponds to $\cG_{y_0}\cong E$ via~(2.2.3). Inverting $\ov{g}$ near $0$
and composing with $\rho$ we get a lift of $\psi_{\cE}$ with the required
properties. \qed
\enddemo
\subhead
2.3. The symplectic variety $\cMt$ 
\endsubhead
Recall that 
$$\cM:=\av^{-1}(0,\wh{0}),\quad\atv:=\av\circ\pitv,\quad
\cMt:=\atv^{-1}(0,\wh{0}),\quad\omt:=\omtv|_{\cMt}.$$ 
Let us prove that $\cMt$ is smooth of pure dimension $6$. A straightforward
computation shows that the map
$$\matrix
\cM\tm J\tm\Jh & \brel q\over\lra & \Mv \\
([F],x,\xi) & \mapsto & [\xi\ot t_x^{*}F],\\
\endmatrix$$
where $t_x\cl J\to J$ is translation by $x$, is the quotient map for the action of
$J[2]\tm\Jh[2]$ given by
$$([F],x,\xi)^{(x_0,\xi_0)}:=
([\xi_0\ot t_{{x_0}}^{*}F],x-x_0,\xi\ot\xi_0^{-1}).$$ 
We claim that there is a $J[2]\tm\Jh[2]$-quotient map 
$$\cMt\tm J\tm\Jh\brel\wt{q}\over\lra \Mtv\tag{2.3.1}$$
covering $q$. To prove it let
$\wh{\cM}:=\pih_{\vv}^{-1}(0,\wh{0})$. There exists
$$\wh{\cM}\tm J\tm\Jh\brel\wh{q}\over\lra \wh{\cM}_{\vv}$$
covering $q$, because $\pih$ is obtained by first blowing
up $\Ov$ and then the strict transform of $\Siv$, and 
$$q^{*}\Ov=\O\tm J\tm\Jh\qquad q^{*}\Siv=\Si\tm J\tm\Jh,$$
where $\O:=\Ov\cap\cM$, $\Si:=\Siv\cap\cM$.
Furthermore, since the loci we blow up are $J[2]\tm\Jh[2]$-stable, 
the action of $J[2]\tm\Jh[2]$ on $\cM\tm J\tm\Jh$ lifts to an
action on $\wh{\cM}\tm J\tm\Jh$, and $\wh{q}$ is the quotient
for this action. Since $\wh{\cM}_{\vv}\to\Mtv$ is the contraction of
$\wh{\O}_{\vv}$, and $\wh{q}^{*}\wh{\O}_{\vv}=\wh{\O}\tm J\tm\Jh$ (here 
$\wh{\O}:=\wh{\O}_{\vv}\cap \wh{\cM}$), the map $\wh{q}$ descends to a map
$\wt{q}$ as in~(2.3.1);  one verifies that $\wt{q}$ is the quotient map for a
lift of the $J[2]\tm \Jh[2]$-action. Now since $\wt{q}$ is \'etale and
by~(2.1.4) the moduli space $\Mtv$ is smooth of pure dimension $10$, 
$$\text{$\cMt$ is smooth of pure dimension
$6$.}\tag"$\text{(2.3.2)}^{*}$"$$ 

\proclaim{(2.3.3)Proposition}
Keep notation as above. Then $\cMt$ and $J\tm\Jh$ are orthogonal for
$\wt{q}^{*}\omtv$.
\endproclaim

\demo{Proof}
It suffices to prove orthogonality on the dense subset $\wt{q}^{-1}\Mv^{st}$.
The composition 
$$\wt{q}^{-1}(\Mv^s)\brel\wt{q}\over\lra\Mtv\brel\pit\over\lra \Mv^s$$
is \'etale hence the differential at 
 $z\in\wt{q}^{-1}\Mv^{st}$ gives an  isomorphism
$$d(\pit\circ\wt{q}(z))\cl 
T_z(\Mtv^0\tm J\tm\Jh)\brel\sim\over\lra \Ext^1(E,E),\tag{2.3.4}$$ 
where $[E]=\pit\circ\wt{q}(z)$. On the other hand we have isomorphisms
$$\wt{q}^{-1}(\Mv^s)=\pit^{-1}\cM^{st}\tm J\tm\Jh
\brel\pit\tm\id\over\lra\cM^{st}\tm J\tm\Jh.$$
Let $(\pit\tm\id)(z)=([F],x,\xi)$ (thus $E\cong\xi\ot t_x^{*}F$); by the above
isomorphism we have a decomposition 
$$T_z(\cMt\tm J\tm\Jh)\cong T_{[F]}\cM\op T_x J\op T_{\xi}\Jh.$$
 Restricting~(2.3.4) to  the direct summands we get 
$$\align
T_{[F]}\cM\op T_x J & \brel\sim\over\lra \Ext^1(E,E)^0, \\
T_{\xi}\Jh& \brel\d\over\lra \Ext^1(E,E), 
\endalign $$
where $\d$ is the ``diagonal embedding'' $H^1(\cO_J)\hra \Ext^1(E,E)$. 
By definition of Mukai's form (see~(1.9) of~\cite{O1}),
$$\wt{q}^{*}\omt(z)(\a,\b)=
\int\limits_J \om\wedge \text{\rm Tr}(d\wt{q}(z)(\a)\wedge d\wt{q}(z)(\b)),$$
where the second wedge denotes Yoneda product.  If
$d\wt{q}(z)(\a)\in\Ext^1(E,E)^0$ and $d\wt{q}(z)(\b)\in\d(H^1(\cO_J))$, then
$$\text{\rm Tr}(d\wt{q}(z)(\a)\wedge d\wt{q}(z)(\b))=0.$$
Thus $T_{[F]}\cM$ is orthogonal to $T_{[\xi]}\Jh$. Now let's prove
orthogonality to $T_x J$. Let $\psi\cl U\to \Mv^{st}$ be a neighborhood of $[E]$ 
in the \'etale topology such that there is a tautological sheaf $\cE$ on $J\tm U$,
and let $0\in U$ be a point mapping to $[E]$. We can represent 
$\psi^{*}\omt$ as the Mumford two-form induced by a representative of
$c_2(\cE)\in CH^2(J\tm U)$. More precisely, let $\cZ$ be a
cycle representing $c_2(\cE)$ and intersecting $J\tm\{0\}$ transversely.
Thus  
$$\cZ|_{J\tm\{0\}}=\sum\limits_{i}\e_i p_i,$$
where the $p_i$ are pairwise distinct, and $\e_i=\pm 1$. Furthermore, by
transversality $\cZ$ defines for each $i$ a map $f_i\cl T_0 U\to T_{p_i}J$.
By~\cite{O3,(2.9)} we have, for $\a,\b\in T_0 U$,  
$$\psi^{*}\omt(\a,\b)=-{1\over 4\pi^2}\sum\limits_{i}\e_i
\om(f_i(\a),f_i(\b)).\tag{2.3.5}$$ 
Let
$$\align
\a & \in d\psi(0)^{-1}(d\wt{q}(z)(T_{[F]}\cM)),\\
\b & \in d\psi(0)^{-1}(d\wt{q}(z)(T_x J)).
\endalign$$
Then $\sum\limits \e_i f_i(\a)=0$ (in a trivialization of the tangent bundle of
$J$), and $f_i(\b)$ is independent of $i$. Hence the right-hand side
of~(2.3.5) vanishes.     
\qed
\enddemo 

\proclaim{$\text{(2.3.6)}^{*}$Corollary}
Keeping notation as above, $\omt$ is a symplectic form on $\Mt$.
\endproclaim

\demo{Proof}
Since $\omtv$ is non-degenerate and $\wt{q}$ is \'etale, the two-form
$\wt{q}^{*}\omtv$ is non degenerate. By~(2.3.3) we get that $\omt$
is non-degenerate. 
\qed
\enddemo

There is an analogue of Proposition~(2.2.10) valid for $\Mt$, which follows
immediately from~(2.2.10). Let $\cE$ be a family of torsion-free simple
semistable sheaves on $J$ parametrized by a scheme $T$, with
$v(\cE_t)=\vv$, $\det(\cE_t)\cong\cO_J$, $\sum c_2(\cE_t)=0$ for all
$t\in T$. Let $\psi_{\cE}\cl T\to\cM$ be the modular map. 

\proclaim{(2.3.7)Proposition}
Keep notation as above. There exists a lift
$$\psit_{\cE}\cl T\to(\Mt\sm\Ot)$$
of $\psi_{\cE}$ with the following properties.
\roster
\item
The induced map of sets is (given Identification~(2.2.3)):
$$\psit_{\cE}(t)= \text{$\wt{S}$-equivalence class of $\cE_t$}.$$
\item
Let $0\in T$. Assume the map from the germ $(T,0)$
to the locus in $\Def(\cE_0)$ parametrizing sheaves with trivial
determinant and $0$ Albanese map is an isomorphism. Then the map
$\wt{\psi}_{\cE}$ is an isomorphism in a neighborhood of $0$ (in the classical
topology).  
\endroster
\endproclaim

\proclaim{(2.3.8)Notation}
{\rm Let $\Sit:=\pit^{-1}\Si$, i.e.
$$\Sit:=
\{x\in\Mt|\ \text{$\pit(x)$ is represented by a strictly semistable sheaf}\}.$$
Let
$$B:=\text{closure of $\{[F]\in\Mt|\ \text{$F$ is singular and stable}\}$,}$$
 and $\Bt\ss\Mt$ be the strict transform of $B$.} 
\endproclaim

\proclaim{(2.3.9)Remark}
{\rm Let $y_0\in(\Sit\sm\Ot)$. The correspondence~(2.2.3) associates to
$y_0$ two isomorphism classes of simple sheaves on $J$, a sheaf $E$ fitting into 
$$0\to L_1\to E\to L_2\to 0,$$
with extension class $e$, and a sheaf $G$ fitting into 
$$0\to L_2\to G\to L_1\to 0,$$
with extension class $e^{\bot}$. Let $T$ and $\cE$ be as in~(2.3.7), and
assume $\cE_0\cong E$. Then by~(2.3.7) a neighborhood of $y_0$ can be
identified with a neighborhood of $0\in T$. On the other hand, by the same
proposition,  a neighborhood of $0\in T$ must also parametrize a family
$\cG$ of deformations of $G$. 
What is the realtion between $\cG$ and $\cE$ ?One passes from one to the other
by means of an elemntary modification.
In order to explain this, let
$$\Sit(T):=\{t\in T|\ \text{$\cE_t$ is strictly semistable}\}.$$
Then (shrinking $T$ in the classical topology, if necessary) there is an exact
sequence
$$0\to\cL_1\to\cE|_{J\tm\Sit(T)}\to \cL_2\to 0,$$
which for every $t\in\Sit(T)$ restricts to the destabilizing sequence for $\cE_t$.
One verifies that $\cG$ is the sheaf on $J\tm T$ fitting into the exact sequence
$$0\to\cG\to\cE\to i_{*}\cL_2\to 0,$$
where $i\cl J\tm\Sit(T)\hra J\tm T$ is the inclusion.}
\endproclaim
\head
3. An application of the generalized Lefschetz Hyperplane Theorem.
\endhead
For $\a\in J$ let
$$\align
Z_{\a}:=&\{[F]\in\cM|\ \text{$F|_{\Ta}$ is not locally-free
semistable}\},\\
\Zta:=&\pit^{-1}(Z_{\a}).\tag{3.1}
\endalign
$$ 
Applying the generalized Lefschetz
hyperplane Theorem~\cite{GM} to  the determinant map of $\cMt$, and
arguing as in~\cite{O1,(3.0.1)} one gets the following result.

\proclaim{$\text{(3.2)}^{*}$Proposition}
Keep notation as above, and let $i\cl\Zta\hra\cMt$ be the inclusion. The
map
$$i_{\#}\cl \pi_q(\Zta)\to\pi_q(\cMt)$$
is an isomorphism for $q\le 1$ and a surjection for $q=2$. In particular
$$H_2(\Zta;\ZZ)\brel i_{*}\over\lra H_2(\cMt;\ZZ)\qquad
H^2(\cMt;\ZZ)\brel i^{*}\over\lra H_2(\Zta;\ZZ)$$
are surjective and injective, respectively.
\endproclaim

Let 
$$\align
\Si_{\a}:= & \Si\cap Z_{\a},\tag{3.3}\\
B_{\a}^0:= & \{[F]\in \cM|\ \text{$F$ is stable and $\ia^{*}F$ is singular}\},
\tag{3.4}\\ 
V_{\a}^0:= & \{[F]\in \cM|\ \text{$\ia^{*}F$ is locally-free, not semistable}\},
\tag{3.5} 
\endalign$$
and let $\Ba$, $\Va$ be the closures of $\Ba^0$ and $\Va^0$ 
respectively. Let $\Sta:=\pit^{-1}(\Si_{\a})$, and let  $\Bta,\Vta\ss\cMt$
be the strict transforms of $B_{\a}$ and $V_{\a}$, respectively. Clearly we
have a decomposition into closed subsets 
$$\Zta=\Sta\cup\Bta\cup\Vta.$$

\proclaim{(3.6)Claim}
Keep notation as above. If
$$\a\notin\bigcup\limits_{x\in J[2]}\T_{x},\tag{3.7}$$
then $\Zta\ss(\cMt\sm\Ot)$.
\endproclaim

\demo{Proof}
Assume $\Zta\cap\Ot\not=\es$. Let $z\in\Zta\cap\Ot$, and let $[F]=\pit(z)$.
By~(2.1.2) we have an exact sequence
$$0\to I_x\ot\xi\to F\to I_x\ot\xi\to 0,$$
for some $x\in J$, $\xi\in\Jh$. Since $[F]\in Z_{\a}$, we must have $x\in\Ta$,
i.e.~$\a\in\T_{x}$. We have $x\in J[2]$ because $\sum c_2(F)=0$; this proves the
claim.
\qed
\enddemo 

Thus in analyzing $\Zta$ we will be able to use the moduli-theoretic
interpretation of $(\cMt\sm\Ot)$ given by Proposition~(2.3.7). 

\proclaim{(3.8)Remark}
{\rm Let $\atv\cl\Mtv\to J\tm\Jh$ be as in \S $1$ and 
$$\align
\Mt^J &:=\atv^{-1}(J\tm\{{\hat 0}\})\\
\Zta^J & :=\wt{\pi}_{\vv}^{-1}\{[F]\in a_{\vv}^{-1}(J\tm\{{\hat 0}\})|\
\text{$F|_{\Ta}$ is not locally-free semistable}\}.
\endalign$$
Applying the generalized Lefschetz Hyperplane Theorem to the determinant
map on $\Mt^J$ we get that the inclusion $\Zta^J\hra\Mt^J$
induces isomorphisms on $\pi_q$ for $q\le 3$, and a surjection on $H_4$.
The restriction of $\atv$ to $\Mt^J$ is a fibration
over $J$ with fiber $\Mt$, and one verifies easily that 
$\Mt^J$ is cohomologically the product of $J$ and
$\Mt$. Hence if $\Zta^J$ were a cohomological product of
$\Zta$ and $J$ we would know that $\Zta\hra\Mt$ induces isomorphims
on $H^q\ot\QQ$ for $q\le 3$. Now $\Zta^J$ is not a
fibration: the generic fiber is homeomorphic to $\Zta$ but special fibers are
not, however the equality
$b_2(\Mt)=b_2(\Zta)$ that we will prove, shows that up to $H^2$
the cohomology of $\Zta^J$ is the product of that of $\Zta$ and of $J$.} 
\endproclaim
\head
4. Analysis of $\Sta$.
\endhead
The $*$-red tags not contained in Subsection~(4.4) are~(4.1.1), (4.3.2)
and~(4.3.10). We choose $\a$ so that~(3.7) holds. (Co)Homology is with
rational coefficients. 
\subhead
4.1. Low-dimensional topology of $\Sta$ 
\endsubhead
Let $f:=\pit|_{\Sta}$.
By Claim~(3.6) we have $\Sta\ss(\cMt\sm\Ot)$, hence $f$
 is a $\PP^1$-fibration over $\Sia$: 
$$\CD
\PP^1 @>>> \Sta \\
@.                  @VV f V \\
{}         @.       \Sia.
\endCD
\tag"$\text{(4.1.1)}^{*}$"$$
Thus 
$$\align
H^q(\Sta) & \cong H^q(\Sia), \quad q=0,1 \tag{4.1.2}\\ 
H^2(\Sta) &\cong H^2(\Sia)\op\QQ c_1(\om_f),\tag{4.1.3}
\endalign$$
where $\om_f$ is the relative cotangent bundle of $f$.   
\subhead
4.2. Low-dimensional topology of $\Sia$  
\endsubhead
The map $\k\cl C\tm \Jh \to \Sia$ defined by
$$\k(p,\wh{x}) := [I_{i_{\a}(p)}\ot\cL_{\wh{x}}\op
I_{-i_{\a}(p)}\ot\cL_{-\wh{x}}] $$ 
is surjective. It is not injective, because of points which are both in
$\Ta$ and $-\Ta$. We claim that
$$\Ta\cap(-\Ta)=\Ta\cap\T_{-\a}=\{\ia(q_1),\ia(q_2)\}
\qquad q_1\not=q_2.\tag{4.2.1}$$ 
In fact $\Ta\not=\T_{-\a}$ because by~(3.7) $2\a\not=0$. Thus
$i_{\a}^{*}(\T_{-\a})=q_1+q_2$, where
$$q_1+q_2\sim K_C-2\a.\tag{4.2.2}$$
It follows from~(3.7) that $q_1\not= q_2$. Clearly
$\k(q_1,\wh{x})=\k(q_2,-\wh{x})$. One easily shows that 
$$\Sia\cong (C\tm \Jh)/\equiv\tag{4.2.3}$$
where $\equiv$ is the equivalence relation generated by setting
$$(q_1,\wh{x})\equiv(q_2,-\wh{x}).\tag{4.2.4}$$
In order to discuss the topology of
$\Sia$ we introduce some notation. Denote by $\ov{C}$ the quotient of $C$ obtained
identifying $q_1$ with $q_2$, and let 
$$\matrix
\ov{C} & \brel \mu\over\lra & \Sia\\
q & \mapsto & [I_{\ia(q)}\op I_{-\ia(q)}].
\endmatrix$$

\proclaim{(4.2.5)Proposition}
Keep notation as above. 
\roster
\item 
$H^0(\Sia;\QQ)\cong\QQ$
\item
The map $\mu^{*}\cl H^1(\Sia;\QQ)\brel\sim\over\lra H^1(\ov{C};\QQ)$ is
an isomorphism. 
\item
The map $\k^{*}\cl H^2(\Sia;\QQ)\to H^2(C\tm \Jh;\QQ)$ is an
isomorphism. 
\endroster
\endproclaim
\demo{Proof}
Item~(1) follows at once from~(4.2.3). Let $C^0:=(C\sm\{q_1,q_2\})$, and
$$\Sia^0:=\k(C^0\tm \Jh)\cong  C^0\tm \Jh.\tag{4.2.6}$$
(Thus $\Sia^0$ is the smooth locus of $\Sia$.) Let
$$\lra H^{q-1}(\Sia^0)\brel d^{q-1}\over\lra H^q(\Sia,\Sia^0)\lra H^q(\Sia)
\lra H^q(\Sia^0)\brel d^q\over\lra\cdots\tag{4.2.7}$$ 
be the long exact sequence of the couple $(\Sia,\Sia^0)$ (with rational
coefficients). Let $\D_1,\D_2\ss C$ be disjoint discs centered at $q_1$, $q_2$
respectively, and let  $B$ be obtained from $\D_1\cup \D_2$ by gluing $q_1$ and
$q_2$. Let $B^0\ss B$ be the complement of the ``center'', 
i.e.~$B^0:=(\D_1\sm\{q_1\})\cup(\D_2\sm\{q_2\})$. By excision and K\"unneth 
$$H^q(\Sia,\Sia^0)\cong H^q(B\tm\Jh,B^0\tm\Jh)
\cong \QQ\e\ot H^{q-1}(\Jh)\op
(\QQ\eta_1\op\QQ\eta_2)\ot H^{q-2}(\Jh),$$
where $\e$ is a generator of $H^1(B,B^0)$ and $\{\eta_1,\eta_2\}$ is
the basis of $H^2(B,B^0)$ dual to the basis of $H_2(B,B^0)$ given
by the push-forward of the local orientation classes of
$(\D_i,\D_i\sm\{q_i\})$. By~(4.2.6) and K\"unneth   
$$H^{q-1}(\Sia^0)\cong H^{q-1}(\Jh)\op H^1(C^0)\ot H^{q-2}(\Jh).$$
Given the above identifications, the boundary map in~(4.2.7) is given by
$$\matrix
H^{q-1}(\Sia^0) & \brel d^{q-1}\over\lra & H^{q}(\Sia,\Sia^0)\\
(\a,\vf\ot\psi) & \mapsto &
(\a+(-1)^q\a,(\int\limits_{\g}\vf)(\eta_1\ot\psi-\eta_2\ot\psi)),
\endmatrix$$
where $\g\in H_1(C^0)$ is the class of a loop around $q_i$ (with
orientation induced by that of $\D_i$). Thus for $q=1$ Exact
Sequence~(4.2.7) gives 
$$0\to H^1(B,B^0)\ot H^0(\Jh) \to H^1(\Sia)\to H^1(C)\to 0,$$
 and this gives Item~(2). For $q=2$ we get
$$0\to \QQ(\eta_1+\eta_2)\ot H^0(\Jh) 
\to H^2(\Sia)\to H^2(\Jh)\op H^1(C)\ot H^1(\Jh)\to 0.$$
This gives Item~(3). 
\qed
\enddemo
\subhead
4.3. The intersection with $\Bta$
\endsubhead
Let $\nu_0\cl \wh{I}_0\to \Jh$ be the blow-up of $\Jh[2]$, i.e.~the case
$\a=0$ of~(1.12). Let $\iota\cl \wh{I}_0\to \wh{I}_0$ be the involution
covering $(-1)$. Let $\approx$ be the equivalence relation on $C\tm
\wh{I}_0$ generated by 
$$(q_1,\wh{x})\approx(q_2,\iota(\wh{x})),\tag{4.3.1}$$
where $q_1,q_2$ are as in~(4.2.1). We will prove the following result.

\proclaim{$\text{(4.3.2)}^{*}$Proposition}
Keep notation as above. Then:
\roster
\item
Given Bijection~(2.2.3),  $\Sta\cap\Bta$ is the subset of $\Sta$
parametrizing $\wt{S}$-equivalence classes of sheaves $F$ which are 
singular at two distinct points.
\item
There is a map ${\tilde \k}\cl(C\tm \wh{I}_0)\to\Sta\cap\Bta$ which
identifies 
 $\Sta\cap\Bta$ with the quotient of $(C\tm \wh{I}_0)$ modulo $\approx$.
\item
Given~(4.2.3) and Item~(2) above, the restriction of $\pit$ to
$\Sta\cap\Bta$ is identified with the map 
$$(C\tm \wh{I}_0/\approx)\lra (C\tm \Jh/\equiv)$$
induced by $\id_C\tm\nu_0$.  
\endroster
\endproclaim 

The proof of Proposition~(4.3.2) will be given after some preliminary results.

\proclaim{(4.3.3)Lemma}
Let $F$ be a simple torsion-free sheaf on $J$ with $v(F)=\vv$, whose
$\wt{S}$-equivalence class belongs to $\Btv$ via~(2.2.3). 
\roster
\item
The double dual $F^{**}$ fits into an exact sequence
$$0\to \xi_1\to F^{**}\to\xi_2\to 0,\tag{4.3.4}$$
where $\xi_i\in\Jh$. In particular $\ell(F^{**}/F)=2$. 
\item
Exact sequence~(4.3.4) is split if and only if
$\xi_1^{-1}\ot\xi_2\not\cong\cO_J$.  
\item
If $[F]\in\Bta$ then $Sing(F)=\{x,-x\}$ where $x\not=-x$.
\endroster
\endproclaim 
\demo{Proof}
Let $E:=F^{**}$. First we prove that 
$$c_2^{hom}(E)=0,\tag{4.3.5}$$
where  $c_q^{hom}$ is the Chern class in cohomology.  Since $E$ is
slope-semistable, Bogomolov's inequality tells us that  $c_2^{hom}(E)\ge 0$.
Thus it suffices to prove that~(4.3.5) holds for sheaves $F$ parametrized by
an open dense subset of $\Bt$. We will prove it for $[F]\in\Bt\sm\Sit$,
i.e.~we assume 
 $F$ is stable. Since $c_2^{hom}(E)<c_2^{hom}(F)$,
because $F$ is singular,  $c_2^{hom}(E)$ equals $0$ or $1$. Let us show
that %
$$\text{if $c_2^{hom}(E)=1$ then  $E$ is not slope-stable.}
\tag{4.3.6}$$
Assume $E$ is slope-stable. Then
$$H^0(E\ot\cL_{\xh})=0=H^0(E^{*}\ot\cL^{-1}_{\xh})^{*}=H^2(E\ot\cL_{\xh}).
\tag{4.3.7}$$
Hence by Riemann-Roch 
$$\text{$\dim H^1(E\ot\cL_{\xh})=1$ for all $\xh\in\Jh$.}\tag{4.3.8}$$
From~(4.3.7) it follows~\cite{Mum,pp.46-55} that
$\phih_{*}(\cL\ot\phi^{*}E)=R^2\phih_{*}(\cL\ot\phi^{*}E)=0$.
Applying Grothendieck-Riemann-Roch we get
$$ch(R^1\phih_{*}(\cL\ot\phi^{*}E))=
-\phih_{*}[ch(\cL)\phi^{*}ch(E)]=1-2\eta_{\Jh},$$
where $\eta_{\Jh}$ is the fundamental  class of $\Jh$. This is absurd because
by~(4.3.8) the sheaf $R^1\phih_{*}(\cL\ot\phi^{*}E)$ is a line-bundle on
$\Jh$. We have proved~(4.3.6). Let
$$0\to \xi_1\to E\to I_Z\ot \xi_2\to 0$$
be a slope-destabilizing sequence, where the $\xi_i$'s are line bundles, and $Z$ is a
zero dimensional subscheme of $J$. Since $E$ is slope-semistable
(because $F$ is) we have $L\cdot \T=0$, so by Hodge index $L\cdot L\le 0$.
From
$$1=c_2^{hom}(E)=\ell(Z)+\int\limits_{J}c_1(\xi_1)\cdot c_1(\xi_2)$$
one easily gets that $\xi_i\cdot \xi_i=0$. Hence by Hodge index $\xi_i\in\Jh$. If
$x$ is the singular point of $F$ the sheaf $I_x\ot \xi_1$ injects into $F$. This
contradicts our assumption that $F$ is stable, and proves~(4.3.5). Let us
prove Item~(1). Since $c_1^{hom}(E)=0=c_2^{hom}(E)$ and $E$ is
slope-semistable, it follows that $E=F^{**}$ fits into Exact Sequence~(4.3.4).
This is well-known ($E$ is a flat vector-bundle), a quick proof is as follows: by 
Grothendieck-Riemann-Roch we get   
$$\sum\limits_{i=0}\limits^{2}(-1)^i ch(R^i\phih_{*}(\cL\ot\phi^{*}E))=
\phih_{*}[ch(\cL)\phi^{*}ch(E)]=2\eta_{\Jh}.$$
On the other hand if $\Hom(\cL_{\wh{x}},E)=0$ for all $x$, then each $R^i$
appearing in the right-hand side of the above formula is zero. Thus  
$\Hom(\cL_{\wh{x_0}},E)\not=0$ for some $x_0\in J$, and this proves Item~(1).  To
prove Item~(2), first notice that if $\xi_1^{-1}\ot\xi_2\not\cong\cO_J$ then every
extension~(4.3.4) is trivial. Thus  we may assume that
$\xi_1\cong\xi_2\cong\xi$, and we must prove that~(4.3.4) is non-split.
Assume the contrary, i.e.~that  $F^{**}\cong\xi\op\xi$. If $x$ is a singular
point of $F$, one easily checks that $I_x\ot\xi$ injects into $F$. Thus, since
$[F]\notin\Ot$ (because $F$ is simple), we must have $-x\not=x$. Since also
$I_{-x}\ot\xi$ injects into $F$, we get that $F\cong I_x\ot\xi\op
I_{-x}\ot\xi$, contradicting the hypothesis that $F$ is simple.  Finally
Item~(3) follows immediately from Assumption~(3.7).    \qed  
\enddemo

Let $[F]\in\Sit\sm\Ot$, where $F$ fits into the exact sequence 
$$0\to I_{x}\ot \xi\brel\a\over\lra F\to I_{-x}\ot \xi^{-1}
\brel\b\over\lra 0,\tag{4.3.9}$$
and is simple. Let $\cE$ be a family of torsion-free simple semistable sheaves on
$J$ parametrized by $T$, satisfying the hypotheses of Proposition~(2.3.7),
with $\cE_0\cong F$. Thus the map $\wt{\psi}_{\cE}\cl T\to\Mt$ identifies a
neighborhood of $0\in T$ with a neighborhood of $[F]\in\Mt$. Let 
$$\align
\D(T) := & \{t\in T|\ \text{$\cE_t$ is singular}\},\\
\Sit(T):= & \{t\in T|\ \text{$\cE_t$ is strictly semistable}\},\\
\Bt(T):= & \text{closure of $\{t\in T|\ \text{$\cE_t$ is singular and stable}\}$.}
\endalign$$
Thus $\D(T)=\Sit(T)\cup \Bt(T)$ (recall Lemma~(2.1.2)). 

\proclaim{$\text{(4.3.10)}^{*}$Lemma}
Keep notation as above. Assume that $x\not= -x$. Then (after shrinking $T$
around $0$ if necessary)  
$$\Sit(T)\cap \Bt(T)=\{t\in \Sit(T)|\ \text{$\cE_t$ is singular at two points}\}.
\tag{4.3.11}$$
Furthermore $\Sit(T)$ and $\Bt(T)$ intersect transversely. 
\endproclaim

\demo{Proof}
Let us prove~(4.3.11). Clearly it suffices to show that 
$$\text{$0\in \Bt(T)$ if and only if $F$ is singular at $-x$.}\tag{4.3.12}$$
The ``only if'' follows immediately from Lemma~(4.3.3). Now assume $F$
is singular at $(-x)$.  From~(4.3.9) we get an exact sequence
$$0\to\xi\to F^{**}\to\xi^{-1}\to 0.\tag{4.3.13}$$
We also have an exact sequence
$$0\to F\lra F^{**}\brel(p,q)\over\lra\CC_{x}\op\CC_{-x}\to 0,$$
where 
$$p(\xi^{\pm 1}(x))\not=0, \qquad q(\xi(-x))=0.$$
(Here $\xi^{\pm 1}(\pm x)$ is the fiber of $\xi^{\pm 1}$ at $\pm x$.) Choose a
surjection 
$$F^{**}(-x)\brel\tau\over\lra\xi(-x),$$
and let $S:=\Hom(\xi(-x),\CC_{-x})$. For $s\in S$ let 
$$q(s):=q+s\circ\tau,$$
and let $F_s$ be the sheaf on $J$ fitting into the exact sequence
$$0\to F_s\lra F^{**}\brel(p,q(s))\over\lra\CC_{x}\op\CC_{-x}\to 0.$$
The sheaves $F_s$ form a family of torsion free simple semistable sheaves on
$J$, with $v(F_s)=\vv$ for all $s$. Furthermore  $F_s$ is
singular for all $s$, isomorphic to $F$ for $s=0$,  and stable for $s\not=0$. This
shows that $0\in \Bt(T)$, and finishes the proof of~(4.3.11). 
Now let's prove that $\Sit(T)\cap \Bt(T)$ intersect transversely. Clearly it
suffices to prove that they are transverse at $0$.   Let $\ov{T}$ be the
universal deformation space of $F$, and $\ov{\cE}$ the family of sheaves on $J$
parametrized by $\ov{T}$. Thus $T\ss\ov{T}$, and $\cE=\ov{\cE}|_{J\tm T}$. Set 
$$\align
\Sit_{\vv}(\ov{T}):= & \{t\in \ov{T}|\ \text{$\ov{\cE}_t$ is strictly semistable}\},\\
\Bt_{\vv}(\ov{T}):= & 
\text{closure of $\{t\in \ov{T}|\ \text{$\ov{\cE}_t$ is singular and stable}\}$.} 
\endalign$$
The map
$$\matrix
\ov{T} & \lra & J\tm\Jh\\
t & \mapsto & (\sum c_2(\ov{\cE}_t),\det(\ov{\cE}_t))
\endmatrix$$
is submersive at $0$ and gives fibrations $\Sit_{\vv}(\ov{T})\to U$
and $\Bt_{\vv}(\ov{T})\to U$,
where $U\ss J\tm\Jh$ is an open subset, with fibers $\Sit(T)$, $\Bt(T)$
respectively. Thus in order to prove that 
$\Sit(T)$ and $\Bt(T)$ iintersect transversely at $0$ it
suffices to show that $\Sit_{\vv}(\ov{T})$ and $\Bt_{\vv}(\ov{T})$ are transverse at
$0$.  By~\cite{O4,(1.17)} we have  
$$T_0\Sit_{\vv}(\ov{T})=\{\e\in \Ext^1(F,F)|\ \b\circ\e\circ\a=0\},$$
where $\a,\b$ are as in~(4.3.9), and we identify $T_0 (\ov{T})$ with 
$\Ext^1(F,F)$ via the Kodaira-Spencer map. Let $S$ be as above, and 
$\k\cl T_0 S\to\Ext^1(F,F)$ be the Kodaira-Spencer map. It is easy to check
that  
$$\b\circ\k({\partial\over\partial s}(0))\circ\a\not=0.$$
Since $S\ss \Bt(T)\ss\Bt_{\vv}(\ov{T})$, we see that $T_0\Sit_{\vv}(\ov{T})\not\ss
T_0 \Bt_{\vv}(\ov{T})$. In particular $\Sit_{\vv}(\ov{T})$ is smooth in $0$ (this also
follows from the fact that $\Sit_{\vv}(\ov{T})$ is the pull-back of $\Sit_{\vv}$ for
the map $\ov{T}\to\Mt_{\vv}$, which is \'etale in $0$.) Hence all that remains to be 
proved is that $\Bt_{\vv}(\ov{T})$ is smooth at $0$. This we do by analyzing the
map  $\rho$ fitting into the exact sequence
$$\Ext^1(F,F)\brel\rho\over\lra H^0(Ext^1(F,F))\to 
H^2(Hom(F,F))\brel\l\over\lra\Ext^2(F,F),$$
a piece of the local-to-global Exact Sequence for $\Ext^{\bu}(F,F)$. Since
$$H^0(Ext^1(F,F))\cong
Ext^1(F\ot\cO_{J,x},F\ot\cO_{J,x})\op Ext^1(F\ot\cO_{J,-x},F\ot\cO_{J,-x}),$$
we can define $\rho_{\pm x}$ as the composition of $\rho$ with   projection to the
first and second summand respectively.  
 In order to prove that $T_0\Bt_{\vv}(\ov{T})\not= T_0 (\ov{T})$ it suffices to show
that $\rho_{-x}$ is surjective. In fact assume this has been proved: by  
Lemma~(4.3.3)  first-order
deformations belonging to $T_0 \Bt_{\vv}(\ov{T})$ do not smooth the singularity
of $F$ at $(-x)$, and hence $T_0 \Bt_{\vv}(\ov{T})$ is a proper subspace of
$\Ext^1(F,F)$. Now let's show that $\rho_{-x}$ is surjective. The transpose of $\l$ is
given by  
$$\Hom(F,F)\brel\l^{*}\over\lra\Hom(Hom(F,F),\cO_J),$$
From
$$0\to Hom(F,F)\to Hom(F^{**},F^{**})\to \CC_{x}\op\CC_{-x}\to 0$$
one gets that $\Hom(Hom(F,F),\cO_J)\cong \Hom(F^{**},F^{**})$. With this
identification the map $\l^{*}$ is the canonical one. We claim that
$$\dim\Hom(F,F)=1,\qquad \dim\Hom(F^{**},F^{**})=2.\tag{4.3.14}$$
The first equation holds because $F$ is simple. To prove the second equation,
look at Exact Sequence~(4.3.13). If $\xi\not\cong\xi^{-1}$ the exact
sequence is split, and the equation follows. If $\xi\cong\xi^{-1}$, the
exact sequence is not split, otherwise we would have $F\cong (I_{x}\op
I_{-x})\ot\xi$, which is not simple: it follows easily that the second equation
holds.  From~(4.3.14) we get that $\dim\ker\l=1$, and hence
$$\cod(\im\rho)=1.\tag{4.3.15}$$
We claim $\rho_x$ is not surjective. Since  $\Ext^1(F,F)$ is spanned by 
$T_0\Sit_{\vv}(\ov{T})$ and $T_0\Bt_{\vv}(\ov{T})$ (because
$T_0\Sit_{\vv}(\ov{T})$ is of codimension $1$ in $\Ext^1(F,F)$ and not equal to
$T_0\Bt_{\vv}(\ov{T})$) it suffices to check that every first-order deformation in
$T_0\Sit_{\vv}(\ov{T})$ or in $T_0\Bt_{\vv}(\ov{T})$ does not smooth the
singularity in $x$: this follows from~(2.1.2) and~(4.3.3).   Thus we get
from~(4.3.15) that $\rho_{-x}$ is surjective. This finishes the proof of the
proposition.    
\qed
\enddemo

\demo{Proof of Proposition~(4.3.2)}
Item~(1) follows at once from Item~(3) of~(4.3.3) and Lemma~(4.3.10).
Now we prove Items~(2)-(3). Let $U\ss\Si$ be the open subset given by
$$U:= \{[I_x\ot\xi\op I_{-x}\ot\xi^{-1}]|\ x\not=-x\},$$
and let $\tU:=\pit^{-1}U$. 
By Proposition~(2.3.7) and Lemma~(4.3.10)   $\Sit\cap\Bt\cap\tU$
is smooth.  Let $J^0:=(J\sm J[2])$, where $J[2]$ is the kernel of  
multiplication by $2$, and set 
$$\tX:=(J^0\tm\Jh)\tm_U(\Sit\cap\Bt\cap\tU),\tag{4.3.16}$$
where $(J^0\tm\Jh)\to U$ is the \'etale degree-two map
$$\matrix
J^0 \tm\Jh & \lra & U\\
(x,\xi) & \mapsto & [I_x\ot\xi\op I_{-x}\ot\xi^{-1}].
\endmatrix$$
 Since $\Sit\cap\Bt\cap\tU$ is smooth and the above map is \'etale, $\tX$ is
smooth.  

\proclaim{(4.3.17)Claim}
Keep notation as above. The projection $h\cl \tX\to(J^0\tm\Jh)$ is the
blow-up of $(J^0\tm\Jh[2])$.
\endproclaim 

\demo{Proof of the claim}
Since $\tX$ is mooth it suffices to show that 
$$h^{-1}(x,\xi)\cong
\cases
\text{one point} & \text{if $\xi\notin \Jh[2]$,}\\
\PP^1 & \text{if $\xi\in \Jh[2]$.}
\endcases
\tag{4.3.18}$$
The exact sequence
$$0\to I_{-x}\ot\xi^{-1}\to\xi^{-1}\to\xi^{-1}|_{-x}\to 0$$
gives an exact sequence
$$0\to \Ext^1(\xi^{-1},I_{x}\ot\xi)\lra
\Ext^1(I_{-x}\ot \xi^{-1},I_{x}\ot\xi),$$
and
$$h^{-1}(x,\xi)= 
\PP(\Ext^1(\xi^{-1},I_{x}\ot\xi))=\PP(H^1(I_{x}\ot\xi^{\ot 2})).
\tag{4.3.19}$$ 
An easy computation gives~(4.3.18), and this proves~(4.3.17). 
\qed
\enddemo

Let $g\cl \tX\to J^0$ be the projection, and let  
$$\tX_{\a}:= g^{-1}(\Ta)\cong C\tm \wh{I}_0,$$
where the second isomorphism  holds because of~(4.3.17).  By~(4.3.16) we
have a natural map 
$$\tX_{\a}\to\Sta\cap\Bta\cap\tU=\Sta\cap\Bta,$$
where the equality follows from Item~(3) of~(4.3.3).  One verifies easily 
that the map above is the quotient map for $\approx$;
proving Item~(2). Item~(3) is immediate  
\qed
\enddemo

\subhead
4.4. Topological results
\endsubhead
First notice that by~(4.1.1) and~(4.2.3) 
$$\text{$\Sta$ is irreducible.}\tag"$\text{(4.4.1)}^{*}$"$$
Furthermore~(4.1.2) and Item~(2) of~(4.2.5) give
$$b^1(\Sta)=5.\tag"$\text{(4.4.2)}^{*}$"$$

Now we pass to $\Sta\cap\Bta$. Given Item~(2) of~(4.3.2), the
low-dimensional cohomology of $\Sta\cap\Bta$ has a description similar to
that of $\Sia$ given in~(4.2.5). More precisely, choose $z\in \wh{I}_0$
belonging to the exceptional divisor of $\k$. The map
$$\matrix
C & \lra & C\tm \wh{I}_0\\
q & \mapsto & (q,z)
\endmatrix$$
induces a map $\wt{\mu}\cl\ov{C}\to\Sta\cap\Bta$. ``Copying'' the proof of
Item~(2) of~(4.2.5) one gets that
$\wt{\mu}^{*}\cl H^1(\Sta\cap\Bta;\QQ)\to H^1(\ov{C};\QQ)$ is
an isomorphism, hence
$$b^1(\Sta\cap\Bta)=5.\tag"$\text{(4.4.3)}^{*}$"$$
Furthermore, letting ${\tilde\k}$ be as in Item~(2) of~(4.3.2), one gets
that 
$$\text{${\tilde\k}^{*}\cl H^2(\Sta\cap\Bta;\QQ)\to H^2(C\tm
\wh{I}_0;\QQ)$ is an isomorphism}\tag{4.4.4}$$
by ``copying'' the proof of Item~(3) of~(4.2.5). Thus ${\tilde\k}^{*}$ gives
an identification 
$$H^2(\Sta\cap\Bta;\QQ)\cong\bigoplus\limits_{i=1}^{16}\QQ[C\tm E_i]
\bigoplus(\id_C\tm\nu_0)^{*}H^2(C\tm\Jh),\tag{4.4.5}$$
where $E_i\ss \wh{I}_0$ are the irreducible components of the exceptional
divisor of $\nu_0$.  

Next we examine the maps on fundamental groups and cohomology induced
by the inclusion $\rho\cl\Sta\cap\Bta\hra\Sta$. 

\proclaim{$\text{(4.4.6)}^{*}$Proposition}
Keeping notation as above 
$\rho_{\#}\cl\pi_1(\Sta\cap\Bta)\to\pi_1(\Sta)$   is an isomorphism. 
\endproclaim

\demo{Proof}
The maps
$$\pit_{\#}\cl\pi_1(\Sta)\to\pi_1(\Sia)
\qquad\pit_{\#}\circ\rho_{\#}\cl\pi_1(\Sta\cap\Bta)\to\pi_1(\Sia)$$
are isomorphisms, and this implies that $\rho_{\#}$ is an isomorphism.
\qed
\enddemo

\proclaim{(4.4.7)Proposition}
The map 
$$H^2(\Sta;\QQ)\overset\rho^{*}_2\to{\lra} H^2(\Sta\cap\Bta;\QQ)$$
induced by Inclusion $\rho$ is injective. Furthermore, referring to~(4.4.5), we
have 
$$\Im({\tilde\k}^{*}\rho^{*}_2)=\QQ[C\tm E]
\bigoplus(\id_C\tm\nu_0)^{*}H^2(C\tm\Jh).\tag"$\text{(4.4.6)}^{*}$"$$
\endproclaim

\demo{Proof}
Given that  
$$(\id_C\tm\k)=f|_{\Sta\cap\Bta},$$
where $f$ is the map of~(4.1.1), and given~(4.1.3) together with Item~(3)
of~(4.2.5), it suffices to notice that 
$$\langle{\tilde\k}^{*}\rho^{*}c_1(\om_f),p\tm E_i\rangle=-2$$
for every $i=1,\ldots,16$ (here $p\in C$). 
\qed
\enddemo

We will need to know exactly the class of
${\tilde\k}^{*}\rho^{*}(c_1(\om_f))$. Let
$\cL$ be the Poincar\'e line-bundle on $J\tm\Jh$, and $\ia$ as in (1.8).
We will prove that 
$${\tilde\k}^{*}\rho^{*}(\om_f)\cong 
(\ia\tm\nu_0)^{*}(\cL^{\ot 4})\ot [2C\tm E],
\tag"$\text{(4.4.9)}^{*}$"$$ 
where $E\ss\wh{I}_0$ is the exceptional divisor of $\nu_0$. 

\demo{Proof of~(4.4.9)}
Let $\Jh^0:=(\Jh\sm\Jh[2])$, and $\ih\cl\Jh^0\hra\Jh$ be the inclusion.
Notice that we also have a natural inclusion $\Jh^0\hra\wh{I}_0$. It suffices
to prove that  
$${\tilde\k}^{*}\rho^{*}(\om_f)|_{C\tm\Jh^0}\cong 
(\ia\tm\ih)^{*}\cL^{\ot 4}.\tag{4.4.10}$$
 In fact, since 
$${\tilde\k}^{*}\rho^{*}(\om_f)|_{p\tm E_i}\cong\cO_{E_i}(-2)$$
for $p\in C$ and $i=1,\ldots,16$, Equality~(4.4.10) immediately
implies~(4.4.9).  To prove~(4.4.10) we must describe explicitely  the inclusion
$$C\tm\Jh^0\hra (C\tm\Jh^0)\tm_{\Sia}\Sta\tag{4.4.11}$$
whose image is dense in $\Sta\cap\Bta$. Let $\D_{\a},-\D_{\a}\ss
J\tm(C\tm\Jh^0)$ be given by 
$$\align
\D_{\a} & :=\{(\ia(p),p,\xi)|\ p\in C,\quad \xi\in\Jh^0\},\\
-\D_{\a} & :=\{(-\ia(p),p,\xi)|\ p\in C,\quad \xi\in\Jh^0\},
\endalign$$
and let $\vf\cl J\tm(C\tm\Jh^0)\to(C\tm\Jh^0)$ be the projection. 
By~(4.3.19) we can identify~(4.4.11) with 
$$\PP R^1\vf_{*}(I_{\D_{\a}}\ot(\cL^0)^{\ot 2})\hra
\PP Ext^1_{\vf}(I_{-\D_{\a}}\ot\cL^0,
I_{\D_{\a}}\ot(\cL^0)^{-1})),\tag{4.4.12}$$
where $\cL^0:=(\id_J\tm\ih)^{*}\cL$. (To simplify notation we omit
pull-back signs whenever possible.) Thus  
$${\tilde\k}^{*}\rho^{*}(\om_f)|_{C\tm\Jh^0}\cong
R^1\vf_{*}(I_{\D_{\a}}\ot(\cL^0)^{\ot 2})\ot 
\cQ^{-1},\tag{4.4.13}$$
where $\cQ$ is the line-bundle fitting into the exact sequence
$$0\to R^1\vf_{*}(\bu)\to
Ext^1_{\vf}(\bu,\bu)\to \cQ\to 0.$$
(Here $R^1$ and  $Ext^1$ are the sheaves appearing in~(4.4.12).) 
We have an isomorphism
$$\cQ\cong Ext^2_{\vf}((\cL^0)^{-1}|_{-\D_{\a}},
I_{\D_{\a}}\ot\cL^0),$$
and applying Serre duality we get
$$\cQ\cong Hom_{\vf}(\cL^0,
(\cL^0)^{-1}|_{-\D_{\a}})^{-1}\cong (\ia\tm\ih)^{*}\cL^{-2}.$$
Another easy computation gives
$$R^1\pi_{*}(I_{\D_{\a}}\ot(\cL^0)^{\ot 2})\cong
(\ia\tm\ih)^{*}\cL^{2}.$$
Using~(4.4.13) we get Equation~(4.4.10).
\qed
\enddemo
\head
5. Analysis of $\Bta$.
\endhead
There is a single $*$-red tag not contained in Subsection~(5.2), i.e.~(5.1.1).
\subhead
5.1. Realization of $\Bta$ as a $\PP^1$-fibration
\endsubhead
We will prove that there is a fibration
$$\CD
\PP^1 @>>> \Bta \\
@.                  @VV g V \\
{}         @.       \ov{C}\tm K^{[2]}\Jh,
\endCD
\tag"$\text{(5.1.1)}^{*}$"$$
where $K^{[2]}\Jh$ is the Kummer surface of $\Jh$, i.e.~the subset of $\Jh^{[2]}$
parametrizing subschemes whose associated cycle sums up to $\wh{0}$. We define
$g$ by giving its two components. First we define $g_1\cl\Bta\to\ov{C}$. Let
$z\in\Bta$. By Item~(3) of~(4.3.3) the sheaf $F_z$ has two distinct
singularities, say $x$ and $-x$. One at least of $x$, $-x$ belongs to $\Ta$.
The set of unordered couples $(x,-x)$ intersecting $\Ta$ is identified with
$\ov{C}$ via the map   
$$\matrix
C & \lra & J\tm J\\
p & \mapsto & (i_{\a}(p),-i_{\a}(p)).
\endmatrix$$
Given the above identification we set
$$\matrix
\Bta & \brel g_1\over\lra & \ov{C}\\
z & \mapsto & sing(F_z).
\endmatrix
\tag{5.1.2}$$
In order to define $g_2\cl\Bta\to K^{[2]}\Jh$ we must analyze the moduli space of
isomorphism classes of semistable rank-two vector-bundles $E$ on $J$ with
$c_1^{hom}=0=c_2^{hom}$ (i.e.~flat vector-bundles), such that $\dim\Hom(E,E)=2$
(the minimal dimension for such bundles).  We claim that this moduli space can
be identified with $\Jh^{[2]}$. This follows from a more general statement proved 
in~\cite{BDL}. We sketch a proof for the reader's convenience. Let $\Mf(2)$
be the set of isomorphism classes of such vector-bundles. First we give a
one-to-one correspondence 
$$\Mf(2)\brel\vf\over\longleftrightarrow\Jh^{[2]}.$$
Let $[E]\in\Mf(2)$; then either 
$E\cong\xi_1\op\xi_2$, where $\xi_1\not\cong\xi_2$, or $E$ fits into a non-split
exact sequence
$$0\to\xi\to E\to\xi\to 0.\tag{5.1.3}$$
(We are excluding the flat bundles isomorphic to $\xi\op\xi$.) In the former
case we set $\vf([E]):=\{\xi_1,\xi_2\}$, in the latter  we associate to $[E]$ the
couple $([\xi],\CC e)$, where $e\in H^1(\cO_J)\cong T_{[\xi]}\Jh$ is the extension
class of~(5.1.3).  (Thus in the latter case $\vf([E])$ is a non-reduced
subscheme of $\Jh$.) Clearly $\vf$ is a bijection.

\proclaim{(5.1.4)Proposition}
Keep notation as above. Let $\cE$ be a family of vector-bundles on $J$
parametrized by (a reduced scheme) $T$, with $[\cE_t]\in\Mf(2)$ for all $t\in T$.
The map $\vf_{\cE}\cl T\to \Mf(2)$ is regular.
\endproclaim

\demo{Proof}
It is sufficient to prove the following. Assume that $0\in T$, and that $T$ is the
deformation space of $E:=\cE_0$: then 
$\vf_{\cE}$ is regular in a neighborhood of $0$. If $E\cong\xi_1\op\xi_2$ the
statement is obvious, so let's assume $E$ fits into Exact Sequence~(5.1.3),
with extension class $e$. Let 
$$\L:=\{t\in T|\ \text{$\cE_t$ is not a direct sum}\}.$$
Then $\L$ is of pure codimension $1$. Let $f\cl\wt{T}\to T$ be the relative
Quot-scheme parametrizing quotients $\cE_t\to \eta$, where $\eta\in\Jh$. Then
$\wt{T}$ is smooth, and $f$ is a double covering, ramified over $\L$. Let
$\wt{\L}:=f^{-1}(\L)$. For $\wt{t}\in\wt{T}$ with $f(\wt{t})=t$, we have a
well-defined exact sequence
$$0\to\l\to\cE_t\to\eta\to 0,$$
where $\l,\eta\in\Jh$. Thus we can define a (regular) map
$$\matrix
\wt{T} & \brel g\over\lra & \Jh\tm\Jh\\
\wt{t} & \mapsto & (\l,\eta).
\endmatrix$$
Clearly $g(\wt{\L})=\wh{\D}$, where $\wh{\D}$ is the diagonal. Thus $dg(\wt{0})$
(here $\wt{0}$ is the point such that $f(\wt{0})=0$) induces a map
$$\ov{dg}(\wt{0})\cl N_{\wt{T},\wt{\L}}(\wt{0})\to 
T_{(\xi,\xi)}(\Jh\tm\Jh)\cong H^1(\cO_J).$$
We claim that 
$$\Im\ov{dg}(\wt{0})=\CC e.\tag{5.1.5}$$
Given the above formula one concludes that $\vf_{\cE}(0)$ is regular in a
neighborhood of $0$ as follows. By~(5.1.5) the pull-back
$g^{*}(I_{\wh{\D}})$ is a locally principal ideal sheaf with zero-locus
$\wt{\L}$, hence $g$ lifts to a regular map 
$$G\cl\wt{T}\to Bl_{\wh{\D}}(\Jh\tm\Jh).$$
Since $\Jh^{[2]}$ is the quotient of $Bl_{\wh{\D}}(\Jh\tm\Jh)$ by the natural
involution, we get a map $H\cl\wt{T}\to\Jh^{[2]}$. This map is invariant by the
involution interchanging the sheets of $f$, hence it descends to a regular map
$h\cl T\to\Jh^{[2]}$. Obviously $h$ coincides with $\vf_{\cE}$ on $T\sm\L$.
Furthermore by~(5.1.5) we have $h(0)=\vf_{\cE}(0)$; since the germ $(T,t)$
is the deformation space of $\cE_t$ for all $t$ in a neighborhood of $0$, we
get that $h(t)=\vf_{\cE}(t)$ for all $t$ near $0$, and proves the
proposition. We are left with the task of proving~(5.1.5). Let $\G\ss\Jh$ be a
curve containing $\xi$, smooth in $\xi$, and such that $T_{\xi}\G=\CC e$. Let
$i\cl J\tm\{0\}\hra J\tm\G$ be the inclusion. Let $\cF$ be the sheaf on 
$J\tm\G$ fitting into the exact sequence 
$$0\to\cF\to\xi\ot\cO_{\G}\op\cL|_{J\tm\G}\brel\a\over\to i_{*}\xi\to 0,$$ 
where $\cL$ is the normalized Poincar\'e line-bundle, and we choose $\a$ so that  
neither $\xi\ot\cO_{\G}$ nor $\cL|_{J\tm\G}$ belongs to $\cF$. Then $\cF$ is a family
of sheaves on $J$ parametrized by $\G$, and $[\cF_{\wh{x}}]\in\Mf(2)$ for all
$\wh{x}\in\G$. Furthermore we have an exact sequence
$$0\to\xi\to\cF_0\to \xi\to 0$$
and one verifies that the extension class is equal to $e$. Thus $\cF_0\cong E$,
hence $\cF$ induces a regular map $m\cl \G\to T$, with $m(0)=0$. Actually $m$
lifts to a map $\wt{m}\cl\G\to\wt{T}$, because the inclusion
$\xi\ot\cO_{\G}(-J\tm\{0\})\hra \cF$ determines a quotient $\cF_{\wh{x}}\to
\eta$, with $\eta\in\Jh$, for all $\wh{x}\in\G$. As is easily verified 
$$\ov{dg}(\wt{0})\circ d\wt{m}(0)=\CC e.$$
This proves~(5.1.5), because $\Im\ov{dg}(\wt{0})$ can not be equal to all of
$H^1(\cO_J)$.   
\qed
\enddemo

Now we are ready to define $g_2\cl\Bta\to K^{[2]}\Jh$. We work
locally. Let $z\in\Bta$; by Proposition~(2.3.7) there exists a neighborhood
$U\ss\Bta$ (in the classical topology) of $z$ parametrizing a family $\cF$ of
torsion-free simple semistable sheaves on $J$, with $v(\cF_t)=\vv$ for all $t\in
U$, such that the $\wt{S}$-equivalence class of $\cF_t$ corresponds to $t$
via~(2.2.3). By Proposition~(4.3.3) we have
 $[\cF_t^{**}]\in\Mf(2)$ for all $t\in U$. Thus by Proposition~(5.1.4) the
family $\cF^{**}$ induces a regular map $h\cl U\to \Jh^{[2]}$. Since
$\det(\cF_t)\cong\cO_J$ for all $t\in U$, the image of $h$ is contained in
$K^{[2]}\Jh$. These maps are independent of the choice of family $\cF$ (see
Remark~(2.3.9)) and they glue together, defining $g_2\cl\Bta\to
K^{[2]}\Jh$. 

\proclaim{Claim}
The map $g:=(g_1,g_2)\cl\Bta\to\ov{C}\tm K^{[2]}\Jh$ is a $\PP^1$-fibration.
\endproclaim

\demo{Proof}
Let $((x,-x),[E])\in \ov{C}\tm K^{[2]}\Jh$. The fiber $g^{-1}((x,-x),[E])$ is
isomorphic to the set of $\wt{S}$-equivalence classes of simple semistable
sheaves $F$ fitting into an exact sequence
$$0\to F\to E\to \CC_{x}\op\CC_{-x}\to 0.$$
One verifies easily that this set is isomorphic to $\PP^1$. We leave it to the reader
to check that $g$ is locally trivial.  
\qed
\enddemo   
\subhead
5.2. Topological results
\endsubhead
From Fibration~(5.1.1) we get that
$$\text{$\Bta$ is irreducible,}\tag"$\text{(5.2.1)}^{*}$"$$
and that
$$\align
H^1(\Bta) & \cong H^1(\ov{C})\cong\QQ^5, \tag"$\text{(5.2.2)}^{*}$"\\ 
H^2(\Bta) &\cong H^2(\ov{C})\op H^2(K^{[2]}\Jh)\op\QQ c_1(\om_g),
\tag{5.2.3}
\endalign$$
where $\om_g$ is the relative cotangent bundle of $g$. 
 
In proving that $\Mt$ is simply connected it will be useful to represent
$\pi_1(\Bta)$ as follows. Let   $E_1\ss \wh{I}_0$ be the component of the
exceptional divisor of $\nu_{0}$  mapping to $\hat 0$ and let 
$$R:=(C\tm E_1/\approx)\cong\ov{C}\tm E_1,\tag{5.2.4}$$
where $\approx$ is as in~(4.3.1), and $j\cl R\hra\Bta$ be the inclusion
given by~(4.3.2). Then
$$\text{$j_{\#}\cl\pi_1(R)\to\pi_1(\Bta)$ is an
isomorphism}\tag"$\text{(5.2.5)}^{*}$"$$ 
because, letting $g_1$ be as in~(5.1.2), the composition $(g_1\circ j)\cl
R\to\ov{C}$ is the projection to  the first factor of~(5.2.4), hence 
$g_{1,\#}\circ j_{\#}$ is an isomorphism, and 
$g_{1,\#}\cl\pi_1(\Bta)\to\pi_1(\ov{C})$ is also an isomorphism.     

For the proof of the main Theorem we will need to know the map on
$2$-cohomology induced by the inclusion $\l\cl(\Sta\cap\Bta)\hra\Bta$.  

\proclaim{$\text{(5.2.6)}^{*}$Proposition}
Keep notation as above. Then 
$$\l^{*}\cl H^2(\Bta)\to H^2(\Bta\cap\Sta)$$ 
is  injective.
Furthermore, given Isomorphism~(4.4.5), 
$$\Im({\tilde\k}^{*}\l^{*}_2)=\bigoplus\limits_{i=1}^{16}\QQ[C\tm E_i]
\bigoplus(\id_C\tm\nu_0)^{*}(H^2(C)\op H^2(\Jh)
\op\QQ c_1((\ia\tm\id_{\Jh})^{*}\cL)),\tag"$\text{(5.2.7)}^{*}$"$$
where ${\tilde\k}$ is as in Item~(2) of~(4.3.2). 
\endproclaim

\demo{Proof}
Given~(5.2.3) and~(4.4.5), in order to prove that $\l^{*}$ is injective and 
that (5.2.7) holds, it suffices to show that
$${\tilde\k}^{*}\l^{*}\om_g=
(\id_C\tm\nu_0)^{*}\cL^{\ot(-4)}\ot[-2C\tm E].\tag{5.2.8}$$ 
(Recall~(4.3.2).) We will show that
$${\tilde\k}^{*}\l^{*}(c_1(\om_g))=-{\tilde\k}^{*}\rho^{*}(c_1(\om_f))
\tag{5.2.9}$$
in the Chow group, and then~(5.2.8) will follow from~(4.4.9). 
Let $\Sit^{*}\ss\Sit$ and $\Bt^{*}\ss\Bt$ be the open subsets given  by
$$\align
\Sit^{*}:= & \pit^{-1}\{[I_x\ot\xi\op I_{-x}\ot\xi^{-1}]|
\ x\not=(-x),\quad \xi\not\cong\xi^{-1}\},\\
\Bt^{*}:= & \{[F]|\ sing(F)=\{x,-x\},\ x\not=(-x),
\quad F^{**}\cong\xi\op\xi^{_1},\ \xi\not\cong\xi^{-1}\}.
\endalign$$
By Lemma~(4.3.10) they intersect transversely, and hence by adjunction
$$K_{\Sit^{*}\cap\Bt^{*}}\sim[\Sit+\Bt]|_{\Sit^{*}\cap\Bt^{*}}.$$
(Because $K_{\Mt}\sim 0$.) Restricting $\pit$ to $\Sit^{*}\cap\Bt^{*}$ we get an
isomorphism $(\Sit^{*}\cap\Bt^{*})\cong (J^0\tm\Jh^0)$, hence 
$K_{\Sit^{*}\cap\Bt^{*}}\sim 0$. Thus
$$[\Sit]|_{\Sit^{*}\cap\Bt^{*}}\sim -[\Bt]|_{\Sit^{*}\cap\Bt^{*}}.
\tag{5.2.10}$$
Restricting $\pit$ to $\Sit^{*}$ we get a $\PP^1$-fibration 
$F\cl \Sit^{*}\to(J^0\tm\Jh^0)$, and thus
$$[\Sit]|_{\Sit^{*}}\sim K_{\Sit^{*}}\sim\om_F.$$
In fact the first equivalence follows from adjunction, and the second
one (where $\om_F$ is the relative cotangent bundle of $F$) from
triviality of the canonical bundle of $(J^0\tm\Jh^0)$. We also have a
$\PP^1$-fibration  
$$\matrix
\Bt^{*} & \brel G\over\lra & (J^0/(-1))\tm(\Jh^0/(-1))\\
[F] & \mapsto & (sing(F),\{\xi,\xi^{-1}\}),
\endmatrix$$
where $F^{**}\cong\xi\op\xi^{-1}$. Thus by a similar argument we get
$$[\Bt]|_{\Bt^{*}}\sim K_{\Bt^{*}}\sim\om_G.$$
The last two equivalences together with~(5.2.10) prove~(5.2.9).
\qed
\enddemo
\head
6. Analysis of $\Vta$.
\endhead
\subhead
6.1. Introduction
\endsubhead
By definition $\Va^0$ is the locus parametrizing sheaves $F$ such
that there is an exact sequence
$$0\to\l\lra\isa^{*}F\brel f\over\lra\xi\to 0\tag{6.1.1}$$
with $\isa^{*}F$ locally-free and $\xi$ a line-bundle of negative degree.   Let $G$
be the torsion-free sheaf on $J$ fitting into the exact sequence
$$0\to G\lra F(\Ta)\brel \wt{f}\over\lra i_{\a,{*}}(\xi\ot K_C)\to 0,
\tag{6.1.2}$$
where $\tilde f$ is obtained tensorizing $f$ with the identity map of $K_C$ (notice
that by adjunction $\isa^{*}(\Ta)\sim K_C$). In other words $G$ is the elementary
modification of $F(\Ta)$ associated to $\tilde f$. We will show
(see Lemma~(6.3.22)) that  $\deg\xi=-1$, and thus by a straightforward 
computation   
$$v(G)=2+\t=:\ww,$$
where  $\t:=c_1^{hom}(\T)$. We will also prove that $G$ is
slope-stable~(see Lemma~(6.3.19)), hence $[G]\in\Mw$. Since~(6.1.1) is
uniquely determined by $F$, the above construction defines a regular map
from $\Va^0$ to $\Mw$.  Notice that $G$ comes with extra structure,
namely an exact sequence  
$$0\to\xi\lra\isa^{*}G\brel\psi\over\lra\l\ot K_C\to 0\tag{6.1.3}$$
which is part of the long exact sequence of Tor's obtained by applying the functor
$\ot\cO_{\Ta}$ to~(6.1.2). As is easily verified, the sheaf $F$ is the
elementary modification of $G$ associated to Exact Sequence~(6.1.3). Hence
by associating to $[F]$ the couple consisting of $[G]$ and Exact
Sequence~(6.1.3) we get an isomorphism between  $\Va^0$ and a
locally-closed subset of a quot-scheme over $\Mw$. Since $\Va^0$ is dense
in $\Va$, this will give a way of relating $\Va$ to $\Mw$. The section is
organized as follows. In~(6.2) following Mukai, Yoshioka~\cite{Muk1,Y2} we
apply the Mukai transform to prove that $\Mw$ is the product $J\tm \Jh$.
In~(6.3) we show that $\Va^0$ is isomorphic to $\Pa^0$, a locally closed
subset of a Quot-scheme over $\Mw$. We denote by $\Pa$ the closure of
$\Pa^0$ in the Quot-scheme over $\Mw$; in Subsections(6.4)-(6.5) we give
results on $\Pa$ and the boundary $\Pa\sm\Pa^0$. We will not construct
directly $\Pa$ (the difficulty lies in analyzing it singular points), instead we
construct a $\PP^1$-bundle over the blow-up of $\Jh$ at $16$ points (we
denote it by $\Pam$) and then prove that it is a desingularization of $\Pa$:
this is done in~(6.6). The $3$-fold $\Pam$ has $16$ $(-1,-1)$ curves: we let
$\Pap$ be the manifold obtained by flopping all of these curves. In~(6.7) we
we show that $\Vta$ is isomorphic to $\Pap$. In the last subsection we
give the results on the topology of $\Vta$ and its intersections with $\Sta$,
$\Bta$ which are needed in the proof of Theorem~(1.4).
\subhead
6.2. The moduli space $\cM_{\ww}$
\endsubhead
For $(x,\yh)\in J\tm\Jh$, let $E(x,\yh)$ be the sheaf on $J$ defined by 
$$E(x,\yh):=\phi_{*}(\cL\ot\phih^{*}(i_{\yh,*}[3p_0+x])).$$
(Here $\phi,\phih$ are as in~(1.9).) 
We will prove the following result.

\proclaim{(6.2.1)Proposition}
Keep notation as above. 
\roster
\item
$E(x,\yh)^{*}$ is a vector-bundle on $J$, and $v(E(x,\yh)^{*})=2+\t$.
Furthermore
$$\align
\det(E(x,\yh)^{*})& \cong[\T]\ot\cL_{-\wh{x}-2\yh}\tag{6.2.2}\\
\sum c_2(E(x,\yh)^{*}) & = -x-y.\tag{6.2.3}
\endalign$$
\item
$E(x,\yh)^{*}$ is slope-stable for every $(x,\yh)$, and the map
$$\matrix
J\tm\Jh & \brel\rho\over\lra & \cM_{\ww} \\
(x,\yh) & \mapsto & [E(x,\yh)^{*}]
\endmatrix$$
is an isomorphism. 
\endroster
\endproclaim

\proclaim{(6.2.4)Remark}
{\rm Item~(2) of Proposition~(6.2.1) is part of a more general result
of Yoshioka~\cite{Y2} on moduli spaces of sheaves on an abelian surface 
in the case when semistability implies stability.
However, since we will need various results about the moduli space which are
not found elsewhere, e.g.~(6.2.5), (6.2.11), and which may be used to
prove~(6.2.1), we have added a proof of~(6.2.1).}  
\endproclaim

\subsubhead
Proof of Item~(1) of~(6.2.1)
\endsubsubhead
Since $h^0(C,[3p_0+x+t])=2$ for all $t\in J$, we
see~\cite{Mum,pp.46-55} that $E(x,\yh)$ is locally-free of rank two. Since the 
higher direct images vanish, the Chern classes of $E(x,\yh)$ are computed by the
Grothendieck-Riemann-Roch formula: one gets $v(E(x,\yh))=2-\t$. Thus
$E(x,\yh)^{*}$ is a rank-two vector-bundle, and its Mukai vector is as stated. 
Before going on with the proof of~(6.2.1) we determine the restriction of
$E(x,\yh)^{*}$ to $\Ta$. 

\proclaim{(6.2.5)Proposition}
Keep notation as above. There is an exact sequence
$$0\to [p_0-\a-x-y]\to
\isa^{*}E(x,\yh)^{*}\to [p_0-y]\to 0.\tag{6.2.6}$$
The extension class is the (unique up to scalars) non-trivial
one if $\a+x\not=0$, and is trivial if $\a+x=0$.
\endproclaim

\demo{Proof}
Let $\D\ss C\tm C$ be the diagonal, and $\pi_i\cl C\tm C\to C$ be the first and
second projection respectively, for $i=1,2$. Then 
$$(\isa\tm\isyh)^{*}\cL\cong
[\D]\ot\pi_1^{*}[-p_0+y]\ot\pi_2^{*}[-p_0+\a].\tag{6.2.7}$$
Thus
$$\isa^{*}E(x,\yh)=[-p_0+y]\ot
\pi_{1,{*}}([\D]\ot\pi_2^{*}[2p_0+\a+x]).\tag{6.2.8}$$
If $\a+x= 0$ one gets immediately that  
$$\isa^{*}E(x,\yh)^{*}\cong \cO_{C}(p_0-y)^{(2)}.$$
If $\a+x\not= 0$, we apply the $\pi_{1,{*}}$-functor to the exact sequence 
$$0\to\pi_2^{*}[2p_0+\a+x]\to[\D]\ot\pi_2^{*}[2p_0+\a+x]\to
([\D]\ot\pi_2^{*}[2p_0+\a+x])|_{\D}\to 0,$$
and we get an exact sequence
$$0\to\cO_C\to
\pi_{1,{*}}([\D]\ot\pi_2^{*}[2p_0+\a+x])\to[\a +x]\to 0.\tag{6.2.9}$$
The extension class is non-trivial, because otherwise we would have a
non-zero section of
$$[\D]\ot\pi_2^{*}[2p_0+\a+x] \ot\pi_1^{*}[-\a-x],$$
which is absurd since $\cO_C(\a+x)$ is non-trivial. 
Tensorizing~(6.2.9) by $[-p_0+y]$ and taking the dual exact sequence we get
the proposition.
\qed
\enddemo
Let us prove~(6.2.2)-(6.2.3). The first formula follows immediately
from~(6.2.6). In order to prove~(6.2.3), choose a canonical divisor 
$D\in |K_C|$, and consider the exact sequence of sheaves on $J\tm \Jh$:
$$0\to \cL\ot \phih^{*}(i_{\yh,*}[p_0+x]) \to \cL\ot \phih^{*}(i_{\yh,*}[p_0+x+D])
\to\cL\ot \phih^{*}(i_{\yh,*}(\cO_{D}))\to 0.$$
Applying the
$\phi_{*}$-functor and recalling that $D\sim 2p_0$, we get an exact sequence of
sheaves on $J$:   
$$0\to E(x,\yh)\to\cL_{\wh{D}}\to 
R^1\phi_{*}(\cL\ot\phih^{*}(i_{\yh,*}[p_0+x]))\to 0,\tag{6.2.10}$$
where $\wh{D}:=i_{\yh,*}(D)$. 

\proclaim{Lemma}
Keep notation as above. The sheaf 
$R^1\phi_{*}(\cL\ot \phih^{*}(i_{\yh,*}[p_0+x]))$ is supported on $\T_{-x}$.
More precisely we have
$$R^1\phi_{*}(\cL\ot \phih^{*}(i_{\yh,*}[p_0+x]))=
i_{-x,*}[p_0+y].\tag{6.2.11}$$
\endproclaim 

\demo{Proof}
By Formula~(6.2.2) and Exact Sequence~(6.2.10) we have (in the Chow ring)
$$-c_1(\T)+c_1(\cL_{\xh+2\yh})=c_1(E(x,\yh))=
c_1(\cL_{\wh{D}})-c_1(R^1\phi_{*}(\cL\ot \phih^{*}(i_{\yh,*}[p_0+x]))).$$
This gives
$$c_1(R^1\phi_{*}(\cL\ot \phih^{*}(i_{\yh,*}[p_0+x])))=c_1(\T_{-x}).$$
One concludes easily that $R^1\phi_{*}(\cL\ot \phih^{*}(i_{\yh,*}[p_0+x]))$ is the
push-forward by $i_{-x}$ of a line-bundle on $C$. 
To compute this line-bundle, pull-back by $(i_{-x}\tm i_{\yh})$ and   apply
Formula~(6.2.7) to get  
$$R^1\phi_{*}(\cL\ot \phih^{*}(i_{\yh,*}[p_0+x]))=
i_{{-x},*}([-p_0+y]\ot R^1\pi_{1,*}([\D])).\tag{6.2.12}$$
Applying the functor $\pi_{1,*}$ to the exact sequence 
$$0\to \cO_{C\tm C}\to [\D]\to[\D]|_{\D}\to 0,$$
we get
$$0\to K_C^{-1}\to H^1(\cO_C)\ot\cO_C\to 
R^1\pi_{1,*}([\D])\to 0.$$
Thus $R^1\pi_{1,*}([\D])\cong K_C$. Plugging this
into~(6.2.12) we get~(6.2.11).
\qed
\enddemo  
Let us prove~(6.2.3). It follows from~(6.2.10) and~(6.2.11) that 
$$c_2(E(x,\yh))=c_2(\cL_{\wh{D}})-i_{-x,*}c_1([-p_0+y]).$$
Since $\sum c_2(\cL_{\wh{D}})=0$, Formula~(6.2.3) follows immediately.  
\subsubhead
Proof of Item~(2) of~(6.2.1)
\endsubsubhead
By~(6.2.5) the restriction of $E(x,\yh)^{*}$ to $\Ta$ is
semistable; since $\Ta$ is numerically equivalent to $\T$ (the ample divisor
defining stability) this implies that $E(x,\yh)^{*}$
is slope-semistable. That $E(x,\yh)^{*}$ is slope-stable is an immediate
consequence of the following result.

\proclaim{(6.2.13)Lemma}
Let $F$ be a torsion-free slope-semistable sheaf on $J$ with $v(F)=\ww$. Then $F$
is slope-stable. In particular every sheaf parametrized by $\Mw$ is slope-stable.
\endproclaim

\demo{Proof}
Assume that $F$ is strictly slope-semistable.
Then there is an exact sequence
$$0\to I_W(\xi)\to F\to I_Z(\Tt-\xi)\to 0,\tag{6.2.14}$$ 
for some $\tau\in J$, where $W$, $Z$ are zero-dimensional subschemes of $J$, and
$\xi$ is a divisor such that 
$$(2\xi-\Tt)\cdot \T=0,\quad (2\xi-\Tt)^2\le 0.\tag{6.2.15}$$
(The inequality follows from the equality and Hodge index.)
A simple manipulation then gives that $\xi\cdot\xi\le 0$.
On the other hand, computing $c_2(F)$ using Exact Sequence~(6.2.14)
one gets that
$$\xi\cdot\xi=\ell(W)+\ell(Z)\ge 0.$$
Thus $\xi\cdot\xi=0$, hence~(6.2.15) gives
$$(2\xi-\Tt)^2=-2.$$
This together with the equality of~(6.2.15) contradicts Assumption~(1.3).
Hence $F$ must be slope-stable. 
\qed
\enddemo

Now we prove that $\rho$ is an isomorphism.  By~(6.2.13) all sheaves
parametrized by $\Mw$ are stable, hence by Mukai~\cite{Muk2,(0.1)} $\Mw$
is  smooth of pure dimension $4$.   The map $\rho$ is injective because
of~(6.2.2)-(6.2.3), and thus  
 $\rho$ is an isomorphism between $J\tm\Jh$ and an
irreducible component of $\Mw$ by Zariski's Main Theorem. We finish by
showing that $\rho$ is surjective.  

\proclaim{(6.2.16)Lemma}  
Assume $[F]\in\Mw$.  Then the W.I.T.~holds for $F$, with index $i(F)=1$,
i.e.~$R^j\phih_{*}(\cL\ot\phi^{*}F)=0$ for $j\not=1$ (see~\cite{Muk1}).
Furthermore the Mukai transform of $F$ is given by
$$\wh{F}:=R^1\phih_{*}(\cL\ot\phi^{*}F)=i_{{\hat v},*}(\d),$$
for a certain ${\hat v}\in\Jh$, where $\d$ is a line-bundle of degree $(-1)$ on $C$.
\endproclaim

\demo{Proof}
By Serre duality and slope-stability of $F$,
$$H^2(\cL_{\hat t}\ot F)=\Hom(\cL_{\hat t}\ot F,\cO_J)^{*}=0.$$
Applying Grothendieck-Riemann-Roch one gets that
$$\sum\limits_{i=0}^{1}(-1)^j c_1^{hom}(R^j\phih_{*}(\cL\ot\phi^{*}F))=-\t.$$
Since $\chi(\cL_{\hat t}\ot F)=0$ for all $t\in J$, this implies that there exist
$\xi_1,\xi_2\in\Jh$, with $\xi_1\not\cong\xi_2$, such that
$$H^0(\xi_1^{-1}\ot F)\not=0\not=H^0(\xi_2^{-1}\ot F).$$
By slope-stability of $F$, the map $(\xi_1\op\xi_2)\to F$ is an
isomorphism at the generic point, hence we have an exact sequence
$$0\to\xi_1\op\xi_2\to F\to i_{u,*}(\eta)\to 0,$$
where $\eta$ is a line-bundle on $C$, of degree $1$. This shows that 
$\phih_{*}(\cL\ot\phi^{*}F)=0$, hence the W.I.T.~holds for $F$, with index $1$.
In order to compute $\wh{F}$ pull-back the above sequence to $J\tm \Jh$,
 tensorize with $\cL$, and apply the functor $\phih_{*}$: one gets an exact
sequence 
$$0\to R^1\phih_{*}(\cL\ot\phi^{*}F)\to 
i_{{\hat v},*}(\l)\to\CC_{[\xi_1]}\op\CC_{[\xi_2]}\to 0,$$
where $\l$ is a degree-one line-bundle on $C$. It follows that $\wh{F}$ is as
claimed. 
\qed
\enddemo

Now let's prove that $F$ is isomorphic to $E(x,\yh)^{*}$ for some $(x,\yh)$. By a
Theorem of  Mukai~\cite{Muk1, Cor.(2.4)} and~(6.2.16) we have 
$$F\cong(-1_J)^{*}{\Hat{\Hat F}}:=
(-1_J)^{*}R^1\phi_{*}(\phih^{*}(\cL\ot i_{{\hat v},*}(\d))).$$
Applying  Serre duality one gets 
$$\split
(-1_J)^{*}R^1\phi_{*}(\phih^{*}(\cL\ot i_{{\hat v},*}(\d))) & \cong  
R^1\phi_{*}(\phih^{*}(\cL^{-1}\ot i_{{\hat v},*}(\d)))\\
& \cong (\phi_{*}(\cL\ot\phih^{*}(i_{{\hat v},*}(2p_0-\d)))^{*}
=E(-p_0-\d,{\hat v})^{*}.
\endsplit$$ 
This finishes the proof of Proposition~(6.2.1).
\subhead
6.3. Description of $\Va^0$
\endsubhead
Let $P_{\a}^0$ be the variety
parametrizing isomorphism classes of couples $(G,\psi)$, where
$[G]\in\cM_{\ww}$, and $\psi$ is a quotient appearing in an
exact sequence  
$$0\to\xi\lra \ia^{*}G\brel\psi\over\lra \zeta\to 0, \tag{6.3.1}$$   
such that the following holds:
$$\gather
\text{$\zeta$ is an invertible sheaf of degree $3$,}\tag{6.3.2}\\
\det G\cong\cO_J(\Ta),\tag{6.3.3}\\
\sum c_2(G)-\sum i_{\a,*}c_1(\xi)=0.\tag{6.3.4}
\endgather$$
Clearly $\Pa^0$ is a locally closed subset of a Quot-scheme $\cQ$ over 
the subvariety of $\Mw$ parametrizing sheaves $G$ such that~(6.3.3)
holds (by~\cite{Muk3,(A.6)} there exists a tautological sheaf on
$J\tm\cM_{\ww}$);  let $\Pa$ be the closure of $\Pa^0$
in $\cQ$.  Thus every point of $\Pa$ is represented by a
couple $(G,\psi)$ satisfying the conditions above, except for
Condition~(6.3.2), which is replaced by ``$\zeta$ is a rank-one sheaf
of degree $3$''.  As explained in the introduction to this section,  
$\Pa$ parametrizes sheaves whose moduli belong to $\Va$, and which
are obtained by elementary modification from a sheaf $G$ as above.   
In order to prove this we introduce some notation. 

\proclaim{(6.3.5)Definition}
{\rm Let $\cEa$ be a rank-two vector-bundle on $J\tm\Jh$ such that
$$\cEa_{\hat y}\cong E(-\a-2y,{\hat y}).\tag{6.3.6}$$
 (The existence of $\cEa$ is guaranteed by the existence of a tautological
sheaf on $J\tm\Mw$.) Let $\cGa:=(\cEa)^{*}$.}
\endproclaim 

By~(6.2.2) the map
$$\matrix
\Jh & \lra & \Mw\\
{\hat y} &\mapsto & [\cG^{\a}_{\hat y}]
\endmatrix$$
is an isomorphism between $\Jh$ and
the subvariety of $\Mw$ parametrizing sheaves with trivial determinant.
Thus $\Pa$ is a Quot-scheme over $\Jh$: let 
$$f\cl\Pa\to\Jh\tag{6.3.7}$$
be the natural map. 

Thus if $t=(G,\psi)\in\Pa$, and ${\hat y}=f(t)$, we have
$$\text{$G\cong\cG^{\a}_{f(t)}\cong E(-\a-2y,{\hat y})^{*}$,
 where ${\hat y}=f(t)$.}\tag{6.3.8}$$
Let $F=\cF^{\a}_t$ be the elementary
modification of $G$ associated to~(6.3.1), i.e.~the sheaf fitting into the
exact sequence  
$$0\to F\brel j\over\to G\to i_{\a,*}(\zeta)\to 0.\tag{6.3.9}$$   
Clearly $\cF^{\a}_t$ is the restriction to $J\tm\{t\}$ of an elementary
modification 
$$0\to \cF^{\a}\lra(\id_J\tm f)^{*}\cG^{\a}\brel\Psi^{\a}\over\lra
(i_{\a}\tm\id_{\Pa})_{*}(\zeta^{\a})\to 0.\tag{6.3.10}$$ 

\proclaim{(6.3.11)Lemma}
Keep notation as above. Then $\cF^{\a}$ is a family of torsion-free 
semistable sheaves on $J$ parametrized by $\Pa$, with
$$v(\cF^{\a}_t)=\vv,\quad
\det\cF^{\a}_t\cong\cO_J,\quad 
\sum c_2(\cF^{\a}_t)=0\tag{6.3.12}$$ 
for all $t\in\Pa$. 
\endproclaim 

\demo{Proof}
That $\cF^{\a}$ is a family of sheaves (i.e.~it is flat over $\cO_{\Pa}$) 
follows from flatness over $\cO_{\Pa}$ of the other two sheaves
appearing in~(6.3.10). Let $F$ and $G$ be as above. Since $F$ is a
subsheaf of $G$ it is torsion-free. The formulae of~(6.3.12) hold by~(6.3.2),
(6.3.3) and (6.3.4) respectively. We finish by showing that $F$ is
semistable.  Suppose $F$ is not semistable, and let   
$$0\to I_Z\ot L\brel f\over\lra F\lra I_W\ot L^{-1}\to 0\tag{6.3.13}$$
be desemistabilizing, where $L$ is a line bundle and $Z,W$ are zero-dimensional
subschemes. Thus $L\cdot\T\ge 0$. Since $j\circ\phi$ is non-zero (here $j$ is as
in~(6.3.9)), and $G$ is slope-stable by~(6.2.13), we get $L\cdot\T=0$. Thus
by Hodge index $L^2\le 0$; we claim $L^2=0$. In fact from~(6.3.13) we get 
$$2=c_2^{hom}(F)=\ell(Z)+\ell(W)-L^2.\tag{6.3.14}$$
Hence if $L^2<0$ then $L^2=(-2)$, contradicting~(1.3). Since $L^2=0$ and
$L\cdot\T=0$, we get by Hodge index that $[L]\in\Jh$. 
Thus by~(6.3.14) we have $\ell(Z)\ge\ell(W)$ unless $\ell(Z)=0$.
This shows that $Z$
is empty, because otherwise $I_Z(L)$
does not desemistabilize. The restriction of $j\circ f$ to $\Ta$ maps
$L|_{\Ta}$ to $\xi$: since $\deg(L|_{\Ta})=0$ and $\deg\xi=-1$,  we
get that  $j\circ f$ vanishes along $\Ta$. Thus $j\circ f$ extends to a
map $L(\Ta)\to G$. This is absurd because $G$ is slope-stable by~(6.2.13). 
\qed
\enddemo

By the above lemma the family $\cF^{\a}$ induces a modular map
$$\mu_{\a}\cl\Pa\to\cM.\tag{6.3.15}$$
Let $\mu_{\a}^0$ be the
restriction of $\mua$ to $\Pa^0$. The main result of this subsection is the
following.

\proclaim{(6.3.16)Proposition}
Keep notation as above. Then $\mua^0$ is an isomorphism between
$\Pa^0$ and $\Va^0$. In particular $\mu_{\a}(\Pa)=\Va$, and $\mu_{\a}$
is a birational morphism from $\Pa$ to $\Va$.   
\endproclaim

\demo{Proof}
Let $t\in\Pa$, and set $F=\cF_t^{\a}$, $G=\cG_{f(t)}^{\a}$. The long
exact sequence of Tor's obtained by applying the functor $\ot\T_{\a}$
to~(6.3.9) gives an exact sequence
$$0\to \zeta\ot K_C^{-1}\to\ia^{*}F\to\xi\to 0.\tag{6.3.17}$$
By~(6.2.1) $G$ is locally-free, hence $\xi$ is locally-free. Thus the above
exact sequence shows that $\ia^{*}F$ is singular if and only if $\zeta$ is not
locally-free. Since $F$ is isomorphic to $G$ outside $\Ta$,  we get that  
$$\text{$\ia^{*}\cF_t^{\a}$ is locally-free iff $\cF_t^{\a}$ is locally-free
iff $t\in\Pa^0$.}\tag{6.3.18}$$  
Now assume $t\in\Pa^0$. Since $\deg\xi=-1$, Exact Sequence~(6.3.17)
shows that $\ia^{*}F$ is not semistable, hence
$[F]\in\Va^0$. Thus $\mua^{0}(\Pa^0)\ss\Va^0$. To finish the proof of
the proposition we must define a regular map $\Va^0\to\Pa^0$ inverse to
$\mua^0$. One  proceeds as explained in the
introduction to this section, i.e.~to $[F]\in\Va^0$ we associate the couple
$(G,\psi)$ where $G$ is the sheaf fitting into the exact sequence~(6.1.2),
and $\psi$ is the map of~(6.1.3). Of course we need to prove the following
two results.

\proclaim{(6.3.19)Lemma}   
Let $[F]\in\Va^0$, and assume~(6.1.1) is the desemistabilizing sequence of
$\ia^{*}F$. The sheaf $G$ fitting into Exact Sequence~(6.1.2) is
locally-free and slope-semistable. If $\deg\xi=-1$ (where $\xi$ is the
destabilizing quotient of $\ia^{*}F$)  then $G$ is slope-stable.   
\endproclaim

\demo{Proof}
Since~(6.1.1) is the desemistabilizing sequence, both $\l$ and $\xi$ are
locally-free: from Exact Sequence~(6.1.3) we get that $G$ is locally-free
along $\Ta$. On the other hand it follows from~(2.1.2) and~(4.3.3) that
$F$ is locally-free: since $G$ is isomorphic to $F$ outside $\Ta$, we get that
$G$ is locally-free. Now  assume that $G$ is not slope-semistable. Then 
there exists an injection $L\hra G$, where $L$ is locally-free of rank one
with 
$$L\cdot\T> slope(G):=
{1\over \rk(G)}c_1^{hom}(G)\cdot\T=1.\tag{6.3.20}$$ 
From~(6.1.2) we have an injection $L(-\Ta)\hra F$, hence we get an exact
sequence 
$$0\to M\to F\to I_Z\ot M^{-1}\to 0,\tag{6.3.21}$$
where $M:=L(-\Ta+D)$, for $D$ an effective divisor, and $Z$ a zero-dimensional
subscheme of $J$. From~(6.3.20) and slope-semistability of $F$
we get  $M\cdot\T= 0$.
Hence by Hodge index $M\cdot M\le 0$.  From~(6.3.21)
we get
$$2=c_2^{hom}(F)=-M\cdot M+\ell(Z).$$
Since $M\cdot M\le 0$, and the intersection from is even, either $M\cdot M=0$ or
$M\cdot M=-2$. In the former case Exact Sequence~(6.3.21) shows that $F$
is not (Gieseker-Maruyama) semistable, which is absurd. Thus $M\cdot M=-2$;
since $M\cdot\T=0$, this contradicts~(1.3).  Thus $G$ is
slope-semistable. To finish the proof of the lemma we notice that if
$\deg\xi=-1$ then $v(G)=\ww$, hence $G$ is slope-stable by~(6.2.13).  
\qed
\enddemo

\proclaim{(6.3.22)Lemma}   
Let $[F]\in\Va^0$, and assume~(6.1.1) is the desemistabilizing sequence.
Then $\deg\xi=-1$.  
\endproclaim

\demo{Proof}
Let $d:=\deg\l$: we know $d\ge 1$. To prove the lemma it suffices to show
that $d\le 1$. From~(6.1.2) one gets that
$$c_1^{hom}(G)=\t,\quad c_2^{hom}(G)=2-d.$$
By Lemma~(6.3.19) $G$ is slope-semistable, hence by Bogomolov's theorem
$$0\le 4c_2^{hom}(G)-c_1^{hom}(G)\cdot c_1^{hom}(G)=6-4d.$$
This gives $d\le 1$.
\qed
\enddemo

Given the above lemmas, the argument described in the introduction to this
section shows that the map
$$\matrix
\Va^0 & \lra & \Pa^0 \\
[F] & \mapsto & [(G,\psi)],
\endmatrix$$
where $\psi\cl\ia^{*}G\to \l\ot K_C$ is the map appearing in~(6.1.3), is the
inverse of $\mua^0$. This finishes the proof of the proposition.
\qed
\enddemo
\subhead
6.4. First analysis of $\Pa$
\endsubhead
Let $\Jh[2]_{-{\hat\a}}$ be as in~(1.11).

\proclaim{(6.4.1)Proposition}
The map $f\cl\Pa\to\Jh$ (see~(6.3.7)) is a $\PP^1$-fibration  
away from $\Jh[2]_{-{\hat\a}}$, and  the remaining fibers are isomorphic
to $\PP^2$.     
\endproclaim

\demo{Proof}
Let $\yh\in\Jh$, and suppose $(G,\psi)\in f^{-1}(\yh)$.  By~(6.3.8) the
isomorphism class of $G$ is determined by $\yh$. Similarly, the
isomorphism class of the line-bundle $\xi$ appearing in~(6.3.1) is
completely determined by Equation~(6.3.4): eplicitely, Condition~(6.3.4)
together with~(6.2.3) gives    
$$\xi\cong[-p_0+2\a+y].\tag{6.4.2}$$ 
Thus     
$$f^{-1}(\yh)\cap\Pa^0=\{[\s]\in\PP H^0([p_0-2\a-y]\ot\ia^{*}\cG^{\a}_{\hat y})|
\ \text{$\s$ has no zeroes}\}.$$ 
Exact Sequence~(6.2.6) gives   
$$\multline
0\to H^0([2p_0-2\a])\to
H^0([p_0-2\a-y]\ot \ia^{*}\cG^{\a}_{\hat y})\\
\brel\k\over\lra H^0([2p_0-2\a-2y])\to 0.
\endmultline
\tag{6.4.3}$$ 
(Recall that $(2p_0-2\a)\not\sim K_C$ by~(3.7).)
One verifies easily that for $y$ generic the generic section in the middle 
$H^0$ has no zeroes, hence 
$$f^{-1}(\yh)=\PP H^0([p_0-2\a-y]\ot\ia^{*}\cG^{\a}_{\hat y}).
\tag{6.4.4}$$
Using~(6.4.3) we get that the fibers of $f$ are as
stated in the proposition, and that $f$ is a $\PP^1$-fibration over $\Jh\sm
(\Jh[2]_{-{\hat\a}})$. 
\qed
\enddemo

\proclaim{(6.4.5)Remark}
{\rm Let $\pi_{\Jh}\cl C\tm\Jh\to\Jh$ be the projection. Let  $\cH$ be a
line-bundle on $C\tm\Jh$ such that for  ${\hat y}\in\Jh$ we have
$$\cH_{\hat y}\cong[-p_0+2\a+y].$$
Let $\cSa$ be the sheaf on $\Jh$ defined by
$$\cSa:=\pi_{\Jh,*}(\cH^{-1}\ot(\ia\tm\id_{\Jh})^{*}\cGa).\tag{6.4.6}$$
In the proof above we have shown that over $(\Jh\sm\Jh[2]_{-\hat\a})$
the sheaf $\cSa$ is locally-free of rank two (in fact this is true over all of
$\Jh$, see~(6.6.1)), and that $f^{-1}(\Jh\sm\Jh[2]_{-\hat\a})$ is naturally
isomorphic to $\PP(\cSa|_{(\Jh\sm\Jh[2]_{-\hat\a})})$.}
\endproclaim
\subhead
6.5. The boundary of $\Pa$
\endsubhead
We  analyze $(\Pa\sm\Pa^0)$. First notice that by~(6.3.18) 
$$\Pa\sm\Pa^0=\mua^{-1}(\Ba\cup\Sia)=\mua^{-1}(B\cup\Si).\tag{6.5.1}$$
For $t=(G,\psi)\in\Pa$   let $\zeta_t$ be the
rank-one sheaf $\zeta$ appearing in Exact Sequence~(6.3.1) and
$Tors(\zeta_t)$ be its torsion subsheaf. Then
$$\ell(Tors(\zeta_t))\le 2.\tag{6.5.2}$$ 
In fact since $\ia^{*}G$ is semistable and we have a surjection
$$\ia^{*}G\to\zeta_t/Tors(\zeta_t),$$
the line bundle on the right has degree at least $1$. 
This implies immediately~(6.5.2).
Let
$$\Wa:=\{t\in\Pa|\ \ell(Tors(\zeta_t))=2\},\qquad
\Ya^0:=\{t\in\Pa |\ \ell(Tors(\zeta_t))=1\}.$$ 
By~(6.5.2) $\Wa$ is closed. Clearly $\Ya^0$ is locally closed; we let
$\Ya$ be its closure. By~(6.3.18) and~(6.5.2) we have
$$\Pa\sm\Pa^0=\Ya\cup\Wa.\tag{6.5.3}$$

\proclaim{(6.5.4)Proposition}
Keep notation as above.
\roster
\item
The restriction of $f$ to $\Wa$ is identified with the blow-up of
$\Jh$ at $\Jh[2]$. 
\item
Let $t\in\Wa$. Then $\cF^{\a}_t$ is singular, with
$$Sêing(\cF^{\a}_t)=\{\ia(q_1),\ia(q_2)\},\tag{6.5.5}$$
where $q_1,q_2$ are given by~(4.2.2). 
\item
If $t$ is generic, $\cF_t^{\a}$ is 
stable. Furthermore $\Wa=\mua^{-1}(\Ba)=\mua^{-1}(B)$.
\endroster
\endproclaim

\demo{Proof}
Item~(1). Let $t=(G,\psi)\in\Wa$, and let $\zeta$ be the rank-one sheaf
appearing in~(6.3.1). Then 
$$\ia^{*}G\brel\psi\over\lra \zeta/Tors(\zeta)$$
is a destabilizing quotient, because $\deg(\zeta/Tors(\zeta))=1$. Conversely, assume
$$0\to\ov{\xi}\to\ia^{*}G\to\ov{\zeta}\to 0\tag{6.5.6}$$
is destabilizing. Thus $\ov{\xi},\ov{\zeta}$ are line-bundles of degree $1$, because
$\ia^{*}G$ is semistable by~(6.2.5). Let $\xi$ be a line-bundle fitting into an
exact sequence
$$0\to\xi\to\ov{\xi}\to\cO_Z\to 0,\tag{6.5.7}$$
where $Z$ is a zero-dimensional subscheme of $C$ of length $2$ such that
$$\a+y-\sum i_{\a,*}c_1(\ov{\xi})+\sum i_{\a,*}(Z)=0,\tag{6.5.8}$$
where ${\hat y}=f(t)$. Composing
$\xi\to\ov{\xi}$ with $\ov{\xi}\to\ia^{*}G$ we get an  injection
$\xi\hra\ia^{*}G$. Let $\zeta:=\ia^{*}G/\xi$.  The exact sequence
$$0\to\xi\lra\ia^{*}G\brel\psi\over\lra\zeta\to 0,$$
satisfies~(6.3.2)-(6.3.3)-(6.3.4), i.e.~$(G,\psi)\in\Pa$. Furthermore,
since $\zeta\cong\cO_Z\op\ov{\zeta}$ we have $(G,\psi)\in\Wa$. Thus
we have defined a one-to-one correspondence between $\Wa$ and the set
of triples consisting of $[G]\in\Mw$ (satisfying~(6.3.3)), Exact
Sequence~(6.5.6) and Exact Sequence~(6.5.7) (with $Z$
satisfiyng~(6.5.8)).  We claim that given~(6.5.6) there is one and only one
$Z$ such that~(6.5.8) holds. In fact from~(6.3.8) and~(6.2.6) we see that
$\ov{\xi}\cong[p_0+y]$, hence~(6.5.8) reads 
$$\sum i_{\a,*}(Z)=0.\tag{6.5.9}$$
The map 
$$\matrix
C^{(2)} & \lra & J\\
Z & \mapsto & \sum i_{\a,*}(Z)
\endmatrix$$
is, up to translation, the Abel-Jacobi map, hence it is surjective. Furthermore the fiber over 
$\sum i_{\a,*}(Z)$ is the single cycle $Z$ unless $Z\in|K_C|$. Since $i_{\a,*}(K_C)=2\a$,
and $2\a\not=0$ by~(3.7), we get that $Z$ is uniquely determined by
$\ov{\xi}$, as claimed. Thus we have defined a one-to-one correspondence
between $\Wa$ and the relative Quot-scheme $D_{\a}$ over $\Jh$
parametrizing couples consisting of $[G]\in\Mw$ (satisfying~(6.3.3)) and a
destabilizing sequence of $\ia^{*}G$. In fact this correspondence is
functorial, hence it gives an isomorphism 
$$\matrix
\Wa & \brel\sim\over\lra & D_{\a}\\
(G,\psi) & \mapsto & (\ia^{*}G\to\zeta/Tors(\zeta)).
\endmatrix
\tag{6.5.10}$$
By~(6.2.5) we get that 
$$f^{-1}({\hat y})\cap\Wa\cong
\cases
\text{a single point} & \text{if $2y\not=0$,}\\
\PP^1 & \text{ if $2y=0$.}
\endcases$$
Let $\nu_0\cl \wh{I}_0\to\Jh$ be the blow up of $\Jh[2]$.  To finish the
proof of Item~(1) it suffices to define a map $\rho\cl \wh{I}_0\to\Wa$ such
that  
$$\gather
f\circ\rho=\nu_0\tag{6.5.11}\\
\text{$\rho |_{\nu_0^{-1}({\hat y}_i)}\cl 
\nu_0^{-1}({\hat y}_i)\to f^{-1}({\hat y}_i)$ is an isomorphism,
for ${\hat y}_i\in\Jh[2]$.}\tag{6.5.12} 
\endgather$$
In fact it follows immediately from~(6.5.11)-(6.5.12) that $\rho$ is bijective,
and that its differential is an isomorphism everywhere, hence $\rho$ is an
isomorphism. Let $\cGa$ be the tautological vector-bundle on
$J\tm\Jh$ defined in~(6.3.5).
 By Isomorphism~(6.5.10), in
order to define $\rho$ it suffices to exhibit a line-bundle $\cH$ on
$\wh{I}_0\tm C$ and an injection
$$ \cH\hra (\nu\tm\ia)^{*}\cGa$$
which restricts to a destabilizing subline-bundle of $\ia^{*}\cGa_{\nu(u)}$ for
every $u\in \wh{I}_0$.  Since $\cGa=(\cEa)^{*}$ (see~(6.3.5)) it is
equivalent to give an injection $\cD\hra(\nu\tm\ia)^{*}\cEa$ which
restricts to a destabilizing subline-bundle of  $\ia^{*}\cEa_{\nu(u)}$ for
every $u\in \wh{I}_0$.  Let 
$$\tau,\pi\cl \wh{I}_0\tm C\tm C\to \wh{I}_0\tm C$$
be the projections which ``forget'' the third and second factor   respectively.
By (6.2.8)  we have an isomorphism (up  to tensoring by a line-bundle on
$\wh{I}_0$) 
$$(\nu_0\tm\ia)^{*}\cEa\cong
\cH^1\ot\tau_{*}([\wh{I}_0\tm\D]\ot\pi^{*}\cH^2),\tag{6.5.13}$$
where $\D\ss C\tm C$ is the diagonal, and $\cH^1$, $\cH^2$ are
line-bundles on $\wh{I}_0\tm C$ such that for $u\in \wh{I}_0$ (with
$\nu_0(u)={\hat y}$) 
$$\cH^1_u\cong[-p_0+y],\quad\cH^2_u\cong[2p_0-2y].$$
Arguing as in the proof of~(6.2.9) we get an exact sequence
$$0\to \tau_{*}(\pi^{*}\cH^2)\brel g\over\to
\tau_{*}([\wh{I}_0\tm\D]\ot\pi^{*}\cH^2)\to 
\cH^3\to\cO_{E\tm C}\to 0,$$
where $\cH^3$ is a line-bundle on $\wh{I}_0\tm C$ such that
$\cH^3_u\cong [2y]$ for $u\in \wh{I}_0$, and $E\ss \wh{I}_0$ is the
exceptional divisor of $\nu_0$.   Since $\tau_{*}(\pi^{*}\cH^2)$ is a
line-bundle, and  the map $g$ vanishes to first order along $E\tm C$,
we get a map 
$$\tau_{*}(\pi^{*}\cH^2)\ot[E\tm
C]\to\tau_{*}([\wh{I}_0\tm\D]\ot\pi^{*}\cH^2)$$ 
which is non-zero at every point. Thus by~(6.5.13) we get a map
$$\cD:=\cH^1\ot \tau_{*}(\pi^{*}\cH^2)\ot[E\tm C]
\brel\g\over\lra(\nu_0\tm\ia)^{*}\cEa,$$
injective at every point. The restriction of $\g$ to $\{u\}\tm C$ is a destabilizing sub
line-bundle of $\ia^{*}\cEa_{\nu_0(u)}$ for every $u\in \wh{I}_0$. As
explained above this defines a map $\rho\cl \wh{I}_0\to\Wa$. 
One verifies easily that~(6.5.11)-(6.5.12) hold.  

\n
Item~(2). Let $G=\cG_{f(t)}^{\a}$ and $F=\cF_t^{\a}$. Exact
Sequence~(6.3.17) shows that $Sing(F)=\ia(Z)$, and by~(6.5.9) we
get~(6.5.5). Assume $F$ is not stable, and thus by~(6.3.11) it is strictly
semistable. Let $I_q\ot L\hra F$ be a destabilizing subsheaf. Thus
by~(2.1.2) $[L]\in\Jh$, and $q\in Sing(F)$. By~(6.5.5) we have
$q=\ia(q_i)$ for $i=1$ or $i=2$. 

\n
Item~(3). The key result is the following.

\proclaim{(6.5.14)Claim}
Keep notation as above. Let $t\in\Wa$, and set ${\hat y}:=f(t)$.
There is an injection 
$$I_{\ia(q_i)}\ot L\hra \cF^{\a}_t\tag{6.5.15}$$
if and only if
$$2y \sim q_{3-i}-r\quad\text{for some $r\in C$}.\tag{6.5.16}$$
Furthermore, if $y$ satisfying~(6.5.16) is given, there exists a unique $t\in
f^{-1}({\hat y})$  such that an injection~(6.5.15) exists, and this injection is
unique up to scalars.  
\endproclaim

\demo{Proof of the claim}
Assume~(6.5.15) exists: composing with the map $j$ of~(6.3.9) we get an
injection  $h\cl I_{q_i}\ot L\hra G$, which comes from a (non-zero) map ${\tilde
h}\cl L\to G$, because $G$ is locally-free.  Notice that by slope-stability of
$G$ and a Chern class computation the map ${\tilde h}$ has a single 
zero, which is simple, hence we get an exact sequence
$$0\to L\brel{\tilde h}\over\lra G
\lra I_p\ot\Ta\ot L^{-1}\to 0.\tag{6.5.17}$$ 
Since $t\in\Wa$, the inclusion $\xi\hra\ia^{*}G$ is
the composition 
$$\cO_C(p_0+y)(-q_1-q_2)\hra\cO_C(p_0+y)\hra\ia^{*}G,$$
where the second inclusion is a destabilizing sub-line-bundle $\ia^{*}G$.
The restriction to $\Ta$ of $h$ has image contained
in $\xi$, and hence the restriction to $\Ta$ of ${\tilde h}$ has image
contained in the destabilizing sub-line-bundle of $\ia^{*}G$, and  ${\tilde h}$
must vanish at $\ia(q_{3-i})$. Hence the point $p$ of~(6.5.17) is equal to
$\ia(q_{3-i})$, and by~(6.2.3) we get the equation
$$\a+y=\sum c_2(G)=i_{\a,*}(q_{3-i}+\ia^{*}L),$$
which gives
$$\ia^{*}L\cong[p_0+y-q_{3-i}].\tag{6.5.18}$$
There is a further constraint coming from the fact that
$h^0(L^{-1}\ot G)>0$. Applying the
$Hom(\bu,\cO_J)$-functor to~(6.2.10) we get the exact sequence
$$0\to\cL_{\wh{D}}^{*}\to E(x,\hat y)^{*}\to 
i_{-x,*}[p_0-y]\to 0.$$
(We have used~(6.2.11).) Thus we have an exact sequence
$$\multline 
0\to H^0(L^{-1}\ot\cL_{\wh{D}}^{*})\to 
H^0(L^{-1}\ot E(x,\hat y)^{*})\\
\to H^0(i_{-x}^{*}L^{-1}\ot[p_0-y])
\to H^1(L^{-1}\ot\cL_{\wh{D}}^{*}).
\endmultline$$ 
From this we get that $h^0(L^{-1}\ot E(x,\hat y)^{*})>0$ if and
only if 
$$\text{$[p_0-y-r]\cong i_{-x}^{*}L$ for some $r\in C$.}\tag{6.5.19}$$
Since $i_{-x}^{*}L\cong \ia^{*}L$, the above isomorphism together
with~(6.5.18) gives an equation for $y$ which turns out to be~(6.5.16).
Now let's prove the viceversa. If~(6.5.16) is satisfied then the argument
just given shows that there exists a map ${\tilde h}\cl L\to G$, which
has a single zero at $\ia(q_{3-i})$.  Since
${\tilde h}$ vanishes at $\ia(q_{3-i})$, the restriction of ${\tilde h}$ to
$\Ta$ must be contained in a destabilizing sub-line-bundle of $\ia^{*}G$.
This destabilizing sub-line-bundle determines a point $t\in
f^{-1}(y)\cap\Wa$ such that ${\tilde h}$ ``comes'' from an injection
$I_{\ia(q_i)}\ot L\hra \cF^{\a}_t$. This argument also shows that $t$ is
unique. Injection~(6.5.15) is unique because otherwise we would have
$\cF^{\a}_t\cong I_{\ia(q_i)}\ot L\op I_{\ia(q_i)}\ot L$, which is
impossible by~(3.6).
\qed
\enddemo  

Let us prove the first statement of Item~(3). By the above claim
$$f(\Wa\cap\mua^{-1}(\Si))=
\bigcup\limits_{i=1,2}\{{\hat y}\in\Jh|
\ 2y\sim q_{3-i}-r,\quad\text{some $r\in C$}\}.\tag{ywsub}$$ 
Since $f(\Wa)=\Jh$ we get that $\cF_t^{\a}$ is stable for $t$ generic. 
Let us prove the second statement of Item~(3). For $t\in\Pa$ the sheaf
$F:=\cF^{\a}_t$ is locally-free outside
$\Ta$, hence $\mua^{-1}(B)=\mua^{-1}(\Ba)$. Let $t\in\mua^{-1}(B)$;
by~(4.3.3) we have $\ell(F^{**}/F)=2$, hence~(6.3.17) gives that
$\ell(Tors(\zeta_t))=2$, i.e.~$t\in\Wa$. This proves that
$\mua^{-1}(B)\ss\Wa$. By~Item(2) an open dense subset of $\Wa$ is
contained in $\mua^{-1}(B)$. Since $\mua^{-1}(B)$ is closed we get that 
$\Wa\ss\mua^{-1}(B)$. This finishes the proof of Item~(3). 
\qed
\enddemo

We finish this subsection by observing that Item~(3) of the above proposition
together with~(6.5.1)-(6.5.3) gives that 
$$\Ya=\mua^{-1}(\Sia)=\mua^{-1}(\Si).\tag{6.5.20}$$
\subhead
6.6. A $\PP^1$-bundle mapping to $\Pa$
\endsubhead
We recall that $f\cl\Pa\to\Jh$ is the natural projection
(see~(6.3.7)). Remark~(6.4.5) identifies
$f^{-1}(\Jh\sm\Jh[2]_{-\hat\a})$ with an explicit $\PP^1$-bundle. Of
course by~(6.4.1) it is not true that $\Pa$ is a $\PP^1$-bundle over $\Jh$;
notice also that by~(6.4.1) there is at least one singular point
of $\Pa$ lying over each point of $\Jh[2]_{-\hat\a}$.  
In this subsection we will construct a $\PP^1$-bundle over the blow-up of
$\Jh$ at $\Jh[2]_{-\hat\a}$ which has a birational regular map to $\Pa$,
contracting sixteen $(-1,-1)$ curves. Then we will prove that $\Pa$ is
isomorphic to the contraction of these curves. 

\proclaim{(6.6.1)Claim}
The sheaf $\cSa$ given by~(6.4.6) is  locally-free of rank two. 
\endproclaim

\demo{Proof}
It has been shown in the proof of~(6.4.1) that $\cSa$ is locally-free outside
$\Jh[2]_{-\hat\a}$. Let $U:=(\Jh\sm\Jh[2])$; by Assumption~(3.7)
$\Jh[2]_{-\hat\a}\ss U$, hence it suffices to prove that $\cS^{\a}$ is
locally-free of rank two on $U$. By~(6.2.5) there is an exact sequence  
$$0\to
\cL^{+}\to
(\ia\tm\id_{\Jh})^{*}\cG^{\a}|_{C\tm U}\to 
\cL^{-}\to 0,$$ 
where $\cL^{+}$, $\cL^{-}$ are line-bundles on $C\tm U$ such that 
$$\cL^{+}_{\hat y}\cong[p_0+y],\quad\cL^{-}_{\hat y}\cong[p_0-y].$$
Letting $\pi_U\cl C\tm U\to U$ be the projection, and $\cH$ be the
line-bundle of~(6.4.5), we get an exact sequence 
$$0\to\pi_{U,*}(\cH^{-1}\ot\cL^{+})\to\cS^{\a}|_U\to 
\pi_{U,*}(\cH^{-1}\ot\cL^{-})\to 0.\tag{6.6.2}$$
The sheaves on the left and the right of $\cS^{\a}|_U$ are locally-free of
rank one, hence $\cS^{\a}|U$ is locally-free of rank two.
\qed
\enddemo

Let $\Ra:=\PP(\cSa)$, and $g\cl\Ra\to\Jh$ be the natural
$\PP^1$-fibration. Let $\eta\hra g^{*}\cSa$ be the tautological 
sub-line-bundle. The corresponding
section of $\eta^{-1}\ot g^{*}(\cSa)$ defines 
$$\Phi\in H^0( C\tm\Ra,\pi_{\Ra}^{*}(\eta^{-1})\ot(\id_C\tm
g)^{*}\cH^{-1} \ot(\ia\tm g)^{*}(\cGa)),$$
where $\pi_{\Ra}\cl C\tm\Ra\to\Ra$ is the projection. We view $\Phi$ as
an injective map fitting into an exact sequence
$$0\to\pi_{\Ra}^{*}(\eta)\ot(\id_C\tm g)^{*}\cH\brel\Phi\over\lra
(\ia\tm g)^{*}(\cGa)\lra \cQ\to 0.\tag{6.6.3}$$
 If ${\hat y}\notin\Jh[2]_{-{\hat\a}}$ then for all $t\in g^{-1}{\hat y}$ the restriction
of~(6.6.3) to $C\tm\{t\}$ is an exact sequence of the kind~(6.3.1),
and hence~(6.6.3) induces a map $(\Ra\sm
g^{-1}(\Jh[2]_{-{\hat\a}}))\to(\Pa\sm f^{-1}(\Jh[2]_{-{\hat\a}}))$: this is the
isomorphism of~(6.4.5). Let $\nua$ and $E_k$ be as in~(1.12) and the
subsequent definition: we set
$${\hat y}_k:=\nua(E_k),\qquad k=1,\ldots,16.\tag{6.6.4}$$
Thus $\Jh[2]_{-{\hat\a}}=\{{\hat y}_1,\ldots,{\hat y}_{16}\}$. 

\proclaim{(6.6.5)Claim}
Keep notation as above. There exist $r_1,\ldots,r_{16}\in\Ra$, with 
$g(r_k)={\hat y}_k$, such that 
the following holds. The restriction of~(6.6.3) to $C\tm\{t\}$ is an exact
sequence of the kind~(6.3.1) for all $t\not=
r_1,\ldots,r_{16}$.  
\endproclaim

\demo{Proof}
Let $\s^{+},\s^{-}$ be local generators near ${\hat y}_k$ of the
invertible sheaves  $\pi_{U,*}(\cH^{-1}\ot\cL^{+})$ and
$\pi_{U,*}(\cH^{-1}\ot\cL^{-})$, respectively. 
Then $\s^{+}$ gives a section of
$\cH^{-1}\ot\cL^{+}$ (in a neighborhood of $C\tm{\hat y}_k$) which restricted to  
$C\tm\{{\hat y}_k\}$ generates $H^0(\cH^{-1}\ot\cL^{+}|_{C\tm\{{\hat y}_k\}})$,
while $\s^{-}$ gives a (local) section of $\cH^{-1}\ot\cL^{-}$ which vanishes identically
on $C\tm\{{\hat y}_k\}$. 
Hence by Exact Sequence~(6.6.2) there exists one and only one $r_k\in
g^{-1}({\hat y}_k)$ with the property that for a local section $\tau$ of
$\cS^{\a}$ near ${\hat y}_k$ the corresponding  section of
$\cH^{-1}\ot(\ia\tm\id_{\Jh})^{*}\cG^{\a}$ (defined in a neighborhood of
$C\tm{\hat y}_k$) vanishes on $C\tm{\hat y}_k$ if and only if 
$$\tau({\hat y}_k)\in r_k.$$
This immediately implies the claim.
\qed 
\enddemo

Thus~(6.6.3) defines a rational map $\Ra\cdots>\Pa$ which is regular
outside $\{r_1,\ldots,r_{16}\}$. Let us construct a birational modification
of $\Ra$ which has a regular map to $\Pa$.  
The inclusion
$$r_k\tm E_k\hra\PP(\cS^{\a}_{{\hat y}_k})\tm E_k=\PP(g^{*}\cS^{\a}|_{E_k})$$
 gives rise to an exact sequence   
$$0\to \cO_{E_k}\lra\nua^{*}\cSa|_{E_k}
\brel\rho_k\over\lra \cO_{E_k}\to 0,\tag{6.6.6}$$ 
where $\PP(\ker\rho_k)=r_k\tm E_k$.  Let $\cTa$ be the rank-two locally-free
sheaf on $\Ia$ fitting into the exact sequence
$$0\to\cTa\lra\nua^{*}\cSa\brel\rho\over\lra
\bigoplus\limits_{k=1}^{16}i_{E_{k},*}\cO_{E_k}\to 0,\tag{6.6.7}$$
where $\rho$ is defined by the $\rho_k$'s, and $i_{E_k}\cl E_k\hra \Ia$ is
the inclusion. We set 
$$\Rat:=\PP(\cTa).\tag{6.6.8}$$
Let us show that $\Rat$ naturally
parametrizes a family of quotients of the kind~(6.3.1). Let
$f_{-}\cl\Rat\to\Ia$ be the natural $\PP^1$-fibration and $\eta\hra
f_{-}^{*}\cTa$ be the tautological  sub-line-bundle. By~(6.6.7) we have a map
$\eta\to  f_{-}^{*}\nua^{*}(\cSa)$. The  corresponding section of
$\eta^{-1}\ot f_{-}^{*}\nua^{*}(\cSa)$ defines a section 
$$\Phi^{-}\in H^0( C\tm\Rat,\pi_{\Rat}^{*}(\eta^{-1})\ot
(\id_C\tm (\nua\circ f_{-}))^{*}\cH^{-1} 
\ot(\ia\tm (\nua\circ f_{-}))^{*}(\cGa)),$$
where $\pi_{\Rat}\cl C\tm\Rat\to\Rat$ is the projection.  
Let $E^{-}_k:=f_{-}^{-1}(E_k)$, and 
$E^{-}:=f_{-}^{-1}(E)$.
The map $\Phi^{-}$ vanishes to order one along
$C\tm E^{-}$, hence it gives an exact sequence
$$0\to\pi_{\Rat}^{*}(\eta\ot[E^{-}])\ot(\id_C\tm
(\nua\circ f_{-}))^{*}\cH \brel\Phi^0\over\lra(\ia\tm
(\nua\circ f_{-}))^{*}(\cGa) \brel\Psi\over\lra \cR\to 0.
\tag{6.6.9}$$
As is easily checked the
restriction of~(6.6.9) to $C\tm \{u\}$ is an exact sequence of 
type~(6.3.1) for all $u\in\Rat$. Thus~(6.6.9) induces a regular map
$\d\cl\Rat\to\Pa$ such that $f\circ\d=\nua\circ f_{-}$. An explicit
description of $\d$ is as follows. First,   $(\Rat\sm E^{-})$  is
canonically isomorphic to $\PP(\cSa|_{(\Jh\sm\Jh[2]_{-\hat\a})})$, and the
latter space is isomorphic to $f^{-1}(\Jh\sm\Jh[2]_{-\a})$ by~(6.4.5); the
restriction of $\d$ to    $(\Rat\sm E^{-})$ is the
composition of these two isomorphisms. It remains to describe $\d$
on $E^{-}_k$, for $k=1,\ldots,16$. 
Let 
$$x\in E_k=\PP(H^1(\cO_C))=\PP(H^0(K_C))^{*}=\PP(H^0(K_C)),$$
and let $\CC\s\in\PP(H^0(K_C))$ correspond to $x$ under the
composition of the above isomorphisms. One easily checks that the
restriction of $\d$ to $f_{-}^{-1}(x)$ is an isomorphism to
$\PP(\k^{-1}(\CC\s))$, where $\k$ is the map in~(6.4.3), and    
$\PP(\k^{-1}(\CC\s))\ss f^{-1}({\hat y}_k)$ by~(6.4.4). In particular letting 
$s_k\in f^{-1}({\hat y}_k)$ be the point corresponding to
$\CC\ker(\k)$, we get a section $\G^{-}_k:=\d^{-1}(s_k)$ of  
$f_{-}|_{E^{-}_k}$ which is contracted by $\d$. The curve $\G^{-}_k$
corresponds to the natural exact sequence 
$$0\to \cO_{E_k}(-E_k)\to \cTa|_{E_k}\to \cO_{E_k}\to 0,$$ 
one gets from~(6.6.6) and~(6.6.7), i.e.~$\G^{-}_k:=\PP(\cO_{E_k}(-E_k))$. 
From the above exact sequence we see that $E^{-}_k$ is isomorphic to
$\FF_1$ and that the normal bundle of $\G^{-}_k$ in $E^{-}_k$ is
isomorphic to  $\cO_{\G^{-}_k}(-1)$. Since we also have isomorphisms
$$\cO_{\G^{-}_k}(E^{-}_k)
\brel df_{-}\over\cong\cO_{\G^{-}_k}(f_{-}^{*}E_k)
\cong\cO_{\G^{-}_k}(-1),$$
we get that 
$$N_{\G^{-}_k/\Pam}
\cong\cO_{\G^{-}_k}(-1)\op\cO_{\G^{-}_k}(-1),
\quad k=1,\ldots,16\tag{6.6.10}$$ 
i.e.~$\G^{-}_k$ is a $(-1,-1)$-curve. Let $Contr(\Rat)$ be the
contraction of the sixteen curves $\G^{-}_k$, where $k=1,\ldots,16$. The
map $\d$ induces a regular bijective map
$\d^{\prime}\cl Contr(\Rat)\to\Pa$.

\proclaim{(6.6.11)Proposition}
Keep notation as above. The map $\d^{\prime}$ is an isomorphism.
\endproclaim

\demo{Proof}
It suffices to prove that $\Pa$ is normal. One checks easily that outside
$\bigcup\limits_{k=1}^{16}\G^{-}_k$ the map $\d$ is injective with injective
differential, hence $\Pa$ is smooth outside $\{s_1,\ldots,s_{16}\}$. On the other hand 
$$\dim T_{s_k}\Pa=4,\quad k=1,\ldots,16,$$
because of the exact sequence (see~(6.4.1))
$$0\to T_{s_k}\PP^2\lra T_{s_k}\Pa\brel df(s_k)\over\lra 
T_{{\hat y}_k}\Jh\to 0.$$
Thus a neighborhood of $s_k$ is a $3$-dimensional hypersurface with an isolated
singularity at $s_k$, hence it is normal.
\qed
\enddemo
\subhead
6.7. A modification of $\Rat$ isomorphic to $\Vta$
\endsubhead
By~(6.6.10) each $\G^{-}_k$ is a $(-1,-1)$-curve, hence the
variety obtained by flopping $\Pam$ at each of the 
$\G^{-}_k$ is a smooth complex manifold, which we denote by $\Pap$.
Explicitely, let 
 $\b_{-}\cl\Pab\to\Pa^{-}$ be the blow-up of 
$\G^{-}$, let $\Sib$ be the
exceptional divisor, and $\Sib_k\ss\Sib$ be the component mapping to
$\G^{-}_k$. By~(6.6.10) each $\Sib_k$ is a copy of
$\PP^1\tm\PP^1$, and $\cO_{\Pab}(\Sib_k)$ has degree $(-1)$ on
the curves of any of its two rulings. Let $\b_{+}\cl\Pab\to\Pap$ 
be the contraction of the $\PP^1$'s of $\Sib_k$ (for $k=1,\ldots,16$) 
belonging to the ruling opposite to that which is contracted by $\b_{-}$. 
 The regular map $\muma:=\mua\circ\d$ defines a rational map
$\e^{-}_{\a}\cl\Pam\cdots>\Vta$. Since $\Pap$ is bimeromorphic to
$\Pam$, we get a meromorphic map $\epa\cl\Pap\cdots>\Vta$. The main
result of this subsection is the following.

\proclaim{(6.7.1)Proposition}
The map $\epa$ is regular and it defines
an isomorphism 
$$\epa\cl\Pap\brel\sim\over\lra\Vta.$$
\endproclaim

The proof goes as follows.  Let $\ov{\cU}^{\a}:=
(\id_J\tm\b_{-})^{*}\cU^{\a}$. We will show~(6.7.2) that the sheaf
$\ov{\cU}^{\a}$ is simple if and only if  $t\notin\Sib$.  Starting from
$\ov{\cU}^{\a}$ we will perform two elementary modifications in order
to get simple sheaves. More precisely, the first elementary modification
will be a family $\cV^{\a}$ of sheaves on $J$ parametrized by
$\Pab$, which coincides with $\ov{\cU}^{\a}$ outside $\Sib$, and
such that $\cV^{\a}_t$ is semistable and simple for all $t$ not
belonging to certain fibers of $\b_{-}$, one for each $\Sib_k$ and 
denoted by $Z_k(1)$. The second elementary modification is a family
$\cW^{\a}(\cN)$ of sheaves on $J$ parametrized by an open
neighborhood $\cN$ of $\bigcup\limits_{k=1}^{16}Z_k(1)$. The
sheaf $\cW^{\a}(\cN)_t$ is semistable and simple for all $t\in\cN$,
and $\wt{S}$-equivalent (see~(2.2.2)) to $\cV^{\a}_t$ if $t\notin 
\bigcup\limits_{k=1}^{16}Z_k(1)$. By~(2.3.7) the sheaf
$\cV^{\a}$ defines a regular map 
$(\Pab\sm\bigcup\limits_{k=1}^{16}Z_k(1))\to\Mt$, and 
$\cW^{\a}(\cN)$ defines a regular map $\cN\to \Mt$. By the  
$\wt{S}$-equivalence property stated above these two maps glue
together and define a regular map $\eba\cl\Pab\to\Vta$. We will
verify that $\eba$ descends to a regular map $\epa\cl\Pap\to\Vta$,
and finally we will prove that $\epa$ is an isomorphism.

An equivalent description of $\Pap$ is the following. By~(6.6.11) the space
$\Pa$ is smooth  except at the points $s_1,\ldots,s_{16}$, where it has
quadratic singularities. Furthermore, by~(6.7.25) the divisor $\Ya$ is not
Cartier at each of the $s_i$, in fact its tangent cone consists of two planes
intersecting only at $s_i$. Hence the blow up of $\Pa$ along $\Ya$ is a
smooth resolution of $\Pa$: it is isomorphic to $\Pap$.
\subsubhead
I.~The map $\muma$ lifts to $\Vta$  outside  $\G^{-}$  
\endsubsubhead
This will be an immediate consequence of the following result. 

\proclaim{(6.7.2)Proposition}    
Keep notation as above. Let $t\in\Rat$, and set $u:=\d(t)$, 
$F:=\cF^{\a}_u$. Then $F$ is simple if and only if
$t\notin\G^{-}$. If $t\in\G^{-}_k$ then
$$F\cong I_{\ia(q_1)}\ot L_k(1)\op I_{\ia(q_2)}\ot L_k(2),\tag{6.7.3}$$ 
where $[L_k(i)]=({\hat y}_k-i_{\hat 0}(q_{3-i}))$.
\endproclaim

\demo{Proof}
Since $\d^{-1}(s_k)=\G^{-}_k$, we must show that $F$ is simple
if and only if $u\notin\{s_1,\ldots,s_{16}\}$, and that if
$u=s_k$ then~(6.7.3) holds. Assume $F$ is not simple. Then $F$ is
properly semistable, hence~(6.7.3) holds for some
$L_k(1),L_k(2)\in\Jh$, and by~(6.5.20) we have $u\in\Ya$. Furthermore 
since $F$ is singular at two points, $u\in\Wa$. Thus
$u\in\Ya\cap\Wa$. Let ${\hat y}:=f(u)$; by Claim~(6.5.14) we get that 
$$2y\in\{q_2-r|\ r\in C\}\cap \{q_1-s|\ s\in C\}.$$
The intersection above equals $\Jh[2]\cup\Jh[2]_{\hat\a}$. 
According to~(6.5.14) there is unique point $u_i\in f^{-1}({\hat y})$ such
that there exists an injection $I_{\ia(q_i)}\ot L_k(i)\hra F$ (for $i=1,2$); we
must show that  $u_1=u_2$, if and only if
${\hat y}\in\Jh[2]_{-{\hat\a}}$. If $u_1=u_2=u$ then we
have Isomorphism~(6.7.3), and in particular $L_k(1)\cong L_k(2)^{-1}$;
thus~(6.5.18) forces ${\hat y}\in\Jh[2]_{-{\hat\a}}$. Now assume  ${\hat
y}\in\Jh[2]_{-{\hat\a}}$; by~(6.5.10) the points $u_1$, $u_2$ correspond
to destabilizing sub-line-bundles of  $\ia^{*}\cG^{\a}_{\hat y}$, and
by~(6.2.5) there is only one such sub-line-bundle, hence $u_1=u_2$.
That the isomorphism class of $L_k(i)$ is as claimed follows from~(6.5.18). 
\qed
\enddemo

Let us show that $\muma$ lifts to $\Vta$
outside $\G^{-}$. Let $\cU^{\a}:=(\id_J\tm\d)^{*}\cF^{\a}$: this is a  
 family of sheaves on $J$ parametrized by $\Rat$. By the above proposition
$\cU^{\a}_t$ is simple for all $t\in\Rat\sm\G^{-}$. Furthermore
$\cU^{\a}_t$ is semistable by~(6.3.11) and Equalities~(6.3.12) hold.
By~(2.3.7) $\cU^{\a}$ induces a regular map
$\e^{-}_{\a}\cl(\Rat\sm\G^{-})\to\cMt$ which lifts the restriction
of $\muma$ to $(\Rat\sm\G^{-})$. Since
$\muma(\Rat\sm\G^{-})$ is dense in $\Va$ by~(6.3.16), 
the image of $\e^{-}_{\a}$ is contained in $\Vta$ (and is dense in it). 
\subsubhead
II.~The elementary modification $\cV^{\a}$  
\endsubsubhead
As is easily checked from the definition, the restriction of $\cU^{\a}$ to
$J\tm\G^{-}_k$ is trivial in the $\G^{-}_k$-direction, i.e.~by~(6.7.2) we 
have an exact sequence 
$$0\to \pi_J^{*}(I_{\ia(q_1)}\ot L_k(1))
\lra\cU^{\a}|_{J\tm\G^{-}_k}\brel\psi_k\over\lra      
\pi_J^{*}(I_{\ia(q_2)}\ot L_k(2))\to 0,$$
where $\pi_J$ is projection to $J$. Pulling back to
$J\tm\Sib_k$ we get a similar exact sequence for
$\ov{\cU}^{\a}|_{J\tm\Sib_k}$; let ${\bar\psi}_k$ be the analogue of
$\psi_k$. We define $\cV^{\a}$ to be the sheaf on $J\tm\Pab$ fitting into
the exact sequence 
$$0\to\cV^{\a}\lra\ov{\cU}^{\a}\brel{\bar\psi}\over\lra
\bigoplus\limits_{k=1}^{16}i_{k,*}(\pi_J^{*}(I_{\ia(q_2)}\ot L_k(2)))
\to 0,\tag{6.7.4}$$ 
where $i_k\cl J\tm\Sib_k\hra J\tm\Pab$ is the inclusion, and ${\bar\psi}$ is
defined by the ${\bar\psi}_k$'s. Before stating the main properties of
$\cV^{\a}$ we introduce some notation. If $p\in C$ we let $p^{\prime}\in
C$ be the point such that 
$$p+p^{\prime}\in|K_C|.\tag{6.7.5}$$
Given the identification
$$\G^{-}_k\cong E_k\cong\PP(H^0(K_C))\tag{6.7.6}$$
let $z_k(1),z_k(2)\in\G^{-}_k$ be given by 
$$z_k(i)  \leftrightarrow q_i+q^{\prime}_i\tag{6.7.7}$$
(where $q_1,q_2$ are defined by~(4.2.2)), and let
$Z_k(i):=\b_{-}^{-1}(z_k(i))$.

\proclaim{(6.7.8)Proposition}
Keeping notation as above, $\cV^{\a}$ is a family of torsion-free
semistable sheaves on $J$ parametrized by $\Pab$, with
$$v(\cV^{\a}_t)=\vv,\quad \det(\cV^{\a}_t)\cong\cO_J,
\quad \sum c_2(\cV^{\a}_t)=0\tag{6.7.9}$$
for all $t\in\Pab$. Furthermore $\cV^{\a}_t$ is simple for all
$t\notin\bigcup\limits_{k=1}^{16}Z_k(1)$. 
If $t\in Z_k(1)$ then
$$\cV^{\a}_t\cong I_{\ia(q_1)}\ot L_k(1)\op I_{\ia(q_2)}\ot L_k(2),$$
where $L_k(i)$ is as in~(6.7.2).
\endproclaim

\demo{Proof}
That $\cV^{\a}$ is a family of torsion-free sheaves on $J$ follows from an
easy local computation. (This computation, away from $h^{-1}(h(z_k(1))$, is essentially
written out below in the proof of~(6.7.19).) If
$t\notin\Sib$ we have
$$\cV^{\a}_t\cong\cU^{\a}_{\b_{-}(t)},\tag{6.7.10}$$
thus~(6.7.9) holds
for all $t\notin\Sib$, hence by continuity it holds for all $t\in\Pab$.
Furthermore by~(6.7.10)  and~(6.7.2) the sheaf  $\cV^{\a}_t$ is semistable
and simple for all $t\notin\Sib$.  Thus we are left with the task of describing
$\cV^{\a}_t$ for $t\in\Sib$. To simplify notation set 
$$\l_k(i):=I_{\ia(q_i)}\ot L_k(i).$$
Applying the $\ot\cO_{J\tm\Sib_k}$-functor to~(6.7.4) we get an exact
sequence %
$$0\to\pi_J^{*}\l_k(2)
\ot\pi_{\Sib_k}^{*}\cO_{\Sib_k}(-\Sib_k)  
\to\cV^{\a}|_{J\tm\Sib_k} \to\pi_J^{*}\l_k(1)\to 0,
\tag{6.7.11}$$  
where $\pi_J$, $\pi_{\Sib_k}$ are the projections to $J$ and $\Sib_k$
respectively. Thus $\cV^{\a}_t$ is semistable also for $t\in\Sib$; it is
simple if and only if the restriction of~(6.7.11) to
$J\tm\{t\}$ is non-split (where $k\in\{1,\ldots,16\}$ is the unique
index such that $t\in\Sib_k$).  Let $e_k$ be the extension class
of~(6.7.11): since
 $\Hom(\l_k(1),\l_k(2))=0$, a
simple spectral sequence argument shows that 
$$e_k\in\Ext^1(\l_k(1),\l_k(2))\ot 
H^0(\cO_{\Sib_k}(-\Sib_k)).\tag{6.7.12}$$
The class $e_k$ is equal to a Kodaira-Spencer map. In fact let $t\in\G^{-}_k$: for
$i=1,2$ there is an exact sequence
$$0\to\l_k(i)\to \cU^{\a}_t\to\l_k(3-i)\to 0.\tag{6.7.13}$$
Composing the associated projection
$$\Ext^1(\cU^{\a}_{t},\cU^{\a}_{t})\to
\Ext^1(\l_k(i),\l_k(3-i))$$
with the Kodaira-Spencer map
of $\cU^{\a}_t$, we get a map of rank-two vector-bundles over $\G^{-}_k$
$$ N_{\G^{-}_k/\Rat}\brel\s_k(i)\over\lra
\cO_{\G^{-}_k}
\ot\Ext^1(\l_k(i),\l_k(3-i)).$$
The geometric meaning of $\s_k(i)$ is contained in 
the following observation~\cite{O4,(1.17)}.

\proclaim{(6.7.14)Remark}
{\rm Let $(x,v)\in N_{\G^{-}_k/\Rat}$, i.e.~$x\in\G^{-}_k$ and $v\in
T_x\Rat/T_{\G^{-}_k}$. Then
$\s_k(i)(x,v)=0$ if and only if the inclusion 
$\l_k(i)\hra \cU^{\a}_x$ 
lifts to first-order in the direction $v$. (Meaning in the direction of any ${\tilde v}\in
T_x\Rat$ representing $v$.)}
\endproclaim

Notice that $\s_k(1)$ is an element of the right-hand side of~(6.7.12). 
By~\cite{O4,(1.12)} we have
$$e_k=\s_k(1).\tag{6.7.15}$$

\proclaim{(6.7.16)Claim}
The map $\s_k(i)$ is an isomorphism away from $z_k(i)$.
\endproclaim

\demo{Proof of the claim}
Let $\Wa^{-}:=\d^{-1}(\Wa)$; thus $t\in\Wa^{-}$ if and only if
$Tors(\cR_t)$ has length $2$, where $\cR$ is the sheaf appearing
in~(6.6.9). Arguing as in the proof of Item~(1) of~(6.5.4) we get that 
$$h|_{\Wa^{-}}\cl \Wa^{-}\to \Ia$$
is the blow-up of $\nua^{-1}(\Jh[2])$. In particular
$\Wa^{-}\cap E^{-}=\G^{-}$, and the intersection is transverse. Thus we
have a direct sum decomposition
$$N_{\G^{-}_k/\Rat}\cong 
N_{\G^{-}_k/\Wa^{-}}\op N_{\G^{-}_k/E^{-}_k}
\cong\cO_{\G^{-}_k}(-1)\op\cO_{\G^{-}_k}(-1).\tag{6.7.17}$$
On the other hand the $E_2$-term of the local-to-global spectral sequence
abutting to $\Ext^{\bu}(\l_k(i),\l_k(3-i))$
gives an exact sequence
$$\multline
0\to H^1(Hom(\l_k(i),\l_k(3-i)))
\lra\Ext^1(\l_k(i),\l_k(3-i))\\
\brel\ell\over\lra
 H^0(Ext^1(\l_k(i),\l_k(3-i)))\to 0.
\endmultline
\tag{6.7.18}$$
Let us show that $\ul{\text{away from $z_k(i)$}}$
$$\align
\ell\circ\s_k(i)(N_{\G^{-}_k/E^{-}_k})&\not=0.\tag{6.7.19}\\
\s_k(i)(N_{\G^{-}_k/\Wa^{-}})&=H^1(Hom(\l_k(i),\l_k(3-i)))\ot\cO_{\G^{-}_k}  
\tag{6.7.20}
\endalign$$
To this end we recall the geometric meaning of the map $\ell$ appearing
in~(6.7.18), i.e.~the local analogue of~(6.7.14). Let
$(x,v)$ be an element of the vector bundle (6.6.10): then
$\ell(\s_k(x,v))=0$ if and only if the singularity of $\cU^{\a}_x$ at
$\ia(q_i)$ deforms to first-order in the direction $v$. Let us prove~(6.7.19).
Let $G:=\cG^{\a}_{{\hat y}_k}$, and assume $x\in(\G^{-}_k\sm\{z_k(i)\})$.
Let $S:=h^{-1}(h(x))$: if $\k$  is the map of~(6.4.3) then
$$ S=\PP(\ker(\k)\op\CC\s),\tag{6.7.21}$$
where $\s\in H^0([p_0-2\a-{\hat y}_k]\ot\ia^{*}G)$ is such that $\k(\s)\in H^0(K_C)$
is non-zero. Since $S\ss E^{-}_k$, and we have an isomorphism
$$T_x S\brel\sim\over\lra(N_{\G^{-}_k/E^{-}_k})_x,\tag{6.7.22}$$
in order to prove~(6.7.19) it suffices to analyze the stalk at $(\ia(q_i),x)$ of
the sheaf $\cA:=\ov{\cU}^{\a}|_{J\tm S}$. By definition we have an exact
sequence 
$$0\to\cA\lra G\ot\cO_S\brel\Psi\over\lra(\ia\tm\id_S)_{*}\cQ\to 0$$
(we have dropped the pull-back signs in the middle term). Here $\cQ$ is the sheaf on
$C\tm S$ fitting into the exact sequence
$$0\to [-p_0+2\a+{\hat y}_k]\ot\cO_S(-1)\brel\Phi\over\lra\ia^{*}G\ot\cO_S
\brel\Psi\over\lra\cQ\to 0,$$
where $\Phi$ is the tautological map (see~(6.7.21)). Let $\{w,z\}$ be a
system of local parameters on $J$ centered at $\ia(q_i)$ such that $w=0$ is
a local equation of $\Ta$, and $s$ be a linear coordinate on $S$ centered at
$x$. Thus $\{z,s\}$ is a system of local parameters on $\Ta\tm S$ centered
at $(\ia(q_i),x)$: abusing notation we denote by the same symbols their
pull-back to $C\tm S$ via $(\ia\tm\id_S)$. We claim that choosing a suitable
local trivialization of $\ia^{*}G\ot\cO_S$ around $(q_i,x)$ we have 
$$\Phi(z,s)=(z,s).\tag{6.7.23}$$
To see why let $(\k(\s))=p+p^{\prime}$. Since $x\not=z_k(i)$ we have
$q_i\notin\{p,p^{\prime}\}$. Thus $\s(q_i)$ is not contained in the sub line-bundle
$[p_0+{\hat y}_k]\hra\ia^{*}G$ (see~(6.2.6)); this implies~(6.7.23) because
a generator of $\ker(\k)$ has a simple zero at $q_i$ and has image contained in the
above sub line-bundle. Given~(6.7.23) one computes immediately the stalk of
$\cA$ at $(\ia(q_i),x)$. Explicitely, letting $\cO$ be the local ring of $J\tm
S$ at $(\ia(q_i),x)$ we get an exact sequence
$$\matrix
0&\to&\cA\ot\cO&\lra&\cO^{(3)}&\lra        &\cO              &\to&0.\\
  &     &                &     &(\a,\b,\g)&\mapsto &w\a-z\b-s\g &      &       
\endmatrix
$$
Let $(\ov{\cO},\ov{m})$ be the local ring of $J$ at $\ia(q_i)$. We may assume that
the stalk at $\ia(q_i)$ of~(6.7.13)
is given by
$$0\to\ov{m}\brel f\over\lra\cA\ot\ov{\cO}\brel g\over\lra\ov{\cO}\to 0.
\tag{6.7.24}$$
Let
$$v:={\partial\over\partial s}\in(N_{\G^{-}_k/E^{-}_k})_x.$$
(Yes, we are abusing notation.) By~(6.7.22)
$$\ell\circ\s_k(i)(x,v)=\e\in Ext^1(\ov{m},\ov{\cO}),$$
where $\e$ is obtained as follows. The first-order deformation of $\cA\ot\ov{\cO}$
defined by $\cA\ot\cO$ gives an exact sequence
$$0\to \cA\ot\ov{\cO}\brel \cdot s\over\lra\cA\ot\cO/(s^2)
\lra\cA\ot\ov{\cO}\to 0$$
with extension class ${\tilde \e}\in Ext^1(\cA\ot\ov{\cO},\cA\ot\ov{\cO})$; by
definition $\e=g\circ{\tilde \e}\circ f$, where $f,g$ are the maps
of~(6.7.24).  An easy computation shows that $\e$ is the (non-zero)
generator of $Ext^1(\ov{m},\ov{\cO})$, and this proves~(6.7.19).
Now let's prove~(6.7.20). The left-hand side of~(6.7.20) is contained in the
right-hand side because as $x$ varies in
$W^{-}_{\a}$ the sheaf $\cU^{\a}_x$ is equisingular. Hence since
both sides of~(6.7.20) are line-bundles it suffices to  show that $\s_k(i)$ is
non-zero on $N_{\G^{-}_k/W^{-}_{\a}}$, away from $z_k(i)$.
Assume $(x,v)\in N_{\G^{-}_k/\Wa^{-}}$ and  $\s_k(i)(x,v)=0$.
By~(6.7.14)  $\d_{*}(v)$ is tangent to $\Ya$ at $s_k$; of course
$\d_{*}(v)$ is also tangent to $\Wa$ at $s_k$. The proof of the first
statement of~(6.5.14) works as well for first-order deformations, hence we
get that 
$$(\nua\circ f_{-})_{*}(v)\in 
T_{{\hat y}_k}\{{\hat y}\in\Jh|\ 
2{\hat y}\sim{\hat q}_{3-i}-{\hat r},
\quad\text{some $r\in C$.}\}.$$
The right-hand side equals $\CC\tau$, where $\tau\in H^0(K_C)$ is a
section such that $(\tau)=q_i+q^{\prime}_i$, and we make the
canonical identification $\PP(T_{{\hat y}_k})\cong \PP(H^0(K_C))$. If
$x=p+p^{\prime}$ then $\CC(\nua\circ f_{-})_{*}(v)=p+p^{\prime}$; hence
$\s_k(i)(x,v)=0$ implies that $p+p^{\prime}=q_i+q^{\prime}_i$, i.e.~$x=z_k(i)$.
This proves~(6.7.20).  Claim~(6.7.16) follows at once from~(6.7.19)
and~(6.7.20). 
\qed
\enddemo

 Let $\Ya^{-}\ss\Rat$ be the closure of  the locally-closed subset of points
$t$ such that $Tors(\cR_t)$ has length $1$, where $\cR$ is the sheaf
appearing in~(6.6.9). Thus for $t\in\Ya^{-}$ the sheaf $\cU^{\a}_t$ is
strictly semistable (see~(6.5.20)).

\proclaim{(6.7.25)Claim}
Keeping notation as above, 
$$\Ya^{-}\cap\G^{-}_k=\{z_k(1),z_k(2)\}.$$
Both intersections are transverse.
\endproclaim

\demo{Proof of the claim}
Assume $x\in\Ya^{-}\cap\G^{-}_k$. Since $\dim\Ya^{-}=2$ there
exists $v\in T_x\Pa^{-}$ ``normal'' to
$\G^{-}_k$ which is tangent to a curve $D\ss\Ya^{-}$. Since $\cU^{\a}_t$ is properly
semistable for all $t\in D$ at least one of the
 inclusions $\l_k(1),\l_k(2)\hra\cU^{\a}_x$ lifts to first-order in the direction
determined by $v$.  By~(6.7.14) and~(6.7.16) we get that 
$x\in\{z_k(1),z_k(2)\}$ and that if $x=z_k(i)$ then
$\l_k(3-i)\hra\cU^{\a}_x$ does not lift.   In particular
$\Ya^{-}\cap\G^{-}_k\ss\{z_k(1),z_k(2)\}$. To prove the reverse inclusion
and the transversality of the intersection we consider $\Ya^{-}\cap E^{-}_k$.
First let's show that the intersection is transverse.  Let $x\in E_k$ correspond
to a reduced divisor $p+p^{\prime}\in K_C$ via~(6.7.6), and assume also
that $x\notin\{h(z_k(1),h(z_k(2)\}$
(i.e.~$p+p^{\prime}\not=q_i+q^{\prime}_i$). Arguing as in the proof
of~(6.7.19) one shows that $f_{-}^{-1}(x)$ intersects $\Ya^{-}$
transversely, and hence $E^{-}_k$ is transverse to $\Ya^{-}$. Next, the
image  
$$\d(\Ya^{-}\cap E^{-}_k)\ss f^{-1}({\hat y}_k)
=\PP(H^0([p_0-2\a-{\hat y}_k]\ot\ia^{*}\cG^{\a}_{{\hat y}_k}))
\cong\PP^2$$
is computed by a Chern class computation. In fact let $\pi_C$, $\pi_{\PP^2}$
be the projections of $C\tm\PP^2$ onto $C$ and $\PP^2$ respectively. Then
$$[\d(\Ya^{-}\cap E^{-}_k)]
=\pi_{\PP^2,*}(c_2(\pi_C^{*}[p_0-2\a-y_k]\ot\pi_C^{*}\ia^{*}\cG^{\a}_{{\hat y}_k}
\ot\pi_{\PP^2}^{*}\cO_{\PP^2}(1))),$$
where the multiplicity of the left-hand side equals one because $\Ya^{-}$ is transverse
to $E^{-}_k$. From the above formula one gets that $\d(\Ya^{-}\cap
E^{-}_k)$ is the image of $C$ for the map associated to the linear system
$|2K_C-2\a|$. Furthermore $s_k=\d(\G^{-}_k)$ is the node of
$\d(\Ya^{-}\cap E^{-}_k)$: since the two tangent directions at the node
correspond to $z_k(1)$ and $z_k(2)$, we get that
$\Ya^{-}\cap\G^{-}_k=\{z_k(1),z_k(2)\}$ and that the intersection is
transverse at each of these points.  \qed
\enddemo

Given claims~(6.7.16) and~(6.7.25)  we can determine $\s_k(i)$: let's show
that 
$$\s_k(i)=\vf_k(i)\ot\tau_k(i),\tag{6.7.26}$$
where
$$\quad \vf_k(i)\in\text{GL}_2(\CC),
\quad \tau_k(i)\in H^0(\cO_{\G^{-}_k}(1)),
\quad\tau_k(i)(z_k(i))=0.\tag{6.7.27}$$
Since $\cU^{\a}_t$ is strictly semistable for every $t\in \Ya^{-}$, either for $j=i$ or for
$j=(3-i)$ the inclusion
$\l_k(j)\hra\cU^{\a}_{z_k(i)}$  lifts for every
deformation inside $\Ya^{-}$. By~(6.7.16) the former inclusion is the one that
lifts.  By~(6.7.25)   every direction normal to $\G^{-}_k$ at $z_k(i)$ is 
represented by a vector tangent to $\Ya^{-}$, and thus
by~(6.7.14) we get that
$\s_k(i)$ vanishes at $z_k(i)$. Hence $\s_k(i)=\psi_k(i)\ot\tau_k(i)$, where
$\tau_k(i)$ is as in~(6.7.27), and $\psi_k(i)\in\End(\cO^{(2)}_{\G^{-}_k})$.
By~(6.7.16) $\det\psi_k(i)$ does not vanish outside $z_k(i)$; since
$\det\psi_k(i)$ is a section of $\cO_{\G^{-}_k}$ we get that it is nowhere
zero, and this proves~(6.7.26)-(6.7.27). Now we can finish the proof
of~(6.7.8). Assume $x\in(\ov{\G}_k\sm Z_k(1))$ (for $k=1,\ldots,16$); we
must prove that  $\cV^{\a}_x$ is simple. This is true because by~(6.7.15)
and~(6.7.26)-(6.7.27) the restriction of~(6.7.11) to $J\tm\{x\}$ is non-split
hence simple. \qed
\enddemo
\subsubhead
The elementary modification $\cW^{\a}(\cN)$
\endsubsubhead
Let $\cN\ss\Pab$ be open and such that 
$$Z_k(1)\ss \cN\quad k=1,\ldots,16.\tag{6.7.28}$$
Let $\Yab:=\b_{-}^{-1}(\Ya^{-})$. Let $t\in\Yab\cap\cN$. If
$t\notin\bigcup\limits_{k=1}^{16}Z_k(1)$ then $\cV^{\a}_{t}$ is
strictly semistable and simple by~(6.7.8), hence it has a unique destabilizing
sequence.  On the other hand if $t\in\bigcup\limits_{k=1}^{16}Z_k(1)$
then Exact Sequence~(6.7.13) for $i=1$ lifts to a destabilizing sequence of
$\cV^{\a}_{\tilde t}$ for every ${\tilde t}\in\Yab$ near $t$, by~(6.7.25)
and~(6.7.16).  Thus there is a well-defined exact sequence  
$$0\to \cS\lra\cV^{\a}|_{J\tm(\Yab\cap\cN)}
\brel g\over\lra\cT\to 0\tag{6.7.29}$$
 which restricts to a destabilizing sequence of $\cV^{\a}_t$ for every
$t\in\Yab\cap\cN$. Let
$\cW^{\a}(\cN)$ be the sheaf on $J\tm\cN$ fitting into the exact sequence
$$0\to\cW^{\a}(\cN)\lra\cV^{\a}|_{J\tm\cN}
\brel{\tilde g}\over\lra i_{*}\cT\to 0,\tag{6.7.30}$$
 where $i\cl(\Yab\cap\cN)\hra\cN$ is the inclusion and ${\tilde g}$ is the
restriction to $J\tm(\Yab\cap\cN)$ followed by the map $g$ of~(6.7.29). 

\proclaim{(6.7.31)Proposition} 
With notation as above, $\cW^{\a}(\cN)$ is a family of torsion-free semistable
sheaves on $J$ parametrized by $\cN$, such that
$$v(\cW^{\a}(\cN)_t)=\vv,\quad\det\cW^{\a}(\cN)_t\cong\cO_J,
\quad\sum c_2(\cW^{\a}(\cN)_t)=0$$
for all $t\in\cN$. If $t\in(\cN\sm\Yab)$ then $\cW^{\a}(\cN)_t\cong \cV^{\a}_t$.
Finally, if $\cN$ is a sufficiently small open set
satisfying~(6.7.28) then $\cW^{\a}(\cN)_t$ is simple for all $t\in\cN$.  
\endproclaim

\demo{Proof}
All the statements except the last are standard. Let us prove that if
$t\in(\Yab\cap\cN)$ the sheaf  
$\cW^{\a}(\cN)_t$ is simple, if $\cN$ is sufficiently small. Applying the
functor $\ot\cO_{J\tm(\Yab\cap\cN)}$ to~(6.7.30) we get an exact
sequence 
$$0\to
\cT\ot\pi_{\cN}^{*}\cO_{\cN}(-\Yab)
\to\cW^{\a}(\cN)|_{J\tm(\Yab\cap\cN)}
\to\cS\to 0,$$ 
where $\pi_{\cN}\cl J\tm\cN\to\cN$ is the projection.
Hence if $t\in(\Yab\cap\cN)$ we get an exact sequence
$$0\to\cT_t\to\cW^{\a}(\cN)_t\to\cS_t\to 0.\tag{6.7.32}$$
Thus it suffices to show that the above extension is non-split. Since $\cN$ can
be chosen to be an arbitrarily small neighborhood of 
$\bigcup\limits_{k=1}^{16}Z_k(1)$, it is enough to prove that~(6.7.32) is
non-split for $t\in\bigcup\limits_{k=1}^{16}Z_k(1)$. This will follow from  an
elementary lemma on extensions of sheaves on a variety $X$. Let $T$ be a
smooth curve, $p\in T$, and $\pi_X$, $\pi_T$ the projections of $X\tm T$
on $X$ and $T$ respectively. Let  
$$0\to \cA_2\lra \cE\brel f\over\lra \cA_1\to 0\tag{6.7.33}$$
be an extension of sheaves on $X\tm T$, with $\cA_1$ flat over $T$. Let
$$0\to \cA_2\to \cF\to \cA_1\ot\pi_T^{*}(-p)\to 0\tag{6.7.34}$$
be the pull-back of~(6.7.33) by the inclusion
$\cA_1\ot\pi_T^{*}(-p)\hra\cA_1$.  Thus if $j\cl X\hra X\tm T$ is the
inclusion defined by $j(x):=(x,p)$ and $A_{\ell}:=j^{*}\cA_{\ell}$, the exact
sequence 
$$0\to A_2\to j^{*}\cF\to A_1\to 0$$
is split, hence we have a surjection 
$$\psi\cl  j^{*}\cF\to A_2.\tag{6.7.35}$$
Let $\cG$ be the sheaf
on $X\tm T$ fitting into the exact sequence
$$0\to \cG\lra\cF\brel\Psi\over\lra j_{*}A_2\to 0,\tag{6.7.36}$$
where $\Psi$ is defined by pulling back to $j^{*}\cF$ and then applying
$\psi$.

\proclaim{(6.7.37)Lemma}
Keeping notation as above, $\cG$ is isomorphic to $\cE\ot\pi_T^{*}(-p)$.
\endproclaim

\demo{Proof of the lemma}
The sheaf $\cF$ is a subsheaf of $\cE$ and the inclusion $\cF\ss\cE$ gives an exact
sequence
$$0\to\cF\lra\cE\brel {\tilde f}\over\lra j_{*}A_1\to 0,\tag{6.7.38}$$
where ${\tilde f}$ is $f$ followed by the pull-back map $\cA_1\to j^{*}A_1$. Applying
the functor $j^{*}$ to~(6.7.38) we get an exact sequence
$$0\to A_1\lra j^{*}\cF\brel\vf\over\lra A_2\to 0.$$
The map $\vf$ is equal to the map $\psi$ of~(6.7.35). Thus $\cG$ is
obtained from $\cF$ by performing the elementary modification which is
``inverse'' to~(6.7.38): as is well-known this implies that
$\cG=\cE\ot\pi_T^{*}(-p)$. \qed
\enddemo

We go back to the proof of the proposition. Let $p\in Z_k(1)$ and
$\L\ss\Sib_k$ be the $\PP^1$ containing $p$ and belonging to the
ruling opposite to that of $Z_k(1)$. Let $\vf_k(1)(\L)=\CC e$ where
$\vf_k(1)$ is as in~(6.7.26). Thus $e$ is the class of a non-trvial extension 
$$0\to\l_k(2)\to\cH\to\l_k(1)\to 0.\tag{6.7.39}$$
Let
$$0\to\pi_J^{*}\l_k(2)\to\pi_J^{*}\cH\to\pi_J^{*}\l_k(1)\to 0
\tag{6.7.40}$$
be its pull-back to $J\tm\L$, where $\pi_J\cl J\tm\L\to J$ is the
projection. By~(6.7.11), (6.7.15), (6.7.26) and transversality of the
intersection $\Yab\cap{\bar\G}$ (see~(6.7.25)) the restriction
$\cV^{\a}|_{J\tm\L}$ is the pull-back of~(6.7.40) by the inclusion
$\pi_J^{*}\l_k(1)(-\pi_{\L}^{*}(-p))\hra\pi_J^{*}\l_k(1)$, where
$\pi_{\L}\cl J\tm\L\to \L$ is the projection. Furthermore
$\cW^{\a}(\cN)|_{J\tm\L}$ is the elementary modification of
$\cV^{\a}|_{J\tm\L}$ given by  
$$0\to\cW^{\a}(\cN)|_{J\tm\L}\to\cV^{\a}|_{J\tm\L}\to
j_{*}(\l_k(2))\to 0,$$ 
where $j\cl J\tm\{p\}\hra J\tm\L$ is the inclusion. Hence we are in the
hypotheses of Lemma~(6.7.37), and we get that 
$$\cW^{\a}(\cN)|_{J\tm\L}\cong\pi_J^{*}\cH\ot\pi_{\L}^{*}(-p).$$
In particular $\cW^{\a}(\cN)_p$ is isomorphic to the non-split
extension~(6.7.39), hence it is simple. This finishes the proof of the
proposition. 
\qed
\enddemo
\subsubhead
The map $\epa$
\endsubsubhead
Let $\cN$ be sufficiently small, so that~(6.7.31) applies. Then the sheaf
$\cW^{\a}(\cN)$ defines a regular map $\mua(\cN)\cl \cN\to\Vta$. On the other hand
by~(6.7.8) the restriction of $\cV^{\a}$ to
$J\tm(\Pab\sm\bigcup\limits_{k=1}^{16}Z_k(1))$ defines a regular map 
$$(\Pab\sm\bigcup\limits_{k=1}^{16}Z_k(1))\to\Vta.$$
By~(6.7.31) the two maps coincide on $(\cN\sm\Yab)$ hence they glue
together and they define a regular map $\eba\cl\Pab\to\Vta$.
By~(6.7.26)-(6.7.27) $\eba$ is constant on every $\PP^1\ss\Sib_k$
belonging to the ruling ``opposite'' to that of $Z_k(1)$, hence $\eba$
descends to a regular map $\epa\cl\Pap\to\Vta$. 
\subsubhead
The map $\epa$ away from $(\epa)^{-1}(\Sta)$
\endsubsubhead
We will prove that
$$\text{the restriction of $\epa$ to $(\Pap\sm\Yap)$ is an isomorphism
onto $(\Vta\sm\Sit)$.}
\tag{6.7.41}$$
 Since $\Yap=(\epa)^{-1}(\Sit)$ and
$\pit\cl\cMt\to\cM$ is an isomorphism outside $\Si=\pit(\Sit)$, it is
equivalent to show that the restriction of $\mua$ to $(\Pa\sm\Ya)$ is an
isomorphism onto $(\Va\sm\Si)$.  One defines an inverse of this map
proceeding as in the proof of~(6.3.16); what is needed is an analogue
of~(6.3.22) valid for all $[F]\in(\Va\sm\Si)$. This is the content of the
following lemma.

\proclaim{(6.7.42)Lemma} 
Let $[F]\in\cM$ and assume $F$ is stable. Then there is at most one exact sequence
$$0\to \l\to\ia^{*}F\to\xi\to 0\tag{6.7.43}$$
such that $\l$, $\xi$ are rank-one sheaves, with  $\deg(\l)=1$, $\deg(\xi)=-1$.
\endproclaim

\demo{Proof of the lemma}
If $\ia^{*}F$ is locally-free the result follows immediately from~(6.3.22),
hence we assume $\ia^{*}F$ is singular. Thus $[F]\in(\Ba\sm\Si)$. Consider
the canonical exact sequence
$$0\to F\lra F^{**}\brel\phi\over\lra Q\to 0.\tag{6.7.44}$$
By~(4.3.3) we have $Q=\CC_{x_1}\op\CC_{x_2}$, where $x_1+x_2=0$ and
$x_1\not=x_2$. By the same proposition
$$F^{**}\cong L\op L^{-1},\quad L^{\ot 2}\not\cong\cO_J\tag{6.7.45}$$
or else there is a non-trivial extension 
$$0\to L\to F^{**}\to L\to 0,\quad L^{\ot 2}\cong\cO_J.\tag{6.7.46}$$
Since $F$ is stable 
$$\ker(\phi_{x_i})\cap L_{x_i}^{\pm 1}=\{0\},\tag{6.7.47}$$
where $\phi_{x_i}$, $L_{x_i}$ are the fibers of $\phi$ and $L$ over
$x_i$. Applying the functor $\ia^{*}$ to~(6.7.44) we get an exact
sequence 
$$0\to\ia^{*}Q\lra\ia^{*}F\lra\ia^{*}F^{**}\brel\ia^{*}\phi\over\lra\ia^{*}Q
\to 0.$$
Let $V:=\ia^{*}F/\ia^{*}Q$; by~(6.7.47) if $\eta$ is a line-bundle on $C$
such that $\Hom(\eta,V)\not=0$ then $\deg(\eta)\le(-1)$ (i.e.~$V$ is
semistable). Any rank-one subsheaf of $\ia^{*}F$ is of the form $\eta\op T$,
where $\eta$ is a sub-line-bundle of $V$ and $T\ss \ia^{*}Q$. Assume that
$\deg(\eta\op T)=1$: since  $\deg\eta\le(-1)$ we get that $\ell(T)=2$.
Thus $x_1,x_2\in\Ta$ so that $x_i=q_i$ (where $q_i$ is given by~(4.2.2)),
$T=\CC_{x_1}\op\CC_{x_2}$, and $\deg\eta=(-1)$. This gives a bijection
between the set of Exact Sequences~(6.7.43) and the set of sub-line-bundles
$\eta\ss V$ with $\deg(\eta)=(-1)$. We will prove that if~(6.7.45) holds
then there is at most one such sub-line-bundle. The proof  of the analogous
property when~(6.7.46) holds is left to the reader. Assume there are two
distinct sub-line-bundles $\eta_1\ss V$, $\eta_2\ss V$ of degree $(-1)$:
since $\deg V=(-2)$ we get that
$$V\cong\eta_1\op\eta_2.\tag{6.7.48}$$
Let $\ell:=\ia^{*}L$; composing  $\eta_j\hra V$ with the injection $V\to\ia^{*}F^{**}$
we get an injection $\eta_j\hra(\ell\op\ell^{-1})$. It follows from~(6.7.47)
that  the projections $\eta_j\to\ell^{\pm 1}$ are non-zero, hence
$$\eta_j\cong\ell(-p_j)\cong\ell^{-1}(-r_j)$$
for some $p_j,r_j$. Thus
$$\ell^{\ot 2}\cong[p_j-r_j].\tag{6.7.49}$$
Since $\ell^{\ot 2}\not\cong\cO_C$ this implies that 
$$r_2=p^{\prime}_1,\quad p_2=r^{\prime}_1.\tag{6.7.50}$$
On the other hand $\det(V)\cong[-q_1-q_2]$ and thus~(6.7.48) gives
$$\ell^{\ot 2}(-p_1-p_2)\cong[-q_1-q_2].$$
Using~(6.7.49) and~(6.7.50) we get that $q_1+q_2\in|K_C|$, a
contradiction. \qed
\enddemo
\subsubhead
The restriction of $\epa$ to $\Yap$
\endsubsubhead
We begin by describing  $\Yab$. 
Let $\pi_{\a}\cl T_{\a}\to C\tm\Jh$ be the $\PP^1$-fibration with fiber
$\PP(\Ext^1(I_{-p}\ot\cL_{\hat w}^{-1},I_p\ot\cL_{\hat w}))$ over
$(p,{\hat w})$. Thus we have a quotient map
$$\z_{\a}\cl T_{\a}\to\Sta,\tag{6.7.51}$$
corresponding to the equivalence relation generated by setting $[e_F]
\equiv[e^{\bot}_F]$ for
$[e_F]\in\PP(\Ext^1(I_{\ia(q_2)}\ot\cL_{\hat
w}^{-1},I_{\ia(q_1)}\ot\cL_{\hat w}))$, and $e^{\bot}_F$ generating the
annihilator of $e_F$ in the Serre duality pairing.  We define a regular map
$$\rho_{\a}\cl\Yab \to T_{\a}$$
as follows. If $t\notin \bigcup\limits_{k=1}^{16}Z_k(1)$ then $\cV^{\a}_t$
is simple and semistable, hence there exists a unique destabilizing exact
sequence 
$$0\to I_{-p}\ot\cL_{\hat w}^{-1}\to \cV^{\a}_t 
\to I_p\ot\cL_{\hat w}\to 0,$$
and this determines a point of $T_{\a}$.  If
$t\in Z_k(1)$ the sheaf $\cW^{\a}(\cN)_t$ is simple and semistable,
hence it  fits into a unique destabilizing exact sequence which determines an
extension class $e_{\cW^{\a}(\cN)_t}\in
\PP(\Ext^1(I_{\ia(q_1)}\ot\cL_{\hat w}^{-1},I_{\ia(q_2)}\ot\cL_{\hat
w}))$ for a suitable $w$. We set 
$\rho_{\a}(t):=e_{\cW^{\a}(\cN)_t}^{\bot}$.  That $\rho_{\a}$ is regular
follows from~(6.7.31). Clearly
$$\eba|_{\Yab}=\z_{\a}\circ\rho_{\a}.\tag{6.7.52}$$

\proclaim{(6.7.53)Lemma}
Keeping notation as above,
$$\Im(\pi_{\a}\circ\rho_{\a})=
\{(p,{\hat w})|\ \text{$p-2w\sim r-2\a$ for some $r\in C$}\}=:\Xi_{\a}.
\tag{6.7.54}$$
Furthermore $\pi_{\a}\circ\rho_{\a}$ is the blow-up of $\Xi_{\a}$ with center 
$$\O_{\a}:=\{(q_i,{\hat y}_k-i_{\hat 0}(q_{3-i})=(q_i,[L_k(i)])\},
\quad i=1,2,\quad k=1,\ldots,16.$$
\endproclaim

\demo{Proof of the lemma}
First we show that $\Im(\pi_{\a}\circ\rho_{\a})\ss\Xi_{\a}$. Assume
$t\in(\Yab\sm\ov{\G})$, and let
$I_p\ot\cL_{\hat w}\hra\cV^{\a}_t$ be destabilizing. Let 
$u:=\delta\circ\b_{-}(t)$. Then we have an injection 
$$I_p\ot\cL_{\hat w}\hra\cF^{\a}_u.\tag{6.7.55}$$
Composing with
$\cF^{\a}_u\hra\cG^{\a}_{f(u)}$ we get that
$$h^0(\cL_{\hat w}^{-1}\ot\cG^{\a}_{f(u)})>0,$$
hence by~(6.5.19) we get that
$$\text{$[w]\cong\ia^{*}\cL_{\hat w}\cong[p_0-\wh{f(u)}-r]$ for some $r\in C$.}
\tag{6.7.56}$$
On the other hand~(6.7.55) gives that the line-bundle $\xi$ appearing
in~(6.3.1) (with $G=\cG^{\a}_{f(u)}$) must be equal to $\cO_C(w-p)$ and
thus from~(6.4.2) we get that
$$w-p\sim -p_0+2\a+\wh{f(u)}.\tag{6.7.57}$$
 From~(6.7.56)-(6.7.57) we get that
$\Im(\pi_{\a}\circ\rho_{\a})\ss\Xi_{\a}$. Since $\Xi_{\a}$ is smooth
irreducible and $\Yab$ is a surface smooth at points of $Z_k(j)$ (for $j=1,2$
and $k=1,\ldots,16$), in order to finish the proof of the lemma it suffices to
show that for   $(p,{\hat w})\in\Im(\pi_{\a}\circ\rho_{\a})$ we have
$$(\pi_{\a}\circ\rho_{\a})^{-1}(p,{\hat w})=
\cases
\text{one point} & \text{if  $(p,{\hat w})\notin\O_{\a}$,}\\
Z_k(i) & \text{if  $(p,{\hat w})\in\O_{\a}$.}\\ 
\endcases$$
Let $t_1,t_2\in(\pi_{\a}\circ\rho_{\a})^{-1}(p,{\hat w})$, and set
$u_j:=\delta\circ\b_{-}(t_j)$. It suffices to show
that $u_1=u_2$. From~(6.7.57) we get that $f(u_1)=f(u_2)$. Thus each of
$u_1$, $u_2$ corresponds to an Exact Sequence~(6.3.1) with
$G=\cG^{\a}_{{\hat y}_0}$, where ${\hat y}_0:=f(u_1)=f(u_2)$. An easy
stability argument shows that $h^0(\cL_{\hat w}^{-1}\ot\cG^{\a}_{{\hat
y}_0})=1$,   and thus there is a unique   Exact Sequence~(6.3.1)
(with $G=\cG^{\a}_{f(u)}$) such that the elementary  modification $F$
fitting into~(6.3.9) is strictly semistable. This shows that $u_1=u_2$. 
 \qed
\enddemo
Now we prove that
$$\text{$\epa|_{\Yap}\cl \Yap\to\Vta\cap\Sit$ is bijective.}
\tag{6.7.58}$$  
The map is surjective by~(6.5.20). Let us prove injectivity.  First notice
that 
$$\Yap=\Yab/\sim,\tag{6.7.59}$$
where $\sim$ is the equivalence relation identifying $t\in Z_k(1)$ with the point
$t^{\prime}\in Z_k(2)$ belonging to the same $\PP^1$ of the ruling opposite to that of
$Z_k(i)$. By~(6.7.52) and~(6.7.59)
injectivity of $\e_{\a}$ restricted to $\Yap$ will follow from
$$\text{$\z_{\a}\circ\rho_{\a}(t_1)=\z_{\a}\circ\rho_{\a}(t_2)$ implies
$t_1\sim t_2$, for $t_1,t_2\in\Yab$.}\tag{6.7.60}$$
If $\rho_{\a}(t_1)=\rho_{\a}(t_2)$  it follows
immediately from~(6.7.53) that~(6.7.60) holds. Thus we may assume that 
$\rho_{\a}(t_1)\equiv\rho_{\a}(t_2)$ but $\rho_{\a}(t_1)\not=\rho_{\a}(t_2)$.
Without loss of generality we may suppose that
$$\pi_{\a}\circ\rho_{\a}(t_j)=(q_j,(-1)^j{\hat w}),\quad {\hat w}\in\Jh.$$
By~(6.7.54) we get that
$$q_j-2(-1)^j w\sim r_j-2\a,\quad j=1,2,\quad r_j\in C,\tag{6.7.61}$$
and hence $r_1+r_2\in |K_C+2\a |$.  Thus 
$$\align
r_j & \sim K_C-q_j,\quad j=1,2\text{ or}\tag{6.7.62}\\
r_j & \sim K_C-q_{3-j}\quad j=1,2.\tag{6.7.63}
\endalign$$
Let $u_j:=\d\circ\b_{-}(t_j)$. If~(6.7.62) holds then
$$2f(u_1)=2f(u_2)=-2\a,$$
hence $t_1,t_2\in{\bar\G}$ and one readily verifies that $t_1\sim t_2$. Now
assume~(6.7.63) holds. Then 
$$\wh{f(u_j)}=w+i_{\a,*}(q_{3-j}-p_0),$$
and thus $2f(u_j)\not=-2\a$. Hence 
$$\cV^{\a}_{t_j}\cong\cF^{\a}_{\d\circ\b_{-}(t_j)}=:F_j,$$ 
Thus we have a non-split exact sequence
$$0\to I_{\ia(q_j)}\ot\cL_{\hat w}\brel\phi_j\over\lra F_j
\lra  I_{\ia(q_{3-j})}\ot\cL_{\hat w}\to 0.\tag{6.7.64}$$
(By~(6.7.63)-(6.7.61) we have $2w=0$.) Let
$e_j\in\Ext^1(I_{\ia(q_{3-j})}\ot\cL_{\hat w},I_{\ia(q_j)}\ot\cL_{\hat w})$ be the
extension class of~(6.7.64): we must prove that 
$$\langle e_1,e_2\rangle\not=0,\tag{6.7.65}$$
where $\langle\ ,\ \rangle$ is the Serre duality pairing. The local-to-global
spectral sequence abutting to $\Ext^{\bu}(\ ,\ )$
gives a natural isomorphism 
$$\g_j\cl \Ext^1(I_{\ia(q_{3-j})}\ot\cL_{\hat w},I_{\ia(q_j)}\ot\cL_{\hat w})
\brel\sim\over\lra H^1(\cO_J).$$
commuting with Serre duality.
Since a one-dimensional subspace $\ell\ss H^1(\cO_J)$ is the Serre-duality
annihilator of itself, we get that~(6.7.65) is equivalent to 
$$[\g_1(e_1)]\not=[\g_2(e_2)].\tag{6.7.66}$$
We claim
that via the composition of canonical isomorphisms
$$\PP(H^1(\cO_J))\cong\PP(H^1(\cO_C))\cong\PP(H^0(K_C)),$$
we have 
$$[\g_j(e_j)]=q_{3-j}+q^{\prime}_{3-j}.\tag{6.7.67}$$
By~(3.7) this will imply~(6.7.66). In order to prove~(6.7.67)   we set
$G_j:=\cG^{\a}_{f(u_j)}$. Composing the map $\phi_j$ of~(6.7.64) with the
inclusion $F_j\hra G_j$ we get an injection $\cL_{\hat w}\hra G_j$, which
gives rise to a (non-split) exact sequence
$$0\to \cL_{\hat w}\to G_j\to I_{\ia(q_{3-j})}\ot\cL_{\hat w}(\Ta)\to 0.$$
The corresponding extension class $f_j$ is a generator of
$\Ext^1(I_{\ia(q_{3-j})}\ot\cL_{\hat w}(\Ta),\cL_{\hat w})$. Applying the
$\Hom(\bu,\cL_{\hat w})$-functor to the exact sequence
$$0\to I_{\ia(q_{3-j})}\ot\cL_{\hat w}\to I_{\ia(q_{3-j})}\ot\cL_{\hat w}(\Ta)
\to  i_{\a,*}(\CC_{q_{3-i}}\op K_C(-q_{3-j}))\ot\cL_{\hat w}\to 0$$
we get an exact sequence
$$\multline
\Ext^1(I_{\ia(q_{3-j})}\ot\cL_{\hat w}(\Ta),\cL_{\hat w})\brel\b_j\over\lra
\Ext^1(I_{\ia(q_{3-j})}\ot\cL_{\hat w},\cL_{\hat w})\\
\lra\Ext^2(i_{\a,*}( \CC_{q_{3-i}}\op K_C(-q_{3-j}))\ot\cL_{\hat w},
\cL_{\hat w})
\brel\l_j\over\lra\Ext^2(I_{\ia(q_{3-j})}\ot\cL_{\hat w}(\Ta),
\cL_{\hat w})\to 0.
\endmultline
\tag{6.7.68}$$
Given the natural isomorphism
$$\Ext^1(I_{\ia(q_{3-j})}\ot\cL_{\hat w},I_{\ia(q_j)}\ot\cL_{\hat w})
\cong  \Ext^1(I_{\ia(q_{3-j})}\ot\cL_{\hat w},\cL_{\hat w})$$
we have $e_j=\b_j(f_j)$. The transpose of $\l_j$ is given by the restriction
map 
$$H^0(I_{\ia(q_{3-j})}(\Ta))
\brel\l^{*}_j\over\lra H^0(C;\CC_{q_{3-i}}\op K_C(-q_{3-j})),$$
which has image equal to $H^0(C;\CC_{q_{3-i}})$. By~(6.7.68) we get that
$$\text{Ann}(\g_j(e_j))=\Im(H^0(C;K_C(-q_{3-j}))\hra H^0(C;K_C)).$$
(We make the identification $H^1(\cO_J)^{*}\cong H^1(\cO_C)^{*}\cong H^0(K_C)$.)
This proves~(6.7.67) and finishes the proof of~(6.7.58). 
\subsubhead
Proof of Proposition~(6.7.1)
\endsubsubhead
By~(6.7.41)-(6.7.58) we see that $\epa$ is a bijection between $\Yap$
and $\Vta$. Since $\Pap$ is smooth, in order to prove the proposition
it suffices to show that $\epa$ has injective differential everywhere.
This is true away from $\Yap$ by~(6.7.41), hence we are left with the task
of proving that
$$\text{$d\epa(t)$ is injective for $t\in\Yap$.}\tag{6.7.69}$$ 
First we notice that we have an equality of Cartier divisors
$$(\epa)^{*}\Sit=\Yap.\tag{6.7.70}$$
In fact by~(6.7.53) $\Yap$ is irreducible, hence it suffices to show that for
some $t\in\Yap$ there exists $v\in T_{t}\Pap$ such that 
$$d\epa(v)\notin T_{\epa(t)}\Sit.\tag{6.7.71}$$
In order to exhibit such a $t$ we let $\G^{+}_k:=\b_{+}({\ov{\G}_k})$ and
$\G^{+}$ be the union of the $\G^{+}_k$'s (for $k=1,\ldots,16$): thus 
$\G^{+}_k$ is a $(-1,-1)$-curve of $\Pap$. By~(6.7.25) and~(6.7.59)  
$$\text{$\G^{+}\ss\Yap$ and is the nodal curve of $\Yap$.}\tag{6.7.72}$$
Assume $\cN$ is
sufficiently small (as in~(6.7.31)): then~(6.7.71) holds for all
$t\in(\b_{+}(\cN)\sm\G^{+})$ because $\cW^{\a}(\cN)_{t}$ is simple. Now
let  $t\in( \Yap\sm\G^{+})$, i.e.~in the smooth
locus of $\Yap$. By~(6.7.70) in order to prove~(6.7.69) it suffices to verify
that the restriction of $d\epa(t)$ to $T_t\Yap$ is injective: this is
straightforward. If instead $t\in\G^{+}$ one
verifies by an explicit computation that~(6.7.69) holds. 
\subhead
6.8. Topological results
\endsubhead
First  notice that $\Vta$ is isomorphic to $\Pam$ by~(6.7.1),  and that
$\Pap$ is birational to $\Pam$, hence $\Vta$ is birational to $\Pam$. Since
$\Pam$ is a $\PP^1$-bundle over $\Jh$ (see~(6.6.8)) it is irreducible, and
thus 
$$\text{$\Vta$ is irreducible.}\tag"$\text{(6.8.1)}^{*}$"$$
Since both $\Vta$ and $\Pam$ are smooth 
we have $H^1(\Vta)\cong H^1(\Pam)$. Composing the isomorphisms
$$H^1(\Jh)\overset\nua^{*}\to{\overset\sim\to{\lra}}H^1(\Ia)
\overset f_{-}^{*}\to{\overset\sim\to{\lra}}H^1(\Pam)$$
we get that  $b^1(\Pam)=4$, hence
$$b^1(\Vta)=4.\tag"$\text{(6.8.2)}^{*}$"$$

Let $\D_k:=\epa(\G^{+}_k)$ (for $k=1,\ldots,16$) 
and $\D:=\epa(\G^{+})$. (Thus  $\D_k$ is a $(-1,-1)$  curve
of $\Vta$.) Let $\L(\Vta)\ss H^2(\Vta;\ZZ)$ be the
subgroup defined by  
$$\L(\Vta):=\{\tau\in H^2(\Vta;\ZZ)|
\ \int\limits_{\D_j}\tau=\int\limits_{\D_k}\tau,\quad j,k=1,\ldots,16\}.$$
We claim that
$$\Im(H^2(\Mt;\ZZ)\to
H^2(\Vta;\ZZ))\ss\L(\Vta),\tag"$\text{(6.8.3)}^{*}$"$$ 
where the map between $H^2$'s is given by restriction. 
To prove~(6.8.3) we notice that by~(6.7.72) and~(6.7.58) each  $\D_k$
belongs to $\Sit$. Furthermore each $\D_k$ is a fiber of 
$\Sit\overset\pit\to\lra\Sigma$. Over $\Sigma\sm\O$ (where $\O$ is the
set of equivalence classes of sheaves $(I_p\ot L)\op(I_p\ot L)$) the map
$\pit$ is a locally-trivial fibration, hence all the fibers represent the same
homology class in $H_2(\cMt;\ZZ)$. Thus the integral over $\D_k$ of a class
in $H^2(\Mt;\ZZ)$ is independent of $k$, and this proves~(6.8.3).
We will need to know that
$$\L(\Vta)\cong\ZZ^{8}.\tag"$\text{(6.8.4)}^{*}$"$$
In order to prove it we first notice that there are canonical isomorphisms
$$H^2(\Vta;\ZZ)\cong H^2(\Vta\sm\D;\ZZ)\cong
H^2(\Pam\sm\G^{-};\ZZ)\cong H^2(\Pam;\ZZ).\tag{6.8.5}$$
In fact the first map is an isomorphism because $\Vta$ is smooth and
$\D$ has (complex) codimension $2$, the second map is an isomorphism
because  $(\Vta\sm\D)\cong(\Pam\sm\G^{-})$, and the third map is an
isomorphism because $\Pam$ is smooth and $\G^{-}$ has
(complex) codimension $2$. The last term on the right of~(6.8.5) is
isomorphic to $\ZZ^{23}$ because $\Pam$ is a $\PP^1$-bundle over
$\wh{I}_{\a}$, hence in order to prove ~(6.8.4) we must
show that the classes in $H_2(\Vta;\ZZ)$ represented by
$\D_1,\ldots,\D_{16}$ are independent. Let $\ov{E}_k\ss\Pab$ be the strict
transform of $E^{-}_k$ (for the blow-up $\b_{-}$), let
$E^{+}_k:=\b_{+}(\ov{E}_k)$, and $D_k:=\epa(E^{+}_k)$. One easily
checks that  
$$\langle c_1(D_i),\D_j\rangle=\d_{i,j},$$ 
where $\d_{i,j}$ is Kronecker's symbol. This shows that the homology
classes represented by $\D_1,\ldots,\D_{16}$ are independent, and
proves~(6.8.4).

Now we give results on $\Vta\cap\Sta$. By~(6.7.59), (6.7.72), (6.7.58)
and~(6.7.53) we get the following result. 

\proclaim{$\text{(6.8.6)}^{*}$Proposition}
The intersection $\Vta\cap\Sta$ contains $\D$   as a double curve, and
is smooth away from $\D$. The normalization of $\Vta\cap\Sta$ is
isomorphic to the blow up of $\Xi_{\a}$ at the points of $\O_{\a}$ (see the
statement of~(6.7.53)). 
\endproclaim

Now we prove that
$$b^1(\Sta\cap\Vta)=24.\tag"$\text{(6.8.7)}^{*}$"$$
We identify $\Sta\cap\Vta$ with
$\Yap$, according to~(6.7.58). Let $E_k(i)\ss\Yab$ be the exceptional
divisor of $\pi_{\a}\circ\rho_{\a}$ mapping to the point of $\O_{\a}$
indicized by $k,i$ (see~(6.7.53)). By~(6.7.59) the smooth locus of $\Yap$
is isomorphic to $(\Yab\sm\bigcup\limits_{k,i} E_k(i))$. The long exact
sequence of the couple $(\Yap,sm(\Yap))$ (we let $sm(X)$ be the smooth
locus of a variety $X$) gives an exact sequence
$$0\to H^1(\Yap,sm(\Yap))\lra H^1(\Yap)\overset\phi\to{\lra} 
H^1(\Yab\sm\bigcup\limits_{k,i} E_k(i))\to 0.\tag{6.8.8}$$  
We have
$$H^1(\Yab\sm\bigcup\limits_{k,i} E_k(i))\cong H^1(\Xi_{\a}\sm\O_{\a})
\cong H^1(\Xi_{\a}),$$
where $\Xi_{\a}$ is as in~(6.7.53). Below we will prove that
$$b^1(\Xi_{\a})=8.\tag{6.8.9}$$
Granting this, we finish the computation of $b^1(\Sta\cap\Vta)$ as
follows. The first term appearing in~(6.8.8) is computed applying excision and
K\"unneth: the result is that 
$$\dim H^1(\Yap,\Yap\sm sing(\Yap))=16.$$
Thus~(6.8.8) gives $b^1(\Sta\cap\Vta)=24$. 

\demo{Proof of~(6.8.9)}
Let $m_2\cl J\to J$  be ``multiplication by two''. Let
$\wt{C}:=m_2^{-1}(\T)$, and
$$\matrix
\wt{C} & \overset f\to\lra & C\\
x & \mapsto & u^{-1}(2x),
\endmatrix$$
where $u$ is the Abel-Jacobi map of~(1.6). Thus $f$ is the Galois
cover of $C$ with $J[2]$ as group of deck transformations. 
The map 
$$\matrix
\wt{C}\tm\wt{C} & \overset q_{\a}\to\lra & \Xi_{\a}\\
(x,y) & \mapsto & (f(x),{\hat x}-{\hat y}+{\hat\a})
\endmatrix
\tag{6.8.10}$$
is the quotient map for the diagonal action of $J[2]$ on
$\wt{C}\tm\wt{C}$. Thus $H^1(\Xi_{\a})$ is isomorphic to the subspace of
$H^1(\wt{C}\tm\wt{C})$ invariant for this action. If $\pi_i\cl
\wt{C}\tm\wt{C}\to\wt{C}$ is the projection to the $i$-th factor,
$$H^1(\wt{C}\tm\wt{C})^{J[2]}=
\pi_1^{*}f^{*}H^1(C)\op\pi_2^{*}f^{*}H^1(C).$$
This proves~(6.8.9).
\qed
\enddemo

The next result will be used to prove that $\Mt$ is simply-connected. Let $R$
be as in~(5.2.4): by Item~(2) of~(4.3.2)  we have an inclusion $i\cl
R\hra\Sta$.

\proclaim{(6.8.11)Lemma}
Keeping notation as above, we have
$$\Im(i_{\#})\ss\Im(\pi_1(\Vta\cap\Sta)\to\pi_1(\Sta)).$$
\endproclaim

\demo{Proof}
Referring to~(6.7.53) let $h_{\a}$ be given by
$$\matrix
C & \brel h_{\a}\over\lra & \Xi_{\a}\\
p & \mapsto & (p,{\hat\a}).
\endmatrix$$
By~(6.7.53) the map $h_{\a}$ lifts to a map ${\bar h}_{\a}\cl C\hra
\ov{Y}_{\a}$. Composing with the quotient map  $\ov{Y}_{\a}\to
Y^{+}_{\a}$ (see~(6.7.59)) we get an inclusion $h^{+}_{\a}\cl C\hra
Y^{+}_{\a}$. Let $\ell^{+}_{\a}\cl\Yap\hra\Sta$ be the inclusion given
by~(6.7.58). As is easily checked the image of  
$$\ell^{+}_{\a,\#}\circ
h^{+}_{\a,\#}\cl\pi_1(C)\to\pi_1(\Sta)\tag{6.8.12}$$ 
is contained in the image of $i_{\#}$. Next let $L_1,L_2\ss\ov{Y}_{\a}$ be
two exceptional divisors of $\pi_{\a}\circ\rho_{\a}$ whose images in
$\O_{\a}$ are indicized by the same $k$ (see~(6.7.53)): thus the
equivalence relation $\sim$ of~(6.7.59) glues together $L_1$ and $L_2$.
Let ${\bar\g}\cl[1,2]\to\ov{Y}_{\a}$ be a continuous path such that ${\bar
\g}(i)\in L_i$ and $\g(1)\sim\g(2)$. Composing ${\bar\g}$ with the
quotient map $\ov{Y}_{\a}\to Y^{+}_{\a}$ we get a loop
$\g^{+}\cl[1,2]\to Y^{+}_{\a}$ such that 
$\ell^{+}_{\a,\#}(\g^{+})$ is in the image of $i_{\#}$.  
Furthermore one easily checks that 
$i_{\#}\pi_1(R)$ is generated by
the image of~(6.8.12) and by $\ell^{+}_{\a,\#}(\g^{+})$.
Since both $\ell^{+}_{\a,\#}\circ h^{+}_{\a,\#}(C)$ and
$\ell^{+}_{\a,\#}(\g^{+})$ are  contained in $\Vta\cap\Sta$, this implies
the lemma.  
\qed
\enddemo

 Next we give results on $\Vta\cap\Bta$.
The composition
$$\Pap\overset\b_{+}^{-1}\to{\cdots>}\Pab\overset\b_{-}\to\lra\Pam
\overset f_{-}\to\lra\Ia$$
is a rational map $f_{+}\cl\Pap\cdots>\Ia$, regular away from $\G^{+}$.
The following result follows immediately from~(6.5.4)
and~(6.7.17). 

\proclaim{$\text{(6.8.13)}^{*}$Proposition}
Keeping notation as above, the intersection of $\Vta\cap\Bta$ and $\D$
inside $\Vta$ is transverse, and  consists of 
$16$ points, one on each $\D_k$. The rational map 
$$f_{+}\circ(\epa)^{-1}|_{\Vta\cap\Bta}\cl \Vta\cap\Bta\cdots>\Ia$$
is obtained by blowing up $\Bta\cap\D$, and then
contracting $16$ disjoint $(-1)$ curves (not intersecting the exceptional
divisor of the blow-up of $\Bta\cap\D$) to $\nua^{-1}\Jh[2]$. 
\endproclaim

\proclaim{$\text{(6.8.14)}^{*}$Corollary}
The map 
$$\pi_1(\Vta\cap\Bta)\to\pi_1(\Vta)\tag{6.8.15}$$
 induced by inclusion is an isomorphism. In particular
$$b^1(\Vta\cap\Bta)=4.\tag"$\text{(6.8.16)}^{*}$"$$
\endproclaim

\demo{Proof of the corollary}
We have a series of isomorphisms
$$\pi_1(\Vta)\cong\pi_1(\Pap)\cong\pi_1(\Pap\sm\G^{+})
\overset f_{+,\#}\to\lra\pi_1(\Ia\sm E)\cong\pi_1(\Ia).\tag{6.8.17}$$
If we compose with the map of~(6.8.15) we get an
isomorphism $\pi_1(\Vta\cap\Bta)\overset\sim\to{\to}\pi_1(\Ia)$,
by~(6.8.13). Thus the map of~(6.8.15) must be an isomorphism.
Formula~(6.8.16) follows from the first statement of the corollary
and~(6.8.2).  
\qed
\enddemo

Another result on $\Vta\cap\Bta$ that will be useful is the following:
$$\text{the map $\pi_1(\Vta\cap\Bta)\to\pi_1(\Bta)$ induced by
inclusion is trivial.}\tag"$\text{(6.8.18)}^{*}$"$$
In fact by~(6.5.4) we have 
$$\Vta\cap\Bta\ss g_1^{-1}(\text{node of $\ov{C}$}),$$
where $g_1$ is the map of~(5.1.2). Since the right-hand side of the above
formula is a $\PP^1$-bundle over the simply-connected
surface $K^{[2]}\Jh$ we get~(6.8.18).

We give a description of the triple intersection
$\Vta\cap\Bta\cap\Sta$. Let $D_i\ss\Jh$ (for $i=1,2$) be the smooth
irreducible curve defined by 
$$D_i:=\{{\hat y}\in\Jh|\ 2y\sim q_{3-i}-r,\quad\text{some $r\in C$}\}.$$
The last statement of~(6.5.14) gives an inclusion $D_i\hra\Pa$. Composing
with the birational map $\Pa\cdots>\Vta$ we get an inclusion $h_i\cl
D_i\hra\Vta$. It follows from~(6.5.14) that 
$$h_1(D_1)\cap h_2(D_2)=\Bta\cap\D.\tag{6.8.19}$$
The following proposition also follows easily from~(6.5.14).

\proclaim{$\text{(6.8.20)}^{*}$Proposition}
Keeping notation as above 
$$\Vta\cap\Bta\cap\Sta=h_1(D_1)\cup h_2(D_2).$$
In particular by~(6.8.19) the curve $\Vta\cap\Bta\cap\Sta$ is
connected.   
\endproclaim

Finally we give some results on restriction maps. 
Let $\Pi(\Sta)\ss H^2(\Sta)$ be given by
$$\Pi(\Sta):=\QQ c_1(\om_f)
\op H^2(C)\op H^2(\Jh)\op\QQ c_1((\ia\tm\id_{\Jh})^{*}\cL),
\tag"$\text{(6.8.21)}^{*}$"$$
where the right-hand side is viewed as a subspace of $H^2(\Sta)$
thanks to~(4.1.3) and Item~(3) of~(4.2.5). Let $\z_{\a}$ 
be as in~(6.7.51): we claim that 
$$\text{$\Pi(\Sta)\to H^2(\z_{\a}^{-1}(\Sta\cap\Vta))$ is injective.}
\tag"$\text{(6.8.22)}^{*}$"$$   
(The map is the restriction to
$\Sta\cap\Vta$ followed by pull-back for $\z_{\a}$.) By~(6.7.58)
and~(6.7.59) we have $\z_{\a}^{-1}(\Sta\cap\Vta)\cong\Yab$, and the
restriction to $\Yab$ of the map $f$ of~(4.1.1) is identified with the
blow-down map $\Yab\to\Xi_{\a}$, where $\Xi_{\a}$ is as in~(6.7.53). Thus
$\om_f$ has degree $(-2)$ on each of the exceptional divisors of   
$\Yab\to\Xi_{\a}$, and all the elements of the other direct summands
of~(6.8.21) have degree zero on the exceptional divisors.   Hence the kernel
of the map of~(6.8.22) is contained in  
$$H^2(C)\op H^2(\Jh)\op\QQ c_1((\ia\tm\id_{\Jh})^{*}\cL).
\tag{6.8.23}$$
Let $q_{\a}\cl\wt{C}\tm\wt{C}\to\Xi_{\a}$ be as in~(6.8.10). In order to
prove~(6.8.22) it suffices to show that the natural map from~(6.8.23) to
$H^2(\wt{C}\tm\wt{C})$ is injective. This is an easy exercise given the
explicit formula for $q_{\a}$ of (6.8.10).  

Let
$$H^1(\Sta\cap\Vta)\overset r_1\to{\lra} H^1(\Sta\cap\Bta\cap\Vta)
\qquad
(H^1(\Bta\cap\Vta)\overset r_2\to{\lra} H^1(\Sta\cap\Bta\cap\Vta)$$
be the maps induced by restriction. We claim that
$$\Im(r_1)+\Im(r_2)\ge 19.\tag"$\text{(6.8.24)}^{*}$"$$
In fact consider the exact sequence   
$$\multline
0\to\QQ\lra H^1(\Sta\cap\Bta\cap\Vta, sm(\Sta\cap\Bta\cap\Vta))\\
\lra H^1(\Sta\cap\Bta\cap\Vta)
\brel\psi\over\lra H^1(sm(\Sta\cap\Bta\cap\Vta))
\endmultline\tag{6.8.25}$$
The kernel of the map $\phi$ of~(6.8.8) maps surjectively to the kernel of
the map $\psi$ of~(6.8.25), hence $\ker(\psi)\ss\Im(r_1)$. On the other
hand  $\psi(\Im(r_2))=4$. Since $\dim\ker(\psi)=15$ we get~(6.8.24)
\head
7. Proof of Theorem~(1.4)
\endhead
By~(2.3.2) and~(2.3.6) we know that $\Mt$ is of pure dimension
$6$ and that $\omt$ is symplectic.  
\subhead
7.1. Proof that $\Mt$ is connected and simply connected
\endsubhead
By~(3.2) it suffices to show that
$$\gather
\text{$\Zta$ is connected,}\tag{7.1.1}\\
\pi_1(\Zta)=\{1\},\tag{7.1.2}
\endgather$$
where 
$$\Zta=\Sta\cup\Bta\cup\Vta.\tag{7.1.3}$$
Let us prove~(7.1.1). By~(4.4.1), (5.2.1), (6.8.1) each of $\Sta$,
$\Bta$, $\Vta$ is irreducible, and by~(4.3.2), (6.8.6)
and~(6.8.13)  every pairwise intersection is non-empty. By~(7.1.3) we get
that $\Zta$ is path-connected. Now we prove~(7.1.2). By~(6.8.20) the
triple intersection $\Sta\cap\Bta\cap\Vta$ is connected, hence
$\pi_1(\Zta)$ is generated by the images of $\pi_1(\Sta)$, $\pi_1(\Bta)$
and $\pi_1(\Vta)$. By~(4.4.6) and~(6.8.14) the images of $\pi_1(\Sta)$
and  $\pi_1(\Vta)$ are contained in the image of $\pi_1(\Bta)$, hence
$\pi_1(\Bta)\to\pi_1(\Zta)$ is surjectiv. By~(5.2.5) we get that
$\pi_1(R)\to\pi_1(\Zta)$ is surjective, where $R$ is as in~(5.2.4).
Thus~(6.8.11) gives  that $\pi_1(\Vta)\to\pi_1(\Zta)$ is surjective, and
hence by~(6.8.14)  $\pi_1(\Vta\cap\Bta)$ generates $\pi_1(\Zta)$.
By~(6.8.18) we conclude that $\Zta$ is
simply-connected.
\subhead
7.2. Proof that $b_2(\cMt)\le 8$
\endsubhead
The key is the following result.

\proclaim{(7.2.1)Proposition}
The map $H^2(\Mt:\ZZ)\to H^2(\Vta;\ZZ)$ induced by restriction is injective.
\endproclaim

The proposition implies the claimed
bound on $b_2(\cMt)$ because by~(6.8.3) the image of the restriction map 
is contained in $\L(\Vta)$, and this group has rank $8$ by Formula~(6.8.4).

The first element in the proof of~(7.2.1) consists of a
Mayer-Vietoris argument. Let $X_1:=\Sta$, $X_2:=\Bta$, $X_3:=\Vta$, and
$$C^p(H^q):=\bigoplus\limits_{1\le i_0<\cdots< i_p\le 3} 
H^q(X_{i_0}\cap\cdots\cap X_{i_p};\QQ).$$
We let $\d\cl C^p(H^q)\to C^{p+1}(H^q)$ be the usual \v Cech cochain 
map, and $Z^p(H^q)\ss C^p(H^q)$ be the group of $\d$-cocycles. 

\proclaim{(7.2.2)Proposition}
The map
$$H^2(\Zta)\to Z^0(H^2)$$
induced by restriction is injective.
\endproclaim

\demo{Proof}
This is proved by considering the Mayer-Vietoris spectral sequence
associated to the decomposition $\Zta=\Sta\cup\Bta\cup\Vta$, and
applying our results on the cohomology maps induced by restriction. More
precisely one can triangulate $\Zta$ so that every intersection of the $X_i$'s
is the support of a sub-triangulation. Let  
$$E_0^{p,q}:=\bigoplus\limits_{1\le i_0<\cdots< i_p\le 3} 
S^q(X_{i_0}\cap\cdots\cap X_{i_p}),$$
where $S^q(X_{i_0}\cap\cdots\cap X_{i_p})$ is the group of cochains
supported on $X_{i_0}\cap\cdots\cap X_{i_p}$. Letting ${\tilde\d}\cl
E_0^{p,q}\to E_0^{p+1,q}$ be the \v Check cochain map, and $d\cl
E_0^{p,q}\to E_0^{p,q+1}$ be the map induced by the coboundary of
simplicial cochains, we get a double complex which computes the
cohomology of $\Zta$, because $\tilde\d$ is exact except at $E_0^{0,q}$,
 where its homology is the group of simplicial
cochains of $\Zta$ (see~\cite{W, p.~202}). On the other hand the
filtration ``by $p$'' generates a spectral sequence 
 with $E_1$-term given by
$$E_1^{p,q}=C^p(H^q),$$
whose differential $E_1^{p,q}\to E_1^{p,q+1}$ is equal to the \v Check
cochain map $\d$ considered above. Thus $E_2^{0,2}=Z^0(H^2)$, and
 in order to prove the lemma it suffices to show that
$$E_2^{2,0}=E_2^{1,1}=0.$$
The first vanishing is trivial: by~(6.8.20) the intersection
$\Sta\cap\Bta\cap\Vta$ is connected, hence $\d\cl E_1^{1,0}\to
E_1^{2,0}$ is surjective. To prove the second vanishing we must show that
the complex
$$C^0(H^1)\brel\d^{0,1}\over\lra C^1(H^1)
\brel\d^{1,1}\over\lra C^2(H^1)\tag{7.2.3}$$
 is exact. By~(4.4.2), (5.2.2) and~(6.8.2) we have
$$\dim C^0(H^1)=14,\tag{7.2.4}$$
and Formulae~(4.4.3), (6.8.7), (6.8.16) give
$$\dim C^1(H^1)=33.\tag{7.2.5}$$
We claim that
$$\text{the map $\d^{0,1}$ is injective.}\tag{7.2.6}$$
This may be verified
directly, or else one may argue that since $E_2^{2,0}=0$ we have 
$H^1(\Zta)=E_2^{0,1}$; by~(7.1.2) we   get
$E_2^{0,1}=0$, i.e.~$\ker(\d^{0,1})=0$. 
By~(6.8.24) we have
$$\rk(\d^{1,1}) \ge 19.\tag{7.2.7}$$
Exactness of~(7.2.3) follows at once from~(7.2.4), (7.2.5), (7.2.6)
and~(7.2.7). 
\qed
\enddemo

We represent elements of $Z^0(H^2)$ as $\g=(\g_1,\g_2,\g_3)$,
where $\g_i\in H^2(X_i;\QQ)$. The second element in the proof of~(7.2.1)
consists of the following result.

\proclaim{(7.2.8)Proposition}
Keep notation as above. The map
$$\matrix
Z^0(H^2) & \lra & H^2(\Vta)\\
(\g_1,\g_2,\g_3) & \mapsto & \g_3
\endmatrix$$
is injective.
\endproclaim

\demo{Proof}
Since  $\g$ is a $\d$-cocycle
$$\g_1|_{\Sta\cap\Bta}=\g_2|_{\Sta\cap\Bta}.\tag{7.2.9}$$
Thus~(4.4.6) and~(5.2.7) give that 
$${\tilde\k}^{*}(\g_1)\in\QQ[C\tm E]
\op(\id_C\tm\nu_0)^{*}(H^2(C)\op H^2(\Jh)
\op\QQ c_1((\ia\tm\id_{\Jh})^{*}\cL)).$$
Applying~(4.4.9) we get that
$\g_1\in\Pi(\Sta)$, where $\Pi(\Sta)$ is as in~(6.8.21).
Now assume $\g_3=0$. 
Since $\g$ is a cocycle we get that the restriction of $\g_1$ to
$\Sta\cap\Vta$ vanishes, hence by
injectivity of~(6.8.22) we get $\g_1=0$. By~(7.2.9) and~(5.2.6) we
also get $\g_2=0$.  
\qed
\enddemo

\demo{Proof of Proposition~(7.2.1)}
Since $\cMt$ is simply-connected $H^2(\cMt;\ZZ)$ has no torsion, and thus 
it suffices to prove that the restriction map
$H^2(\cMt;\QQ)\to H^2(\Vta;\QQ)$ is injective. This map is equal to the
composition of three maps:   
$$H^2(\cMt)\to H^2(\Zta)\to Z^0(H^2)\to H^2(\Vta).$$
The first map is injective by~(3.2), the second map is injective
by~(7.2.2) and the third map is injective by~(7.2.8). Hence the composition
is injective. 
\enddemo

\subhead
7.3. Generators of $H^2(\cMt;\QQ)$
\endsubhead
There exists a compactification of the moduli space of slope-stable
vector-bundles on a projective surface which is ``smaller'' than the
moduli space of Gieseker-Maruyama semistable sheaves, namely Uhlenbeck's
compactification~\cite{L2,Mor}. Let $\Mv^U$ be the Uhlenbeck
compactification of the moduli space of slope-stable vector-bundles $F$ on
$J$ with $v(F)=\vv$: it is a projective variety and there
exists a regular map $\vf\cl \Mv\to\Mv^U$ which is an isomorphism on the
subset parametrizing slope-stable vector-bundles. Let $\cM^U:=\vf(\cM)$.
We have Donaldson's homomorphism~\cite{L2,Mor} 
$$\mu\cl H^2(J;\ZZ)\to H^2(\cM^U;\ZZ),\tag{7.3.1}$$
characterized by the following property: if $\cF$ is a family of semistable
sheaves on $J$ parametrized by $S$, with moduli belonging to $\cM$,
and $f\cl S\to\cM$ is the the modular map,
then 
$$f^{*}\circ\vf^{*}\circ\mu(\a)=p_{S,*}(c_2(\cF)\cup\a).\tag{7.3.2}$$
(Here $p_S\cl J\tm S\to S$ is the projection.)

\proclaim{(7.3.3)Proposition}
The homomorphism $\pit^{*}\circ\vf^{*}\circ\mu\cl H^2(J;\ZZ)\to
H^2(\Mt;\ZZ)$ is injective, and 
$$\pit^{*}\circ\vf^{*}\circ\mu(H^2(J;\ZZ))\ot\QQ,\quad\QQ c_1(\Sit),
\quad \QQ c_1(\Bt)$$
are linearly independent subspaces of $H^2(\Mt;\QQ)$. 
\endproclaim

\demo{Proof}
Let $x\in (J\sm J[2])$. The sheaf on $\Jh$ given by
$$Ext^1_{{\hat\phi}}(\phi^{*}(I_x)\ot\cL,
\phi^{*}(I_{-x})\ot\cL^{-1})\tag{7.3.4}$$
is locally-free of rank two because $x\not=-x$. (Here $\phi,\phih$ are as
in~(1.9).) Let $T$ be the projectivization of~(7.3.4): it parametrizes
semistable simple sheaves on $J$ whose moduli belong to $\cM$. In fact,
letting $h\cl T\to\Jh$ be the $\PP^1$-fibration, $\pi_J,\pi_T$ the
projections of $J\tm T$ to $J$ and $T$ respectively, and $\xi$ the
tautological sub-line-bundle on $T$, we have a tautological exact sequence of
sheaves on $J\tm T$ 
$$0\to \pi_J^{*}(I_{-x})\ot(\pi_J\tm h)^{*}\cL^{-1}\to\cE
\to\pi_J^{*}(I_{x})\ot(\pi_J\tm h)^{*}\cL\ot\pi_T^{*}\xi\to 0.$$
Since $\cE$ is a family of simple semistable sheaves on $J$ parametrized by
$T$, with $v(\cE_t)=\vv$ and $\det(\cE_t)\cong\cO_J$, $\sum
c_2(\cE_t)=0$, it induces by~(2.3.7) a regular map $\g\cl T\to\Mt$.
Applying~(7.3.2) one easily gets that  
$$\g^{*}\circ\pit^{*}\circ\vf^{*}\circ\mu(\a)=
h^{*}\circ(\d^{-1})^{*}(\a),$$
where $\d\cl J\to\Jh$ is the isomorphism~(1.7). This shows that 
$\pit^{*}\circ\vf^{*}\circ\mu$ is injective. Now we prove the second
statement of the proposition. 
Let $\G\ss\Sta$ and $\L\ss\Bta$ be generic fibers of the fibrations $f$
of~(4.1.1) and $g$ of~(5.1.1), respectively. Let us prove that
the intersection numbers of $c_1(\Sit)$ and $c_1(\Bt)$  with $\G$ and $\L$
are given by the entries of the following intersection matrix
$$\matrix\format\c\quad&\r\quad&\r\\
      & c_1(\Sit) & c_1(\Bt)\\
\G  &  -2        &     1       \\
\L   &  2            &    -2
\endmatrix\tag{7.3.5}$$
The restriction of $\pit$ to
$(\Sit\sm\Ot)$ is a $\PP^1$-fibration over $(\Si\sm\O)$, of which $\G$ is
a fiber: this gives the top left entry because by adjunction 
$K_{\Sit}\cong\cO_{\Sit}(\Sit)$. To get the other diagonal entry
one proceeds similarly: a dense open subset $\wt{U}\ss\Bt$ containing
$\L$ is the total space of a $\PP^1$-fibration $G\cl\wt{U}\to U$   defined
similarly to the map $g$ of~(5.1.1), and $\L$
is a fiber of $G$. Since by adjunction $K_{\Bt}\cong\cO_{\Bt}(\Bt)$, we get
the bottom right entry.   To get the off-diagonal entries first notice that
by~(4.3.10) $\Sit$ and $\Bt$ intersect transversely outside
$\pit^{-1}(\{[I_x\ot L\op I_x\ot L^{-1}\})$, hence 
$$\align
\langle c_1(\Sit),\L\rangle= & \#\{\l\in\L|\ \text{$\cG_{\l}$ is strictly
semistable}\},\\
\langle c_1(\Bt),\G\rangle= & \#\{\g\in\G|\ \text{$\cF_{\g}$ has two
singular points}\}.\\
\endalign$$
Here $\cG$ and $\cF$ are families of simple semistable sheaves on $J$
parametrized by $\L$ and $\G$ respectively which induce the inclusion maps
of $\L$ and $\G$ in $\Mt$ respectively. (We have used~(4.3.2) to get the
second equation.) Computing the right-hand side of the above equalities we
get the off-diagonal entries of~(7.3.5). Now assume that
$$x c_1(\Sit)+y c_1(\Bt)+u=0,$$
where $x,y\in\QQ$
and $u\in\pit^{*}\circ\vf^{*}\circ\mu(H^2(J;\ZZ))\ot\QQ$. We intersect with
$\G$ and $\L$, and notice that  
$$\langle u,\G\rangle=\langle u,\L\rangle=0.$$
In fact the first intersection number vanishes because $\pit$ contracts $\G$,
and the second intersection number vanishes because $\pit(\L)$ is
contracted by $\vf$ (see~\cite{L2}). Since the intersection matrix
of~(7.3.5) is non-singular we get that $x=y=0$, and hence also $u=0$. 
\qed
\enddemo

Now we can finish the proof of Theorem~(1.4). By the previous subsection
we have $b_2(\Mt)\le 8$, and by the above proposition $b_2(\Mt)\ge 8$,
hence $b_2(\Mt)=8$. Finally, since $c_1(\Sit)$ and $c_1(\Bt)$ are of type
$(1,1)$ and the map $\pit^{*}\circ\vf^{*}\circ\mu$ is a morphism of Hodge
structures (this follows from~(7.3.2)), we get that
$$h^{2,0}(\Mt)=h^{2,0}(J)=1.$$ 


\Refs
\widestnumber\key{Muk2}
%
\ref \key B
\by A. Beauville  
\paper Vari\'et\'es K\"ahl\'eriennes dont la premi\`ere classe de Chern est
nulle 
\jour J. Diff. Geom
\vol 18 \yr 1983 \pages 755-782
\endref
%
\ref \key BDL 
\by J. Bryan - R. Donagi - Naichung Conan Leung
\paper $G$-Bundles on abelian surfaces, hyperk\"ahler manifolds, and stringy
Hodge numbers 
\jour AG/0004159
\yr 2000
\endref
%
\ref \key HG 
\by D. Huybrechts - L. G\"ottsche 
\paper Hodge numbers of moduli spaces of stable bundles on $K3$ surfaces 
\jour Internat. J. of Math. 
\vol 7\yr 1996\pages 359-372 
\endref
%
\ref \key  G
\by D. Gieseker 
\paper On the moduli of vector bundles on an algebraic surface
\jour Ann. of Math
\vol 106 \yr 1977 \pages 45-60
\endref
%
\ref \key GM 
\by M. Goresky - R. MacPherson
\book Stratified Morse Theory
\bookinfo Ergeb. Math. Grenzgeb. (3. Folge) 14
\publ Springer
\yr 1988
\endref
%
\ref \key L1
\by J. Li 
\paper The first two Betti numbers of the moduli space of vector bundles on
surfaces  
\jour Comm. Anal. Geom.
\vol 5 \yr 1997 \pages 625-684
\endref
%
\ref \key L2
\by J. Li 
\paper Algebraic geometric interpretation of Donaldson's polynomial
invariants  
\jour J. Differential Geometry
\vol 37 \yr 1993 \pages 417-466
\endref
%
\ref \key Mor
\by J. Morgan 
\paper Comparison of the Donaldson polynomial
invariants with tehir algebro geometric analogues  
\jour Topology
\vol 32 \yr 1993 \pages 449-488
\endref
%
\ref \key Muk2
\by S. Mukai 
\paper Symplectic structure of the moduli space of sheaves on an abelian or
$K3$ surface 
\jour Invent. math
\vol 77 \yr 1984 \pages 101-116
\endref
%
\ref \key Muk3
\by S. Mukai 
\paper On the moduli space of bundles on
$K3$ surfaces 
\inbook Vector bundles on Algebraic Varieties
\bookinfo Tata Institute Studies in Mathematics
\publ Oxford University Press
\finalinfo 1987
\endref
%
\ref \key Mum 
\by D. Mumford 
\book Abelian Varieties
\bookinfo Tata Institute Studies in Mathematics
\publ Oxford University Press
\finalinfo 1970
\endref
%
\ref \key O1 
\by K. O'Grady
\paper Desingularized moduli spaces of sheaves on a $K3$
\jour J. reine angew. Math.
\vol 512 \yr 1999 \pages 49-117
\endref
%
\ref \key O2
\by K. O'Grady
\paper The weight-two Hodge  structure  of moduli
spaces of sheaves on a $K3$ surface
\jour J. of Algebraic Geom.
\vol 6\yr 1997\pages 599-644
\endref
%
\ref \key O3
\by K. O'Grady
\paper Moduli of vector-bundles on surfaces
\inbook Algebraic Geometry Santa Cruz 1995
\bookinfo Proc. Symp. Pure Math. vol. 62
\publ Amer. Math. Soc.
\yr 1997
\pages 101-126
\endref
%
\ref \key O4
\by K. O'Grady
\paper Relations among Donaldson polynomials of
certain algebraic surfaces, I
\jour Forum Math.
\vol 8\yr 1996\pages 1-61
\endref
%
\ref \key Y1
\by K. Yoshioka
\paper Some examples of Mukai's reflections on $K3$ surfaces
\jour J. Reine Angew. Math.
\vol 515\yr 1999\pages 97-123
\endref
%
\ref \key Y2
\by K. Yoshioka
\paper Moduli spaces of stable sheaves on abelian surfaces
\jour AG/0009001
\yr 2000
\endref
%
\ref \key W
\by F. W. Warner 
\book Foundations of Differentiable Manifolds and Lie Groups
\bookinfo Graduate Texts in Mathematics 94
\publ Springer-Verlag
\finalinfo 1983
\endref

%
\endRefs
\enddocument